\documentclass[a4paper]{amsart}
\usepackage{amscd, amssymb, booktabs, chemarrow}
\usepackage{myown}

\title[Lusztig Parametrizations for Classical Groups]
{On the Unicity and the Ambiguity of Lusztig Parametrizations for Finite Classical Groups}
\author{Shu-Yen Pan}
\address{Department of Mathematics,
National Tsing Hua University, Hsinchu 300, Taiwan}
\email{sypan@math.nthu.edu.tw}

\keywords{rank, unipotent character, reductive dual pair}
\subjclass[2010]{Primary: 20C33; Secondary: 22E50}

\date{\today}

\begin{document}

\begin{abstract}
The Lusztig correspondence is a bijective mapping between the Lusztig series indexed by the conjugacy class
of a semisimple element $s$ in the connected component $(G^*)^0$ of the dual group of $G$
and the set of irreducible unipotent characters of the centralizer of $s$ in $G^*$.
In this article we discuss the unicity and ambiguity of such a bijective correspondence.
In particular, we show that the Lusztig correspondence for a classical group can be made to be unique
if we require it to be compatible with the parabolic induction and the finite theta correspondence.
\end{abstract}

\maketitle
\tableofcontents

\section{Introduction}

\subsection{}\label{0101}
Let $\bfG$ be a classical group defined over a finite field $\bfF_q$ of odd characteristic,
and let $F$ be the corresponding Frobenius endomorphism.
Let $G=\bfG^F$ denote the finite subgroup of rational points,
and let $\cale(\bfG)$ denote the set of irreducible characters
(i.e., the characters of irreducible representations) of $G$.

Let $R^\bfG_{\bfT^*,s}$ denote the \emph{Deligne-Lusztig virtual characters}
(\cf.~\cite{carter-finite}, \cite{GM-guide}) indexed by conjugacy class of pair
$(\bfT^*,s)$ where $\bfT^*$ is a rational maximal torus in the dual group $\bfG^*$ and
$s$ is a rational semisimple element contained in $\bfT^*$.
Let $\calv(\bfG)$ be the space of (complex valued) class function on $G$,
and let $\calv(\bfG)^\sharp$ denote the subspace spanned by Deligne-Lusztig virtual characters.
Note that $\calv(\bfG)$ is an inner product space with the inner product $\langle,\rangle_\bfG$ given by
\[
\langle f_1,f_2\rangle_\bfG=\frac{1}{|G|}\sum_{g\in G}f_1(g)\overline{f_2(g)}
\]
for $f_1,f_2\in\calv(\bfG)$, and $\overline{f_2(g)}$ denotes the complex conjugate of the value $f_2(g)$.
For $f\in\calv(\bfG)$ the orthogonal projection of $f$ onto $\calv(\bfG)^\sharp$
is denoted by $f^\sharp$ and called the \emph{uniform projection}.
A class function $f\in\calv(\bfG)$ is called \emph{uniform} if $f^\sharp=f$.

A natural question is how much an irreducible character $\rho$ of a classical group
$G$ can be determined by its uniform projection $\rho^\sharp$?

If $\bfG$ is a general linear group $\GL_n$ or a unitary group $\rmU_n$,
then every irreducible character is uniform (i.e., $\rho=\rho^\sharp$) and the above question is trivial.
However, if $\bfG$ is a symplectic group or an orthogonal group,
the question is more subtle.
Recall that $\rho\in\cale(\bfG)$ is called \emph{unipotent} if $\langle\rho,R^\bfG_{\bfT^*,1}\rangle_\bfG\neq 0$
for some $\bfT^*$.
If $\rho$ is unipotent and $\bfG$ is connected,
it is known that $\rho$ is uniquely determined by its uniform projection,
i.e., $\rho'^\sharp=\rho^\sharp$ if and only if $\rho'=\rho$
(\cf.~\cite{DM-Lusztig} proposition 6.3 and \cite{GM-guide} theorem 4.4.23).
If $\rho$ is unipotent and $\bfG=\rmO^\epsilon_{2n}$,
it is also known that $\rho'^\sharp=\rho^\sharp$ if and only if $\rho'=\rho$ or $\rho'=\rho\cdot\sgn$
(\cf.~\cite{pan-Lusztig-correspondence} proposition~3.6).

In this paper, the above question will be answered completely
(\cf.~Corollary~\ref{0607}, Corollary~\ref{0608} and Corollary~\ref{0702}) for classical groups:
\begin{thm}
Let $\bfG$ be a symplectic group or an orthogonal group, and let $\rho,\rho'\in\cale(\bfG)$.
Then $\rho'^\sharp=\rho^\sharp$ if and only of
\begin{equation}\label{0103}
\rho'=\begin{cases}
\rho, & \text{if\/ $\bfG=\SO_{2n+1}$};\\
\rho,\rho^c & \text{if\/ $\bfG=\Sp_{2n}$};\\
\rho,\rho^c,\rho\cdot\sgn,\rho^c\cdot\sgn & \text{if\/ $\bfG=\rmO^\epsilon_{2n}$}.
\end{cases}
\end{equation}
\end{thm}
Here ``$\sgn$'' denotes the sign character of an orthogonal group,
and ``$\rho^c$'' denotes the character obtained from $\rho$ by conjugating an element in
the corresponding similitude group (\cf.~\cite{waldspurger} \S 4.3, \S 4.10).
Note that if $\rho$ is unipotent then $\rho^c=\rho$, and then (\ref{0103}) is reduced to
the above known result.

\subsection{}
In fact, the above theorem is a consequence of a more precise result on the ambiguity
of the Lusztig parametrization of irreducible characters of finite classical groups.
From now on, we assume that $\bfG$ is a symplectic group or an orthogonal group.
It is known by Lusztig that there is a partition
\[
\cale(\bfG)=\bigcup_{(s)}\cale(\bfG,s)
\]
where the union $\bigcup_{(s)}$ runs over conjugacy classes of semisimple elements in the connected component
of the dual group $G^*$ and
\[
\cale(\bfG,s)=\{\,\rho\in\cale(\bfG)\mid\langle\rho,R^\bfG_{\bfT^*,s}\rangle_\bfG\neq0\text{ for some }\bfT^*
\text{ containing }s\,\}.
\]
Each $\cale(\bfG,s)$ is called a \emph{Lusztig series}.
In particular, the subset $\cale(\bfG,1)$ consists of unipotent characters.

Now we first focus on the set of unipotent characters.
Let $\calv(\bfG,1)$ denote the subspace spanned by elements in $\cale(\bfG,1)$,
and let $\calv(\bfG,1)^\sharp$ denote the uniform projection of $\calv(\bfG,1)$.
Following from Lusztig (\cf.~\cite{lg-symplectic}, \cite{lg-orthogonal}) we can define
a set $\cals_\bfG^\sharp$ (\cf.~(\ref{0335})) and class functions $R_\Sigma^\bfG$ (\cf.~(\ref{0330}))
for $\Sigma\in\cals_\bfG^\sharp$ such that the set $\{\,R_\Sigma^\bfG\mid\Sigma\in\cals_\bfG^\sharp\,\}$
forms an orthonormal basis for $\calv(\bfG,1)^\sharp$.
Lusztig constructs a set $\cals_\bfG$ (\cf.~(\ref{0210})) of similar classes of \emph{symbols} and
a bijective mapping $\call_1^\bfG\colon\cals_\bfG\rightarrow\cale(\bfG,1)$ denoted by
$\Lambda\mapsto\rho_\Lambda$ such that the value $\langle\rho_\Lambda,R^\bfG_\Sigma\rangle_\bfG$
is specified (\cf.~Proposition~\ref{0301}).
Such a mapping $\call_1^\bfG$ is called a \emph{Lusztig parametrization} of unipotent characters of $\bfG$.
Because the uniform projection $\rho_\Lambda^\sharp$ can be obtained by the values
$\langle\rho_\Lambda,R^\bfG_\Sigma\rangle_\bfG$,
the problem of the uniqueness of $\call_1^\bfG$ is equivalent to the problem whether a unipotent character
$\rho_\Lambda$ can be uniquely determined by its uniform projection $\rho_\Lambda^\sharp$.

As described in Subsection~\ref{0101}, $\call_1^\bfG$ is known to be unique if $\bfG=\SO_{2n+1}$ or $\Sp_{2n}$.
However, due to the disconnectedness of $\rmO^\epsilon_{2n}$,
the mapping $\call_1^{\rmO^\epsilon_{2n}}$ is not uniquely determined.
By using the result of ``\emph{cells}'' by Lusztig, we determine how ambiguous a Lusztig
parametrization $\call_1^{\rmO^\epsilon_{2n}}$ could be and we show in Proposition~\ref{0305}
that $\call_1^{\rmO^\epsilon_{2n}}$ can be chosen to be unique if we require
it to be
\begin{itemize}
\item \emph{compatible with parabolic induction} and

\item \emph{compatible with theta correspondence} on cuspidal characters or ${\bf 1}_{\rmO^+_2},\sgn_{\rmO^+_2}$.
\end{itemize}
In particular, $\call_1^{\rmO^\epsilon_{2n}}$ (and $\call_1^{\Sp_{2n'}}$) can be given so that
$(\rho_\Lambda,\rho'_{\Lambda'})\in\Theta_{\bfG,\bfG'}$ if and only if
$(\Lambda,\Lambda')\in\calb_{\bfG,\bfG'}$ for $(\bfG,\bfG')=(\rmO^\epsilon_{2n},\Sp_{2n'})$
where $\calb_{\bfG,\bfG'}$ is a relation between $\cals_\bfG$ and $\cals_{\bfG'}$
defined in Subsection~\ref{0403}.

\subsection{}
For a general Lusztig series $\cale(\bfG,s)$,
Lusztig shows (\cf.~\cite{lusztig-book}) that there exists a bijection
\[
\grL_s\colon\cale(\bfG,s)\rightarrow\cale(C_{\bfG^*}(s),1)
\]
satisfying
\begin{equation}\label{0201}
\langle\rho,\epsilon_\bfG R^{\bfG}_{\bfT^*,s}\rangle_\bfG
=\langle\grL_s(\rho),\epsilon_{C_{\bfG^*}(s)}R^{C_{\bfG^*}(s)}_{\bfT^*,1}\rangle_{C_{\bfG^*}(s)}
\end{equation}
where $\epsilon_\bfG=(-1)^{\kappa(\bfG)}$, $\kappa(\bfG)$ denotes the rational rank of\/ $\bfG$,
and $C_{\bfG^*}(s)$ denotes the centralizer of $s$ in $\bfG^*$.
Such a bijection $\grL_s$ will be called a \emph{Lusztig correspondence} in this article
(it is called a \emph{Jordan decomposition} in \cite{GM-guide}).

Now the question is to understand whether the Lusztig correspondence $\grL_s$ is uniquely determined by (\ref{0201}).
If the answer is negative, then we want to know what kind of conditions need to be enforced to make $\grL_s$ unique.
Some discussion on this problem can be founded in \cite{GM-guide} appendix A.5.
If $\bfG$ is a connected group with connected center, it is shown in \cite{DM-Lusztig} theorem 7.1 that $\grL_s$
can be uniquely determined by (\ref{0201}) and some extra conditions.

For $s\in G^*$, we have a decomposition $s=s^{(0)}\times s^{(1)}\times s^{(2)}$
where $s^{(1)}$ (resp.~$s^{(2)}$) is the part whose eigenvalues are all equal to $-1$ (resp.~$1$),
and $s^{(0)}$ is the part whose eigenvalues do not contain $1$ or $-1$.
Then we can define groups $\bfG^{(0)}(s)$, $\bfG^{(-)}(s)$ and $\bfG^{(+)}(s)$ (\cf.~(\ref{0517}))
so that there is a bijective mapping
\[
\cale(C_{\bfG^*}(s),1)\rightarrow\cale(\bfG^{(0)}(s)\times\bfG^{(-)}(s)\times\bfG^{(+)}(s),1).
\]
Combining Lusztig parametrization $\call_1$ of unipotent characters for $\bfG^{(0)}$, $\bfG^{(-)}(s)$, $\bfG^{(+)}(s)$,
and above bijection and the inverse $\grL_s^{-1}$ of a Lusztig correspondence,
we obtain a bijective mapping
\[
\call_s\colon\cals_{\bfG^{(0)}(s)}\times\cals_{\bfG^{(-)}(s)}\times\cals_{\bfG^{(+)}(s)}\rightarrow\cale(\bfG,s)
\]
denoted by $(x,\Lambda_1,\Lambda_2)\mapsto\rho_{x,\Lambda_1,\Lambda_2}$
which is called a \emph{modified Lusztig correspondence}.
Then we prove the following results on the unicity and ambiguity of $\call_s$
(or equivalently, the unicity and ambiguity of $\grL_s$) for classical groups:

\begin{enumerate}
\item Suppose that $\bfG=\SO_{2n+1}$.
There is a unique modified Lusztig correspondence $\call_s$
(\cf.~Theorem~\ref{0507}).
Note that $\SO_{2n+1}$ is a connected group with connected center,
so this case is covered by \cite{DM-Lusztig} theorem 7.1.
However, for this group, no extra condition other than (\ref{0201}) is needed to ensure
the uniqueness of $\call_s$.

\item Suppose that $\bfG=\rmO^\epsilon_{2n}$ where $\epsilon=+$ or $-$.
There exists a unique modified Lusztig correspondence $\call_s$
which is compatible with the parabolic induction and some other conditions (on \emph{basic characters}).
(\cf.~Theorem~\ref{0508}).

\item Suppose that $\bfG=\Sp_{2n}$.
Here we provide two choices of the modified Lusztig correspondence $\call_s$:
\begin{enumerate}
\item There exists a unique modified Lusztig correspondence $\call_s$
which is compatible with the parabolic induction and compatible with the theta correspondence for
the dual pair $(\bfG,\bfG')=(\Sp_{2n},\SO_{2n'+1})$, i.e., we show that
$(\rho_{x,\Lambda_1,\Lambda_2},\rho_{x',\Lambda_1',\Lambda_2'})\in\Theta_{\bfG,\bfG'}^\psi$
if and only if
\begin{itemize}
\item $s^{(0)}=-s'^{(0)}$ and $x=x'$,

\item $\Lambda_2=\Lambda'_1$, and

\item $(\Lambda_1,\Lambda'_2)\in\calb_{\bfG^{(-)}(s),\bfG^{(+)}(s')}$
\end{itemize}
where $\rho_{x',\Lambda_1',\Lambda_2'}$ is given by the unique modified Lusztig correspondence in (1)
(\cf.~Theorem~\ref{0506} and Theorem~\ref{0610}).

\item There exists a unique modified Lusztig correspondence $\call_s$ which is compatible with
the parabolic induction and compatible with the theta correspondence for the dual pair $(\bfG,\bfG')=(\Sp_{2n},\rmO^\epsilon_{2n'})$,
i.e., $(\rho_{x,\Lambda_1,\Lambda_2},\rho_{x',\Lambda_1',\Lambda_2'})\in \Theta_{\bfG,\bfG'}^\psi$
if and only if
\begin{itemize}
\item $s^{(0)}=s'^{(0)}$ and $x=x'$,

\item $\Lambda_1=\Lambda'_1$, and

\item $(\Lambda_2,\Lambda'_2)\in\calb_{\bfG^{(+)}(s),\bfG^{(+)}(s')}$
\end{itemize}
where $\rho_{x',\Lambda_1',\Lambda_2'}$ is given by the unique modified Lusztig correspondence
in (2) (\cf.~Theorem~\ref{0903} and Theorem~\ref{0801}).
\end{enumerate}
It seems that two choices in (a) and (b) of the modified Lusztig correspondence
$\call_s$ for $\Sp_{2n}$ should be the same.
However, we do not know how to obtain the conclusion yet.
\end{enumerate}

\subsection{} The contents of this article are as follows.
In Section 2, we recall the notion and some basic result of ``\emph{symbols}''
and ``\emph{special symbols}'' by Lusztig from \cite{lg}.
In Section 3 we first recall the notion of ``\emph{almost characters}'' by Lusztig
from \cite{lg-symplectic} and \cite{lg-orthogonal}.
Then we record some results of cells from \cite{pan-finite-unipotent}.
In Section 4 we show the uniqueness of $\call_1^\bfG$ for $\bfG=\Sp_{2n}$ and $\SO_{2n+1}$
by using the results in the previous section.
Moreover, we also discuss the ambiguity of $\call_1^\bfG$ for $\bfG=\rmO^\epsilon_{2n}$.
In Section 5 we discuss the relation between the theta correspondence
$\Theta^\psi_{\bfG,\bfG'}$ on unipotent characters for $(\bfG,\bfG')=(\Sp_{2n},\rmO^\epsilon_{2n'})$ and
Lusztig parametrizations $\call_1^\bfG\colon\cals_\bfG\rightarrow\cale(\bfG,1)$
for $\bfG=\Sp_{2n},\rmO^\epsilon_{2n}$.
In Section 6 we discuss the relation between the theta correspondence
$\Theta^\psi_{\bfG,\bfG'}$ on certain Lusztig series for $(\bfG,\bfG')=(\Sp_{2n},\SO_{2n'+1})$
or $(\Sp_{2n},\rmO^\epsilon_{2n'})$ and the Lusztig correspondence
$\grL_s\colon\cale(\bfG,s)\rightarrow\cale(C_{\bfG^*}(s),1)$.
In Section 7 we show that the Lusztig correspondence
$\grL_s\colon\cale(\bfG,s)\rightarrow\cale(C_{\bfG^*}(s),1)$ is unique for $\bfG=\SO_{2n+1}$.
In Section 8 we show that $\grL_s$ can be chosen to be unique for $\bfG=\rmO^\epsilon_{2n}$
if we require $\grL_s$ to be compatible with the parabolic induction and some other conditions
on ``basic characters''.
In the final section we discuss the uniqueness of the Lusztig correspondence $\grL_s$ for
$\bfG=\Sp_{2n}$.
In particular, we show that a unique $\grL_s$ can be chosen to be compatible with
the theta correspondence for the dual pair $(\Sp_{2n},\SO_{2n+1})$ 
or for the dual pair $(\Sp_{2n},\rmO^\epsilon_{2n})$.

\section{Symbols and Special symbols}

\subsection{Irreducible characters of Weyl groups}\label{0209}
Let $\calp(n)$ denote the set partitions of $n$.
It is known that the set of irreducible characters $\cale(S_n)$ of the symmetric group $S_n$
is parametrized by $\calp(n)$.
For $\lambda\in\calp(n)$, we write $|\lambda|=n$,
and the corresponding irreducible character of $S_n$ is denoted by $\varphi_\lambda$.

Let $W_n$ denote the Coxeter group of type $B_n$, i.e.,
$W_n$ consists of all permutations on $\{1,2,\ldots,n,n^*,(n-1)^*,\ldots,1^*\}$ which
commutes with the involution
\[
(1,1^*)(2,2^*)\cdots(n,n^*)
\]
where $(i,j)$ denote the transposition of $i,j$.
For $i=1,\ldots,n-1$, let
\[
s_i=(i,i+1)(i^*,(i+1)^*)\quad \text{and}\quad \sigma_n=(n,n^*).
\]
It is known that $W_n$ is generated by $\{s_1,\ldots,s_{n-1},\sigma_n\}$.
Each element of $W_n$ induces a permutation of $\{1,2,\ldots,n\}$.
So we have a surjective homomorphism $W_n\rightarrow S_n$ with kernel isomorphic to $(\bbZ/2\bbZ)^n$.
Therefore, $\varphi_\lambda$ can be regarded as an irreducible character of $W_n$.

An ordered pair $\sqbinom{\mu}{\nu}$ of two partitions is called a \emph{bi-partitions}.
Let $\calp_2(n)$ denote the set of bipartitions of $n$, i.e.,
\[
\textstyle
\calp_2(n)=\left\{\,\sqbinom{\mu}{\nu}\mid |\mu|+|\nu|=n\,\right\}.
\]
For a bi-partition $\sqbinom{\mu}{\nu}$,
we define its \emph{transpose} by $\sqbinom{\mu}{\nu}^\rmt=\sqbinom{\nu}{\mu}$.
A bi-partition $\sqbinom{\mu}{\nu}$ is called \emph{degenerate} if $\mu=\nu$,
and it is called \emph{non-degenerate} otherwise.
For $\sqbinom{\mu}{\nu}\in\calp_2(n)$ such that $\mu\in\calp(k)$ and $\nu\in\calp(l)$
where $k+l=n$, we define
\begin{equation}\label{0208}
\varphi_{\sqbinom{\mu}{\nu}}
=\Ind^{W_n}_{W_k\times W_l}(\varphi_\mu\otimes(\varepsilon_l\varphi_\nu))
\end{equation}
where $\varepsilon_l\colon W_l\rightarrow\{\pm1\}$ is given by $s_i\mapsto 1$ and $\sigma_l\mapsto -1$.
It is known that $\varphi_{\sqbinom{\mu}{\nu}}$ is an irreducible character of $W_n$,
and the mapping $\calp_2(n)\rightarrow\cale(W_n)$ by $\sqbinom{\mu}{\nu}\mapsto\varphi_{\sqbinom{\mu}{\nu}}$
gives a parametrization of $\cale(W_n)$ such that
\begin{itemize}
\item $\varphi_{\sqbinom{n}{-}}={\bf 1}_{W_n}$

\item $\varphi_{\sqbinom{\mu}{\nu}}\cdot\varepsilon_n=\varphi_{\sqbinom{\nu}{\mu}}$,
in particular $\varphi_{\sqbinom{-}{n}}=\varepsilon_n$
\end{itemize}
(\cf.~\cite{Geck-Pfeiffer} theorem 5.5.6).

The kernel $W_n^+$ of $\varepsilon_n$ is a subgroup of index two generated by
$\{s_1,\ldots,s_{n-1},\sigma_n s_{n-1}\sigma_n\}$.
Let $W_n^-=W_n\smallsetminus W_n^+$.
It is well known that if $\sqbinom{\mu}{\nu}$ is non-degenerate,
then $\varphi_{\sqbinom{\mu}{\nu}}|_{W_n^+}=\varphi_{\sqbinom{\nu}{\mu}}|_{W_n^+}$
which is an irreducible character of $W_n^+$;
if $\sqbinom{\mu}{\nu}$ is degenerate,
then $\varphi_{\sqbinom{\mu}{\nu}}|_{W_n^+}$ is a sum of two irreducible characters of $W_n^+$.

\subsection{Lusztig's symbols}
In this subsection, we recall some basic notations of ``symbols'' from \cite{lg}.
A \emph{symbol} $\Lambda$ is an ordered pair
\begin{equation}\label{0212}
\Lambda=\binom{A}{B}=\binom{a_1,a_2,\ldots,a_{m_1}}{b_1,b_2,\ldots,b_{m_2}}
\end{equation}
of two finite sets $A,B$ (possibly empty) of non-negative integers.
The sets $A,B$ are also denoted by $\Lambda^*,\Lambda_*$ and called the
\emph{first row}, the \emph{second row} of $\Lambda$ respectively.
A symbol $\Lambda$ is called \emph{reduced} if $0\not\in A\cap B$.
If $\Lambda=\binom{A}{B}$, then we define its \emph{transpose} by $\Lambda^\rmt=\binom{B}{A}$.
We denote $\Lambda_1\subset\Lambda_2$ and call $\Lambda_1$ a \emph{subsymbol} of $\Lambda_2$
if both $\Lambda_1^*\subset\Lambda_2^*$ and $(\Lambda_1)_*\subset(\Lambda_2)_*$.
If $\Lambda_1\subset\Lambda_2$, their difference is defined by
$\Lambda_2\smallsetminus\Lambda_1=\binom{\Lambda_2^*\smallsetminus\Lambda_1^*}{(\Lambda_2)_*\smallsetminus(\Lambda_1)_*}$.
If both $\Lambda_1^*\cap\Lambda_2^*=\emptyset$ and $(\Lambda_1)_*\cap(\Lambda_2)_*=\emptyset$,
we define $\Lambda_1\cup\Lambda_2=\binom{\Lambda_1^*\cup\Lambda_2^*}{(\Lambda_1)_*\cup(\Lambda_2)_*}$.

For a symbol $\Lambda$ given in (\ref{0212}), its \emph{rank} and \emph{defect} are defined by
\begin{align}
\begin{split}
{\rm rk}(\Lambda) &=\sum_{i=1}^{m_i}a_i+\sum_{j=1}^{m_2}b_j
-\left\lfloor\left(\frac{|A|+|B|-1}{2}\right)^2\right\rfloor, \\
{\rm def}(\Lambda) &=|A|-|B|.
\end{split}
\end{align}
From the definition, it is not difficult to check that
\begin{equation}\label{0204}
{\rm rk}(\Lambda)\geq\left\lfloor\left(\frac{{\rm def}(\Lambda)}{2}\right)^2\right\rfloor.
\end{equation}
A symbol $\Lambda$ is called \emph{degenerate} if $\Lambda^\rmt=\Lambda$.
If $\Lambda$ is degenerate, then ${\rm rk}(\Lambda)$ is even and ${\rm def}(\Lambda)=0$.

We define an equivalence relation ``$\sim$'' on symbols generated by
\[
\binom{a_1,,a_2,\ldots,a_{m_1}}{b_1,b_2,\ldots,b_{m_2}}\sim
\binom{a_1+1,,a_2+1,\ldots,a_{m_1}+1,0}{b_1+1,b_2+1,\ldots,b_{m_2}+1,0}.
\]
If $\Lambda_1\sim\Lambda_2$, two symbols $\Lambda_1,\Lambda_2$ are called \emph{similar}.
It is not difficult to see that two symbols in the same similar class have the same rank and the same defect,
and each similar class contains a unique reduced symbol.
Let $\cals$ denote the set of similar classes of symbols,
and let $\cals_{n,\delta}\subset\cals$ denote the subset of similar classes of symbols of rank $n$
and defect $\delta$.

A symbol $\Lambda$ is called \emph{cuspidal} if (\ref{0204}) is an equality.
It is not difficult to see that a symbol is cuspidal if and only if it is similar to a symbol of the forms
$\binom{k,k-1,\ldots,0}{-}$ or $\binom{-}{k,k-1,\ldots,0}$ for some non-negative integer $k$.
Note that $\binom{A}{-}$ means that the second row of the symbol is the empty set.

A mapping $\Upsilon$ from symbols to bi-partitions is given by
\begin{equation}
\Upsilon\colon\binom{a_1,a_2,\ldots,a_{m_1}}{b_1,b_2,\ldots,b_{m_2}}\mapsto
\sqbinom{a_1-(m_1-1),a_2-(m_1-2),\ldots,a_{m_1-1}-1,a_{m_1}}{b_1-(m_2-1),b_2-(m_2-2),\ldots,b_{m_2}-1,b_{m_2}}.
\end{equation}
If $\Upsilon(\Lambda)=\sqbinom{\mu}{\nu}$, we write $\Upsilon(\Lambda)^*=\mu$ and
$\Upsilon(\Lambda)_*=\nu$ to denote the first row and the second row of the bi-partition $\Upsilon(\Lambda)$.
We can check that $\Upsilon(\Lambda_1)=\Upsilon(\Lambda_2)$ if and only if $\Lambda_1\sim\Lambda_2$,
and then $\Upsilon$ gives a bijection
\begin{equation}\label{0207}
\cals_{n,\delta}\rightarrow\calp_2\left(n-\bigl\lfloor\bigl(\tfrac{\delta}{2}\bigr)^2\bigr\rfloor\right).
\end{equation}

Modified from Lusztig, we define
\begin{align}\label{0210}
\begin{split}
\cals_{\rmO^+_{2n}} &=\{\,\Lambda\in\cals\mid {\rm rk}(\Lambda)=n,\ {\rm def}(\Lambda)\equiv 0\pmod 4\,\}, \\
\cals_{\Sp_{2n}} &=\{\,\Lambda\in\cals\mid {\rm rk}(\Lambda)=n,\ {\rm def}(\Lambda)\equiv 1\pmod 4\,\}, \\
\cals_{\rmO^-_{2n}} &=\{\,\Lambda\in\cals\mid {\rm rk}(\Lambda)=n,\ {\rm def}(\Lambda)\equiv 2\pmod 4\,\}, \\
\cals_{\SO_{2n+1}} &=\{\,\Lambda\in\cals\mid {\rm rk}(\Lambda)=n,\ {\rm def}(\Lambda)\equiv 3\pmod 4\,\}.
\end{split}
\end{align}
Note that $\Lambda\in\cals_{\rmO^\epsilon_{2n}}$ if and only if $\Lambda^\rmt\in\cals_{\rmO^\epsilon_{2n}}$
where $\epsilon=+$ or $-$.
Then we define
\begin{equation}\label{0206}
\cals_{\SO^\epsilon_{2n}}
=\{\,\Lambda\in\cals_{\rmO^\epsilon_{2n}}\mid\Lambda\neq\Lambda^\rmt\,\}/\{\Lambda,\Lambda^\rmt\}
\cup\{\,\Lambda^{\rm I},\Lambda^{\rm II}\mid\Lambda\in\cals_{\rmO^\epsilon_{2n}},\ \Lambda=\Lambda^\rmt\,\},
\end{equation}
i.e., in $\cals_{\SO^\epsilon_{2n}}$ a non-degenerate symbol is identified with its transpose,
and a degenerated symbol $\Lambda$ occurs with multiplicity $2$ and the two copies are denoted by
$\Lambda^{\rm I},\Lambda^{\rm II}$ respectively.
Note that $\cals_{\rmO^-_{2n}}$ does not contain any degenerate symbols and so
$\cals_{\SO^-_{2n}}=\cals_{\rmO^-_{2n}}/\{\Lambda,\Lambda^\rmt\}$.

\subsection{Special symbols}\label{0213}
Let $\bfG=\Sp_{2n}$, $\SO_{2n+1}$, $\rmO^\epsilon_{2n}$ or $\SO^\epsilon_{2n}$ where $\epsilon=+$ or $-$.
A symbol $Z=\binom{a_1,a_2,\ldots,a_{m+1}}{b_1,b_2,\ldots,b_m}$ of defect $1$ is called \emph{special}
if $a_1\geq b_1\geq a_2\geq b_2\geq\cdots\geq a_m\geq b_m\geq a_{m+1}$;
similarly, a symbol $Z=\binom{a_1,a_2,\ldots,a_m}{b_1,b_2,\ldots,b_m}$ of defect $0$ is called \emph{special}
if $a_1\geq b_1\geq a_2\geq b_2\geq\cdots\geq a_m\geq b_m$.
Define
\[
\delta_0=\begin{cases}
1, & \text{if $\bfG=\Sp_{2n}$ or $\SO_{2n+1}$}; \\
0, & \text{if $\bfG=\SO^\epsilon_{2n}$ or $\rmO^\epsilon_{2n}$}.
\end{cases}
\]
For a special symbol $Z$ of defect $\delta_0$, we define
\begin{align*}
\cals_Z &=\{\,\Lambda\in\cals\mid\Lambda^*\cup\Lambda_*=Z^*\cup Z_*,\ \Lambda^*\cap\Lambda_*=Z^*\cap Z_*\,\}, \\
\cals_{Z,\delta_0} &=\{\,\Lambda\in\cals_Z\mid{\rm def}(\Lambda)=\delta_0\,\}, \\
\cals_Z^\bfG &=\cals_Z\cap \cals_\bfG,
\end{align*}
i.e., $\cals_Z$ is the subset of $\cals$ consisting of the symbols of the exactly the same entries of $Z$.
Note that in the above definition of $\cals_Z^\bfG$, the special symbol $Z$ is not required to be in $\cals_\bfG$.
It is clear that
\begin{equation}\label{0211}
\cals_\bfG=\bigcup_{Z}\cals_Z^\bfG
\end{equation}
where $\bigcup_Z$ runs over
\begin{itemize}
\item special symbols of rank $n$ and defect $1$ if $\bfG=\Sp_{2n}$ or $\SO_{2n+1}$;

\item special symbols of rank $n$ and defect $0$ if $\bfG=\SO^+_{2n}$ or $\rmO^+_{2n}$;

\item non-degenerate special symbols of rank $n$ and defect $0$ if $\bfG=\SO^-_{2n}$ or $\rmO^-_{2n}$.
\end{itemize}

\begin{exam}
\begin{enumerate}
\item Suppose that $Z=\binom{1}{0}\in\cals_{1,0}$.
Then we have
\[
\textstyle
\cals_Z^{\rmO^+_2}=\left\{\binom{1}{0},\binom{0}{1}\right\},\qquad
\cals_Z^{\rmO^-_2}=\left\{\binom{-}{1,0},\binom{1,0}{-}\right\}.
\]

\item Suppose that $Z=\binom{2,0}{1}\in\cals_{2,1}$.
Then we have
\[
\textstyle
\cals_Z^{\Sp_4}=\left\{\binom{2,0}{1},\binom{2,1}{0},\binom{1,0}{2},\binom{-}{2,1,0}\right\},\qquad
\cals_Z^{\SO_5}=\left\{\binom{1}{2,0},\binom{0}{2,1},\binom{2}{1,0},\binom{2,1,0}{-}\right\}.
\]
\end{enumerate}
\end{exam}

For a special symbol $Z$, let $Z_\rmI$ denote the subsymbol consisting of ``\emph{singles}'', i.e.,
\[
Z_\rmI=Z\smallsetminus\binom{Z^*\cap Z_*}{Z^*\cap Z_*},
\]
and we define the \emph{degree} of $Z$ by
\begin{equation}
\deg(Z)=\frac{|Z_\rmI|-{\rm def}(Z)}{2}
\end{equation}
where $|Z_\rmI|$ denotes the number of entries in $Z_\rmI$, i.e.,
$|Z_\rmI|=|(Z_\rmI)^*|+|(Z_\rmI)_*|$.
Note that $\deg(Z)$ is always a non-negative integer.
For a subsymbol $M\subset Z_\rmI$, we define a symbol $\Lambda_M\in\cals_Z$ by
\begin{equation}
\Lambda_M=(Z\smallsetminus M)\cup M^\rmt.
\end{equation}
It is not difficult to see that
\[
\cals_Z^\bfG=\begin{cases}
\{\,\Lambda_M\mid M\subset Z_\rmI,\ |M|\ \text{even}\,\}, & \text{if $\bfG=\Sp_{2n}$ and ${\rm def}(Z)=1$};\\
\{\,\Lambda_M\mid M\subset Z_\rmI,\ |M|\ \text{odd}\,\}, & \text{if $\bfG=\SO_{2n+1}$ and ${\rm def}(Z)=1$};\\
\{\,\Lambda_M\mid M\subset Z_\rmI,\ |M|\ \text{even}\,\}, & \text{if $\bfG=\rmO^+_{2n}$ and ${\rm def}(Z)=0$};\\
\{\,\Lambda_M\mid M\subset Z_\rmI,\ |M|\ \text{odd}\,\}, & \text{if $\bfG=\rmO^-_{2n}$ and ${\rm def}(Z)=0$}.
\end{cases}
\]
Note that
\[
\cals_Z^{\SO^\epsilon_{2n}}=\begin{cases}
\{Z^{\rm I},Z^{\rm II}\}, & \text{if $Z$ is degenerate};\\
\cals_Z^{\rmO^\epsilon_{2n}}/\{\Lambda,\Lambda^\rmt\}, & \text{if $Z$ is non-degenerate}
\end{cases}
\]
where $Z^{\rm I},Z^{\rm II}$ are both equal to $Z$ but are regarded as two elements.
Then we have
\[
|\cals_Z^\bfG|=\begin{cases}
2^{2\deg(Z)}, & \text{if $\bfG=\Sp_{2n}$ or $\SO_{2n+1}$};\\
2^{2\deg(Z)-1}, & \text{if $\bfG=\rmO^\epsilon_{2n}$ and $\deg(Z)>0$};\\
2^{2\deg(Z)-2}, & \text{if $\bfG=\SO^\epsilon_{2n}$ and $\deg(Z)>0$};\\
1, & \text{if $\bfG=\rmO^+_{2n}$ and $\deg(Z)=0$};\\
2, & \text{if $\bfG=\SO^+_{2n}$ and $\deg(Z)=0$}.
\end{cases}
\]
Note that the family $\cals_Z^\bfG$ is empty if $\bfG=\rmO^-_{2n},\SO^-_{2n}$ and $\deg(Z)=0$.

Finally we define a pairing $\langle,\rangle\colon\cals_Z^\bfG\times\cals_{Z,\delta_0}\rightarrow\bfF_2$ by
\begin{equation}\label{0205}
\langle\Lambda_{M_1},\Lambda_{M_2}\rangle=|M_1\cap M_2|\pmod 2.
\end{equation}

\begin{lem}\label{0203}
Let $Z$ be a special symbol of defect $\delta_0$.
Then for any $\Lambda\in\cals_Z^\bfG$ and $\Sigma\in\cals_{Z,\delta_0}$,
we have
\begin{align*}
\langle\Lambda,\Sigma\rangle &=\langle\Lambda^\rmt,\Sigma\rangle, \\
\langle\Lambda,\Sigma\rangle &\begin{cases}
=\langle\Lambda,\Sigma^\rmt\rangle, & \text{if\/ $\bfG=\SO^+_{2n},\rmO^+_{2n}$};\\
\neq\langle\Lambda,\Sigma^\rmt\rangle, & \text{if\/ $\bfG=\SO^-_{2n},\rmO^-_{2n}$}.
\end{cases}
\end{align*}
\end{lem}
\begin{proof}
Suppose that $\Lambda=\Lambda_{M_1}$ and $\Sigma=\Lambda_{M_2}$
for some $M_1,M_2\subset Z_\rmI$.
It is clear that $\Lambda^\rmt=\Lambda_{Z_\rmI\smallsetminus M_1}$.
The assumption $\Sigma\in\cals_{Z,\delta_0}$ implies that $|M_2|$ is even.
Then
\[
|M_1\cap M_2|+|(Z_\rmI\smallsetminus M_1)\cap M_2|=|M_2|\equiv 0\pmod 2.
\]
Hence the first equality is obtained.

Now suppose that $\bfG=\SO^\epsilon_{2n}$ or $\rmO^\epsilon_{2n}$.
Note that $|M_1|$ is even if $\epsilon=+$, and $|M_1|$ is odd if $\epsilon=-$.
Then the remaining assertion is obtained by the analogous argument.
\end{proof}

\section{Lusztig Parametrization of Unipotent Characters}

\subsection{Unipotent characters of $\GL_n$ or $\rmU_n$}\label{0326}
In this subsection, let $\bfG$ be a general linear group $\GL_n$ or a unitary group $\rmU_n$.
It is well known that the Weyl group $W_\bfG=S_n$,
and $\cale(S_n)$ is parametrized by $\calp(n)$.
For $\lambda\in\calp(n)$, we define
\[
R_\lambda^\bfG=\frac{1}{|S_n|}\sum_{w\in S_n}\varphi_\lambda(w)R^\bfG_{\bfT_w,1}
\]
where $\varphi_\lambda$ denotes the irreducible character of $S_n$ corresponding to $\lambda$.
It is known that $R_\lambda^\bfG$ is an irreducible unipotent character of $G$,
i.e., $R_\lambda^\bfG\in\cale(\bfG,1)$.
Let $\cals_{\GL_n}=\cals_{\rmU_n}=\calp(n)$.
Then the mapping
\[
\call_1\colon\cals_\bfG\rightarrow\cale(\bfG,1)\quad \text{ given by }
\lambda\mapsto\rho_\lambda:=R_\lambda^\bfG
\]
is a bijection.
To distinguish $\cals_{\rmU_n}$ from $\cals_{\GL_n}$,
an element in $\cals_{\rmU_n}$ will be denoted by $[\bar\lambda_1,\ldots,\bar\lambda_m]$
for $[\lambda_1,\ldots,\lambda_m]\in\calp(n)$.
An element of the form $[\bar k,\bar{k-1},\ldots,\bar 1]\in\cals_{\rmU_{k(k+1)/2}}$ for some $k$ is
called \emph{cuspidal}.
For a unitary group,
it is well known that $\lambda$ is cuspidal if and only if $\rho_\lambda$ is cuspidal.

It is known that the parametrization $\call_1$ above is compatible with the parabolic induction
on unipotent characters.
More precisely,
let $\bfG_n=\GL_n$, $\rmU_{2n}$ or $\rmU_{2n+1}$.
For $\rho\in\cale(\bfG_n,1)$, define
\begin{equation}\label{0323}
\Omega(\rho)=\left\{\,\rho'\in\cale(\bfG_{n+1},1)\mid
\left\langle\rho',R^{\bfG_{n+1}}_{\bfG_n\times\GL^\dag_1}(\rho\otimes{\bf1})\right\rangle_{\!\bfG_{n+1}}\neq 0\,\right\}
\end{equation}
where $R^{\bfG_{n+1}}_{\bfG_n\times\GL^\dag_1}$ is the standard parabolic induction,
$\GL^\dag_1=\GL_1$ defined over $\bfF_q$ if $\bfG_n=\GL_n$;
and $\GL^\dag_1$ is the restriction to $\bfF_q$ of $\GL_1$ defined over a quadratic extension of $\bfF_q$
if $\bfG_n=\rmU_{2n}$ or $\rmU_{2n+1}$.
For $\lambda\in\cals_{\bfG_n}$,
we define $\Omega(\lambda)$ to be a subset of $\cals_{\bfG_{n+1}}$ consisting of partitions of the following types:
\begin{itemize}
\item If $\bfG_n=\GL_n$ and $\lambda\in\cals_{\bfG_n}$,
then $\Omega(\lambda)$ consists of all partitions $\lambda'\in\cals_{\bfG_{n+1}}$
whose Young diagrams are obtained by adding a box to the Young diagram of $\lambda$.

\item If $\bfG_n=\rmU_{2n}$ or $\rmU_{2n+1}$ and $\lambda\in\cals_{\bfG_n}$,
then $\Omega(\lambda)$ consists of all partitions $\lambda'\in\cals_{\bfG_{n+1}}$
whose Young diagrams are obtained by adding two boxes in the same row or in the same column to the
Young diagram of $\lambda$.
\end{itemize}

\begin{exam}
\begin{enumerate}
\item Suppose that $\lambda=[3,1,1]\in\cals_{\GL_5}$.
Then
\[
\Omega(\lambda)=\{[4,1,1],[3,2,1],[3,1,1,1]\}\subset\cals_{\GL_6}.
\]

\item Suppose that $\lambda=[\bar 3,\bar 1,\bar 1]\in\cals_{\rmU_5}$.
Then
\[
\Omega(\lambda)=\{[\bar 5,\bar 1,\bar 1],[\bar 3,\bar 3,\bar 1],[\bar 3,\bar 2,\bar 2],
[\bar 3,\bar 1,\bar 1,\bar 1,\bar 1]\}\subset\cals_{\rmU_7}.
\]
\end{enumerate}
\end{exam}

Now the parametrization $\call_1\colon\cals_\bfG\rightarrow\cale(\bfG,1)$ by $\lambda\mapsto\rho_\lambda$
is said to be \emph{compatible with parabolic induction} if the following diagram
\begin{equation}\label{0324}
\begin{CD}
\cals_{\bfG_n} @> \Omega >> \cals_{\bfG_{n+1}} \\
@V \call_1 VV @VV \call_1 V \\
\cale(\bfG_n,1) @> \Omega >> \cale(\bfG_{n+1},1)
\end{CD}
\end{equation}
commutes, i.e., $\Omega(\rho_\lambda)=\{\,\rho_{\lambda'}\mid\lambda'\in\Omega(\lambda)\,\}$
for any $\lambda\in\cals_{\bfG_n}$.

\subsection{Uniform almost characters}
First suppose $\bfG$ is a connected classical group $\SO_{2n+1}$, $\Sp_{2n}$, or $\SO^\epsilon_{2n}$.
For a rational maximal torus $\bfT^*$ in the dual group $\bfG^*$ and a rational element $s\in T^*$,
the \emph{Deligne-Lusztig virtual characters} $R^\bfG_{\bfT^*,s}$ is defined in \cite{dl}
(see also \cite{carter-finite}).
For the disconnected group $\rmO^\epsilon_{2n}$,
$R^{\rmO^\epsilon_{2n}}_{\bfT^*,s}$ is defined by
\begin{equation}\label{0303}
R^{\rmO^\epsilon_{2n}}_{\bfT^*,s}=\Ind_{\SO^\epsilon_{2n}}^{\rmO^\epsilon_{2n}}R^{\SO^\epsilon_{2n}}_{\bfT^*,s}.
\end{equation}
For rational maximal tori in $\bfG$, the following is well known (\cf.~\cite{srinivasan} \S 3):
\begin{itemize}
\item If $\bfG=\Sp_{2n}$ or $\SO_{2n+1}$,
it is known that any rational maximal torus in $\bfG$ is conjugate (under $G$) to $\bfT_w$ for some $w\in W_n$.
Moreover, $\bfT_{w}$ and $\bfT_{w'}$ are conjugate if and only if $w,w'$ are conjugate under $W_n$.

\item If $\bfG=\SO^\epsilon_{2n}$ or $\rmO^\epsilon_{2n}$,
it is known that any rational maximal torus is conjugate (under $G$) to $\bfT_w$ for some $w\in W_n^\epsilon$
(\cf.~Subsection~\ref{0209}).
Moreover, $\bfT_{w}$ and $\bfT_{w'}$ are conjugate if and only if $w,w'$ are conjugate under $W_n$.
\end{itemize}

We recall some definitions from \cite{lg-CBMS} \S 3.17 and \S 3.19
(see also \cite{pan-Lusztig-correspondence} \S 2.3).
For $\Sigma\in\cals_{n,\delta_0}$ (\cf.~Subsection~\ref{0213}),
we define a uniform unipotent class function $R_\Sigma^\bfG\in\calv(\bfG,1)^\sharp$ by
\begin{align}\label{0330}
\begin{split}
R^\bfG_\Sigma &=
\begin{cases}
\displaystyle\frac{1}{\sqrt 2|W_n^+|}\sum_{w\in W_n^+}\varphi_{\Upsilon(\Sigma)}(w) R^\bfG_{\bfT_w,1},
& \text{if $\bfG=\SO^+_{2n}$ and $\Sigma$ degenerate};\\
\displaystyle\frac{1}{|W_n^\epsilon|}\sum_{w\in W_n^\epsilon}\varphi_{\Upsilon(\Sigma)}(w) R^\bfG_{\bfT_w,1},
& \text{if $\bfG=\SO^\epsilon_{2n}$ and $\Sigma$ non-degenerate};\\
\displaystyle\frac{1}{|W_n|}\sum_{w\in W_n}\varphi_{\Upsilon(\Sigma)}(w) R^\bfG_{\bfT_w,1},
& \text{if $\bfG=\Sp_{2n}$};\\
\displaystyle\frac{1}{|W_n|}\sum_{w\in W_n}\varphi_{\Upsilon(\Sigma^\rmt)}(w) R^\bfG_{\bfT_w,1},
& \text{if $\bfG=\SO_{2n+1}$}.
\end{cases} \\
R_\Sigma^{\rmO^\epsilon_{2n}} &=
\frac{1}{\sqrt 2}\Ind_{\SO^\epsilon_{2n}}^{\rmO^\epsilon_{2n}}(R^{\SO^\epsilon_{2n}}_\Sigma).
\end{split}
\end{align}
Note that $\varphi_{\Upsilon(\Sigma)}$ is the irreducible character of $W_n$ given in (\ref{0208})
and $\Upsilon$ is the bijection $\cals_{n,\delta_0}\rightarrow\calp_2(n)$ given in (\ref{0207}).
The class function $R^\bfG_\Sigma$ is called an \emph{almost character} of $G$.

\begin{lem}
Suppose that $\bfG=\rmO^\epsilon_{2n}$ where $\epsilon=+$ or $-$.
Then $R^\bfG_{\Sigma^\rmt}=\epsilon R^\bfG_\Sigma$ for any $\Sigma\in\cals_{n,0}$.
\end{lem}
\begin{proof}
From (\ref{0330}), we know that
\[
R^{\rmO^\epsilon_{2n}}_{\Sigma^\rmt}=
\begin{cases}
\displaystyle\frac{1}{2|W_n^+|}\sum_{w\in W_n^+}\varphi_{\Upsilon(\Sigma^\rmt)}(w) R^{\rmO^\epsilon_{2n}}_{\bfT_w,1},
& \text{if $\epsilon=+$ and $\Sigma$ degenerate};\\
\displaystyle\frac{1}{\sqrt 2|W_n^\epsilon|}\sum_{w\in W_n^\epsilon}\varphi_{\Upsilon(\Sigma^\rmt)}(w)
R^{\rmO^\epsilon_{2n}}_{\bfT_w,1},
& \text{if $\Sigma$ non-degenerate}.
\end{cases}
\]
Moreover, from Subsection~\ref{0209}, we have
$\varphi_{\Upsilon(\Sigma^\rmt)}(w)=\epsilon \varphi_{\Upsilon(\Sigma)}(w)$ for
$w\in W^\epsilon_n$.
Then the lemma is proved.
\end{proof}

Now we define a set $\cals_\bfG^\sharp$ as follows.
\begin{equation}\label{0335}
\cals_\bfG^\sharp=\begin{cases}
\cals_{n,1},
& \text{if $\bfG=\Sp_{2n},\SO_{2n+1}$};\\
\cals_{n,0}/\{\Sigma,\Sigma^\rmt\},
& \text{if $\bfG=\SO^+_{2n},\rmO^+_{2n}$};\\
\{\,\Sigma\in\cals_{n,0}\mid\Sigma\neq\Sigma^\rmt\,\}/\{\Sigma,\Sigma^\rmt\},
& \text{if $\bfG=\SO^-_{2n},\rmO^-_{2n}$}.
\end{cases}
\end{equation}
Here $\cals_{n,0}/\{\Sigma,\Sigma^\rmt\}$ means that we choose a representative
for each subset $\{\Sigma,\Sigma^\rmt\}$ in $\cals_{n,0}$.

\begin{exam}
We know that $\cals_{2,0}=\left\{\binom{2}{0},\binom{0}{2},\binom{2,1}{1,0},\binom{1,0}{2,1},\binom{1}{1}\right\}$.
The following are all four possible choices of $\cals_{\SO^-_4}^\sharp$ (or $\cals_{\rmO^-_4}^\sharp$):
\[
\textstyle
\left\{\binom{2}{0},\binom{2,1}{1,0}\right\},\quad
\left\{\binom{2}{0}\binom{1,0}{2,1}\right\},\quad
\left\{\binom{0}{2},\binom{2,1}{1,0}\right\},\quad
\left\{\binom{0}{2},\binom{1,0}{2,1}\right\}.
\]
Moreover, for each above choice $\cals_{\SO^-_4}^\sharp$,
we see that $\cals_{2,0}\smallsetminus\cals_{\SO^-_4}^\sharp$ is a choice for $\cals_{\SO^+_4}^\sharp$
(or $\cals_{\rmO^+_4}^\sharp$).
\end{exam}

Similar to (\ref{0211}), we have the decomposition
\[
\cals_\bfG^\sharp=\bigcup_Z{\cals_Z^\bfG}^\sharp
\]
where the union $\bigcup_Z$ is as in (\ref{0211}) and
\begin{equation}
{\cals_Z^\bfG}^\sharp=\begin{cases}
\cals_Z\cap\cals_{n,1}, & \text{if $\bfG=\Sp_{2n},\SO_{2n+1}$};\\
(\cals_Z\cap\cals_{n,0})/\{\Sigma,\Sigma^\rmt\}, & \text{if $\bfG=\SO^\epsilon_{2n},\rmO^\epsilon_{2n}$}.
\end{cases}
\end{equation}

\begin{lem}\label{0315}
Let $\bfG=\Sp_{2n}$, $\SO_{2n+1}$, $\SO^\epsilon_{2n}$ or $\rmO^\epsilon_{2n}$.
Then the set $\{\,R_\Sigma^\bfG\mid\Sigma\in\cals_\bfG^\sharp\,\}$
forms an orthonormal basis for $\calv(\bfG,1)^\sharp$.
\end{lem}
\begin{proof}
First we show that the set $\{\,R_\Sigma^\bfG\mid\Sigma\in\cals_\bfG^\sharp\,\}$ is orthonormal.
\begin{enumerate}
\item Suppose that $\bfG=\Sp_{2n}$ or $\SO_{2n+1}$.
From \cite{carter-finite} proposition 7.3.6, it is not difficult
to see that
\[
\langle R^\bfG_\Sigma, R^\bfG_{\Sigma'}\rangle_\bfG
=\langle\varphi_{\Upsilon(\Sigma)},\varphi_{\Upsilon(\Sigma')}\rangle_{W_n}
=\begin{cases}
1, & \text{if $\Sigma=\Sigma'$};\\
0, & \text{if $\Sigma\neq\Sigma'$}.
\end{cases}
\]

\item Suppose that $\bfG=\SO^+_{2n}$.
Then
\[
\langle R^\bfG_\Sigma, R^\bfG_{\Sigma'}\rangle_\bfG
=\begin{cases}
\langle\varphi_{\Upsilon(\Sigma)},\varphi_{\Upsilon(\Sigma')}\rangle_{W_n^+},
& \text{if both $\Sigma,\Sigma'$ are non-degenerate};\\
\frac{1}{\sqrt 2}\langle\varphi_{\Upsilon(\Sigma)},\varphi_{\Upsilon(\Sigma')}\rangle_{W_n^+},
& \text{if exactly one of $\Sigma,\Sigma'$ is degenerate};\\
\frac{1}{2}\langle\varphi_{\Upsilon(\Sigma)},\varphi_{\Upsilon(\Sigma')}\rangle_{W_n^+},
& \text{if both $\Sigma,\Sigma'$ are degenerate}.
\end{cases}
\]
We know that
\[
\langle\varphi_{\Upsilon(\Sigma)},\varphi_{\Upsilon(\Sigma)}\rangle_{W_n^+}=\begin{cases}
1, & \text{if $\Sigma$ non-degenerate};\\
2, & \text{if $\Sigma$ degenerate}.
\end{cases}
\]
Therefore, we conclude that
\[
\langle R^\bfG_\Sigma, R^\bfG_{\Sigma'}\rangle_\bfG
=\begin{cases}
1, & \text{if $\Sigma=\Sigma'$};\\
0, & \text{if $\Sigma\neq\Sigma',\Sigma'^\rmt$}.
\end{cases}
\]

\item Suppose that $\bfG=\SO^-_{2n}$.
Fix an element $\sigma_0\in W_n^-$, so an element in $W_n^-$ can be written as
$\sigma_0 w$ for some $w\in W_n^+$.
We know that $\varphi_{\Upsilon(\Sigma)}(\sigma_0w)$ is either equal to
$\varphi_{\Upsilon(\Sigma)}(w)$ or $-\varphi_{\Upsilon(\Sigma)}(w)$.
Therefore, by the proof in (2) we again conclude that
\[
\langle R^\bfG_\Sigma, R^\bfG_{\Sigma'}\rangle_\bfG
=\begin{cases}
1, & \text{if $\Sigma=\Sigma'$};\\
0, & \text{if $\Sigma\neq\Sigma',\Sigma'^\rmt$}.
\end{cases}
\]

\item Suppose that $\bfG=\rmO^\epsilon_{2n}$.
From (\ref{0330}), we know that
\[
R^{\rmO^\epsilon_{2n}}_{\Sigma}(g)=\begin{cases}
\sqrt 2R^{\SO^\epsilon_{2n}}_{\Sigma}(g), & \text{if $g\in\SO^\epsilon_{2n}(q)$};\\
0, & \text{if $g\in\rmO^\epsilon_{2n}(q)\smallsetminus\SO^\epsilon_{2n}(q)$}
\end{cases}
\]
for any $\Sigma\in\cals_\bfG^\sharp$.
Because $|\rmO^\epsilon_{2n}(q)|=2|\SO^\epsilon_{2n}(q)|$, we have
\[
\langle R^{\rmO^\epsilon_{2n}}_\Sigma, R^{\rmO^\epsilon_{2n}}_{\Sigma'}\rangle_{\rmO^\epsilon_{2n}}
=\langle R^{\SO^\epsilon_{2n}}_\Sigma, R^{\SO^\epsilon_{2n}}_{\Sigma'}\rangle_{\SO^\epsilon_{2n}}.
\]
\end{enumerate}

Next we show that the set $\{\,R_\Sigma^\bfG\mid\Sigma\in\cals_\bfG^\sharp\,\}$
spans $\calv(\bfG,1)^\sharp$.
It is well known that
\begin{equation}\label{0331}
\sum_{\varphi\in\cale(W_n)}\overline{\varphi(w)}\varphi(w')=\begin{cases}
\frac{|W_n|}{|\calo_w|}, & \text{if $w,w'$ are conjugate};\\
0, & \text{if $w,w'$ are not conjugate}
\end{cases}
\end{equation}
where $\calo_w$ denotes the conjugacy class of $w$ (\cf.~\cite{sr-finite} proposition 7).
\begin{enumerate}
\item Suppose that $\bfG=\Sp_{2n}$.
For $w\in W_n$, by (\ref{0330}) and (\ref{0331}) we have
\begin{align*}
\sum_{\Sigma\in\cals_\bfG^\sharp}\overline{\varphi_{\Upsilon(\Sigma)}(w)}R^\bfG_\Sigma
&=\frac{1}{|W_n|}\sum_{w'\in W_n}
\left[\sum_{\Sigma\in\cals_{n,1}}\overline{\varphi_{\Upsilon(\Sigma)}(w)}
\varphi_{\Upsilon(\Sigma)}(w')\right] R^\bfG_{\bfT_{w'},1} \\
&= \frac{1}{|\calo_w|}\sum_{w'\in W_n,\ w'\text{ conjugate to }w}R^\bfG_{\bfT_{w'},1} \\
&= R_{\bfT_w,1}^\bfG.
\end{align*}
The proof for $\bfG=\SO_{2n+1}$ is similar.

\item Suppose that $\bfG=\SO^-_{2n}$.
For $w\in W_n^-$, by (\ref{0330}) and (\ref{0331}) we have
\begin{align*}
\sum_{\Sigma\in\cals_\bfG^\sharp}\overline{\varphi_{\Upsilon(\Sigma)}(w)}R^\bfG_\Sigma
&=\frac{1}{2}\sum_{\Sigma\in\cals_{n,0}}\overline{\varphi_{\Upsilon(\Sigma)}(w)}R^\bfG_\Sigma \\
&=\frac{1}{2|W_n^-|}\sum_{w'\in W_n^-}
\left[\sum_{\Sigma\in\cals_{n,0}}\overline{\varphi_{\Upsilon(\Sigma)}(w)}
\varphi_{\Upsilon(\Sigma)}(w')\right] R^\bfG_{\bfT_{w'},1} \\
&= \frac{1}{|\calo_w|}\sum_{w'\in W_n^-,\ w'\text{ conjugate to }w}R^\bfG_{\bfT_{w'},1} \\
&= R_{\bfT_w,1}^\bfG.
\end{align*}
Note that for $w\in W_n^-$, the whole conjugacy class $\calo_w$ is contained in $W_n^-$.
The proof for $\bfG=\rmO^-_{2n}$ is similar.

\item Suppose that $\bfG=\SO^+_{2n}$.
For $w\in W_n^+$, we have
\begin{align*}
& \sum_{\Sigma\in\cals_\bfG^\sharp,\ \Sigma\text{ non-degenerate}}
\overline{\varphi_{\Upsilon(\Sigma)}(w)}R^\bfG_\Sigma
+\sum_{\Sigma\in\cals_\bfG^\sharp,\ \Sigma\text{ degenerate}}
\frac{1}{\sqrt 2}\overline{\varphi_{\Upsilon(\Sigma)}(w)}R^\bfG_\Sigma \\
&=\frac{1}{2|W_n^+|}\sum_{w'\in W_n^+}
\left[\sum_{\Sigma\in\cals_{n,0}}\overline{\varphi_{\Upsilon(\Sigma)}(w)}
\varphi_{\Upsilon(\Sigma)}(w')\right] R^\bfG_{\bfT_{w'},1} \\
&= R_{\bfT_w,1}^\bfG.
\end{align*}
The proof for $\bfG=\rmO^+_{2n}$ is similar.
\end{enumerate}

Therefore for all cases, we show that $\{\,R_\Sigma^\bfG\mid\Sigma\in\cals_\bfG^\sharp\,\}$
forms an orthonormal basis for $\calv(\bfG,1)^\sharp$.
\end{proof}

\subsection{Lusztig's parametrization of unipotent characters}\label{0325}

\begin{prop}[Lusztig]\label{0301}
Let $\bfG=\Sp_{2n}$, $\SO_{2n+1}$, $\SO^\epsilon_{2n}$, or $\rmO^\epsilon_{2n}$ where $\epsilon=+$ or $-$.
There exists a bijection $\call_1=\call_1^\bfG\colon\cals_\bfG\rightarrow\cale(\bfG,1)$
by $\Lambda\mapsto\rho_\Lambda$ satisfying

\begin{equation}\label{0332}
\langle\rho_\Lambda,R_\Sigma^\bfG\rangle_\bfG
=\begin{cases}
\frac{1}{c_Z}(-1)^{\langle\Lambda,\Sigma\rangle},
& \text{if $\Lambda\in\cals_Z^\bfG,\Sigma\in{\cals_Z^\bfG}^\sharp$ for some special $Z$};\\
0, & \text{otherwise},
\end{cases}
\end{equation}
where
\[
c_Z=\begin{cases}
2^{\deg(Z)}, & \text{if\/ $\bfG=\Sp_{2n},\SO_{2n+1}$}; \\
2^{\deg(Z)}, & \text{if\/ $\bfG=\rmO^+_{2n}$ and $Z$ degenerate}; \\
2^{\deg(Z)+1/2}, & \text{if\/ $\bfG=\SO^+_{2n}$ and $Z$ degenerate}; \\
2^{\deg(Z)-1/2}, & \text{if\/ $\bfG=\SO^\epsilon_{2n},\rmO^\epsilon_{2n}$ and $Z$ non-degenerate}.
\end{cases}
\]
\end{prop}
\begin{proof}
If $\bfG=\Sp_{2n}$ $\SO_{2n+1}$ or $\SO^\epsilon_{2n}$,
the result is from \cite{lg-symplectic} theorem 5.8 and \cite{lg-orthogonal} theorem 3.15.
If $\bfG=\rmO^\epsilon_{2n}$,
a proof can be found in \cite{pan-uniform} proposition 3.6.
Note that the definition of $R^\bfG_\Sigma$ here is slightly different from that in \cite{pan-uniform} \S 3.4.
\end{proof}

A bijective mapping $\cals_\bfG\rightarrow\cale(\bfG,1)$ satisfying (\ref{0332})
is called a \emph{Lusztig parametrization} of unipotent characters of $G$.
\begin{rem}
For $\bfG=\Sp_{2n}$, our definition $\cals_\bfG$ in (\ref{0210}) is slightly different from the original definition
by Lusztig in \cite{lg} p.134.
The definition here is more convenient for finite theta correspondence on unipotent characters.
Moreover, it is easy to see that two definitions are equivalent (\cf.~\cite{pan-uniform} (3.5)).
\end{rem}

\begin{lem}\label{0337}
Suppose that $\bfG=\Sp_{2n}$, $\SO_{2n+1}$, $\SO^\epsilon_{2n}$, or $\rmO^\epsilon_{2n}$ where $\epsilon=+$ or $-$.
Let $\call_1,\call'_1\colon\cals_\bfG\rightarrow\cale(\bfG,1)$ be two Lusztig parametrizations of
unipotent characters of $G$.
Then $\call_1(\Lambda)^\sharp=\call'_1(\Lambda)^\sharp$.
\end{lem}
\begin{proof}
Write $\call_1(\Lambda)=\rho_\Lambda$ and $\call'_1(\Lambda)=\rho'_\Lambda$ for $\Lambda\in\cals_\bfG$.
From (\ref{0332}), we see that
\[
\langle\rho_\Lambda,R^\bfG_\Sigma\rangle_\bfG=\langle\rho'_\Lambda,R^\bfG_\Sigma\rangle_\bfG
\]
for any $\Sigma\in\cals_\bfG^\sharp$.
Moreover, by Lemma~\ref{0315} we see that
\begin{equation}\label{0307}
\rho_\Lambda^\sharp=\sum_{\Sigma\in\cals_\bfG^\sharp}
\langle\rho_\Lambda,R^\bfG_\Sigma\rangle_\bfG R^\bfG_\Sigma.
\end{equation}
Then the lemma is proved.
\end{proof}

Suppose that $\bfG=\prod_{i=1}^k\bfG_i$ where each $\bfG_i$ is a classical group.
It is clear that
\begin{align*}
\cale(\bfG,1) &=\{\,\rho_1\otimes\cdots\otimes\rho_k\mid\rho_i\in\cale(\bfG_i,1)\,\}, \\
\calv(\bfG,1) &=\bigotimes_{i=1}^k\calv(\bfG_i,1), \\
\calv(\bfG,1)^\sharp &=\bigotimes_{i=1}^k\calv(\bfG_i,1)^\sharp.
\end{align*}
Then we define
\[
\cals_{\bfG}=\prod_{i=1}^k\cals_{\bfG_i},\qquad
\cals_{\bfG}^\sharp=\prod_{i=1}^k\cals_{\bfG_i}^\sharp,
\]
and then any Lusztig parametrization $\call_1\colon\cals_\bfG\rightarrow\cale(\bfG,1)$ is of the form
\begin{align*}
\call_1\colon\prod_{i=1}^k\cals_{\bfG_i} &\rightarrow\cale\Biggl(\prod_{i=1}^k\bfG_i,1\Biggr) \\
(\Lambda_1,\ldots,\Lambda_k) &\mapsto \rho_{\Lambda_1}\otimes\cdots\otimes\rho_{\Lambda_k}
\end{align*}
where $\Lambda_i\mapsto\rho_{\Lambda_i}$ is given by a Lusztig parametrization
$\cals_{\bfG_i}\rightarrow\cale(\bfG_i,1)$.

The parametrization $\call_1$ given by Lusztig in Proposition~\ref{0301} is compatible with parabolic induction
as follows.
Let $\bfG_n=\Sp_{2n}$, $\SO_{2n+1}$, $\SO^\epsilon_{2n}$ or $\rmO^\epsilon_{2n}$.
For $\rho\in\cale(\bfG_n,1)$, let $\Omega(\rho)\subset\cale(\bfG_{n+1},1)$ be defined similarly as
in (\ref{0323}).
For $\Lambda\in\cals_{\bfG_n}$,
then $\Omega(\Lambda)$ consists of all symbols $\Lambda'\in\cals_{\bfG_{n+1}}$ such that
\begin{itemize}
\item ${\rm def}(\Lambda')={\rm def}(\Lambda)$, and

\item $\Upsilon(\Lambda')$
is obtained from $\Upsilon(\Lambda)$ by adding a box to the Young diagram of $\Upsilon(\Lambda)^*$ or
$\Upsilon(\Lambda)_*$.
\end{itemize}

\begin{exam}
Suppose that $\bfG_n=\Sp_4$ and $\Lambda=\binom{2,0}{1}$,
then $\Upsilon(\Lambda)=\sqbinom{1}{1}$, and so
$\Upsilon(\Lambda')$ is equal to $\sqbinom{2}{1}$, $\sqbinom{1,1}{1}$, $\sqbinom{1}{2}$, or $\sqbinom{1}{1,1}$.
Therefore
\[
\Omega(\Lambda)=\left\{\textstyle\binom{3,0}{1},\binom{2,1}{1},\binom{2,0}{2},\binom{3,1,0}{2,1}\right\}
\subset\cals_{\Sp_6}.
\]
\end{exam}

For $\Sigma\in\cals_{\bfG_n}^\sharp$,
it is known that (\cf.~\cite{Geck-Pfeiffer} \S 6.1.9)
\[
\Ind_{W_n\times S_1}^{W_{n+1}}(\varphi_{\Upsilon(\Sigma)}\otimes{\bf 1})
=\sum_{\Sigma'\in\Omega(\Sigma)}\varphi_{\Upsilon(\Sigma')}.
\]
By direct computation (\cf.~\cite{lg-symplectic} (4.6.3)), we have
\[
R_{\bfG_n\times\GL_1}^{\bfG_{n+1}}(R_\Sigma^{\bfG_n})=\sum_{\Sigma'\in\Omega(\Sigma)}R_{\Sigma'}^{\bfG_{n+1}}.
\]

We say that the parametrization $\call_1\colon\cals_\bfG\rightarrow\cale(\bfG,1)$ is
\emph{compatible with parabolic induction} if the diagram analogous to (\ref{0324}) commutes, i.e.,
\[
\Omega(\rho_\Lambda)=\{\,\rho_{\Lambda'}\mid\Lambda'\in\Omega(\Lambda)\,\}.
\]
Note that ${\rm def}(\Lambda')={\rm def}(\Lambda)$ for any $\Lambda'\in\Omega(\Lambda)$.
This means that under the parametrization by Lusztig the defects of symbols
are preserved by parabolic induction on unipotent characters.
Therefore, if ${\rm def}(\Lambda')\neq 0$ and $\Lambda'\in\Omega(\Lambda)$,
then $\Lambda'^\rmt\not\in\Omega(\Lambda)$.

\begin{lem}\label{0328}
Let $\Lambda\in\cals_{\rmO^+_{2n}}$ such that ${\rm def}(\Lambda)=0$ and $\Lambda\neq\Lambda^\rmt$.
Suppose that $n\geq 2$.
Then there exists $\Lambda_1\in\cals_{\rmO^+_{2(n-1)}}$ such that
$\Lambda\in\Omega(\Lambda_1)$ and $\Lambda^\rmt\not\in\Omega(\Lambda_1)$.
\end{lem}
\begin{proof}
Write $\Lambda=\binom{a_1,\ldots,a_m}{b_1,\ldots,b_m}$ where $a_m,b_m$ are not both zero.
Let $i$ be the largest index such that $a_i\neq b_i$.
Such an index $i$ exists because we assume that $\Lambda\neq\Lambda^\rmt$.
Now we consider the following cases:
\begin{itemize}
\item Suppose that $i=m=1$.
We know that $\Lambda\neq\binom{1}{0},\binom{0}{1}$ because $\binom{1}{0},\binom{0}{1}\in\cals_{\rmO^+_2}$
which contradicts to the assumption $n\geq 2$.
\begin{itemize}
\item If either $a_1\geq 2$ and $b_1=0$; or $b_1>a_1\geq 1$,
let $\Lambda_1=\binom{a_1-1}{b_1}$.

\item If either $a_1=0$ and $b_1\geq 2$; or $a_1>b_1\geq 1$,
let $\Lambda_1=\binom{a_1}{b_1-1}$.
\end{itemize}

\item Suppose that $i=m>1$.
\begin{itemize}
\item If $a_m>b_m$ and $b_{m-1}\geq a_{m-1}$,
then we have $b_{m-1}> b_m+1$ and let
$\Lambda_1=\binom{a_1,\ldots,a_m}{b_1,\ldots,b_{m-2},b_{m-1}-1,b_m}$.

\item If $a_m>b_m$ and $a_{m-1}>b_{m-1}$,
let $\Lambda_1=\binom{a_1,\ldots,a_{m-1},a_m-1}{b_1,\ldots,b_m}$.

\item If $b_m>a_m$ and $a_{m-1}\geq b_{m-1}$,
let $\Lambda_1=\binom{a_1,\ldots,a_{m-2},a_{m-1}-1,a_m}{b_1,\ldots,b_m}$.

\item If $b_m>a_m$ and $b_{m-1}>a_{m-1}$,
let $\Lambda_1=\binom{a_1,\ldots,a_m}{b_1,\ldots,b_{m-1},b_m-1}$.
\end{itemize}

\item Suppose that $i<m$.
\begin{itemize}
\item If $a_i>b_i$, then $a_i>b_i>b_{i+1}=a_{i+1}$ and let
$\Lambda_1=\binom{a_1,\ldots,a_{i-1},a_i-1,a_{i+1},\ldots,a_m}{b_1,\ldots,b_m}$.

\item If $b_i>a_i$, then $b_i>a_i>a_{i+1}=b_{i+1}$ and let
$\Lambda_1=\binom{a_1,\ldots,a_m}{b_1,\ldots,b_{i-1},b_i-1,b_{i+1},\ldots,b_m}$.
\end{itemize}
\end{itemize}
For all cases, it is not difficult to check $\Lambda_1\in\cals_{\rmO^+_{2(n-1)}}$,
$\Lambda\in\Omega(\Lambda_1)$ and $\Lambda^\rmt\not\in\Omega(\Lambda_1)$.
\end{proof}

\begin{rem}\label{0336}
It is obvious that the statement in above lemma is not true without the assumption $n\geq 2$.
Note that $\Lambda=\binom{1}{0}\in\cals_{\rmO^+_2}$, $\Lambda\neq\Lambda^\rmt$.
However, $\cals_{\rmO^+_0}=\bigl\{\binom{-}{-}\bigr\}$ and $\Omega(\binom{-}{-})=\bigl\{\binom{1}{0},\binom{0}{1}\bigr\}$.
\end{rem}

\subsection{Cells in a family of unipotent characters}
In this subsection, we recall some result on \emph{cells} by Lusztig
(\cf.~\cite{lg-symplectic}, \cite{lg-orthogonal}).
Some details and examples can be found in \cite{pan-finite-unipotent} \S 4.
Let $Z$ be a special symbol of rank $n$.
An \emph{arrangement} $\Phi$ of $Z$ is defined as follows:
\begin{itemize}
\item if ${\rm def}(Z)=1$,
$\Phi$ is a partition of $Z_\rmI$ into $\deg(Z)$ pairs and one isolated element such that
each pair contains one element in $(Z_\rmI)^*$ and one element in $(Z_\rmI)_*$;

\item if ${\rm def}(Z)=0$,
$\Phi$ is a partition of $Z_\rmI$ into $\deg(Z)$ pairs such that
each pair contains one element in $(Z_\rmI)^*$ and one element in $(Z_\rmI)_*$.
\end{itemize}
A subset of pairs $\Psi$ in an arrangement $\Phi$ is denoted by $\Psi\leq\Phi$.
For such an arrangement $\Phi$ and a subset of pairs $\Psi$, we define a subset
$C_{\Phi,\Psi}^\bfG=C_{\Phi,\Psi}$ of $\cals_\bfG$
as follows:
\begin{itemize}
\item if ${\rm def}(Z)=1$, we define
\[
C_{\Phi,\Psi}=\{\,\Lambda_M\mid M\subset Z_\rmI,\ |M|\text{ even},\
|M\cap\Psi'|\equiv|(\Phi\smallsetminus\Psi)\cap\Psi'^*|\pmod{2}\text{ for all }\Psi'\leq\Phi\,\};
\]

\item if ${\rm def}(Z)=0$, we define
\[
C_{\Phi,\Psi}=\{\,\Lambda_M\mid M\subset Z_\rmI,\
|M\cap\Psi'|\equiv|(\Phi\smallsetminus\Psi)\cap\Psi'^*|\pmod{2}\text{ for all }\Psi'\leq\Phi\,\}
\]
\end{itemize}
It is not difficult to see that for ,
\[
C_{\Phi,\Psi}\subset\begin{cases}
\cals_Z^{\Sp_{2n}}, & \text{if ${\rm def}(Z)=1$};\\
\cals_Z^{\rmO^+_{2n}}, & \text{if ${\rm def}(Z)=0$ and $\#(\Phi\smallsetminus\Psi)$ even};\\
\cals_Z^{\rmO^-_{2n}}, & \text{if ${\rm def}(Z)=0$ and $\#(\Phi\smallsetminus\Psi)$ odd}.
\end{cases}
\]
Here $\#(\Phi\smallsetminus\Psi)$ denotes the number of pairs in $\Phi\smallsetminus\Psi$.
For $\bfG=\rmO^\epsilon_{2n}$, a special symbol $Z$ of rank $n$ and defect $0$,
and an arrangement $\Phi$ of $Z$,
a subset of pairs $\Psi$ is called \emph{admissible} if $\#(\Phi\smallsetminus\Psi)$ is even
when $\epsilon=+$; and $\#(\Phi\smallsetminus\Psi)$ is odd when $\epsilon=-$.

The following lemmas are from \cite{pan-finite-unipotent} lemma~4.17, proposition~4.18, lemma~4.34,
and proposition~4.35:

\begin{lem}\label{0317}
Let $\bfG=\Sp_{2n}$,
and let $\Lambda\mapsto\rho_\Lambda$ be a Lusztig parametrization of unipotent characters.
Let $Z$ be a special symbol of rank $n$ and defect $1$,
$\Phi$ an arrangement of $Z$, $\Psi$ a subset of pairs of\/ $\Phi$.
\begin{enumerate}
\item[(i)] The class function $\sum_{\Lambda\in C_{\Phi,\Psi}}\rho_\Lambda$ is uniform.

\item[(ii)] For any two distinct symbols $\Lambda_1,\Lambda_2\in\cals_Z^\bfG$,
there exists an arrangement $\Phi$ of $Z$ with two subsets of pairs $\Psi_1,\Psi_2$ such that
$\Lambda_i\in C_{\Phi,\Psi_i}$ for $i=1,2$ and $C_{\Phi,\Psi_1}\cap C_{\Phi,\Psi_2}=\emptyset$.
\end{enumerate}
\end{lem}

\begin{lem}\label{0318}
Let $\bfG=\rmO^\epsilon_{2n}$ where $\epsilon=+$ or $-$,
and let $\Lambda\mapsto\rho_\Lambda$ be a Lusztig parametrization of unipotent characters.
Let $Z$ be a special symbol of rank $n$ and defect $0$,
$\Phi$ an arrangement of $Z$, $\Psi$ an admissible subset of pairs of\/ $\Phi$.
\begin{enumerate}
\item[(i)] The class function $\sum_{\Lambda\in C_{\Phi,\Psi}}\rho_\Lambda$ is uniform.

\item[(ii)] $\Lambda\in C_{\Phi,\Psi}$ if and only if $\Lambda^\rmt\in C_{\Phi,\Psi}$.

\item[(iii)] For any two symbols $\Lambda_1,\Lambda_2\in\cals_Z^\bfG$ such that
$\Lambda_1\neq\Lambda_2,\Lambda_2^\rmt$,
there exists an arrangement $\Phi$ of $Z$ with subsets of pairs $\Psi_1,\Psi_2$ such that
$\Lambda_i\in C_{\Phi,\Psi_i}$ for $i=1,2$ and $C_{\Phi,\Psi_1}\cap C_{\Phi,\Psi_2}=\emptyset$.
\end{enumerate}
\end{lem}

\section{Uniqueness of the Lusztig Parametrizations}

\subsection{Unipotent characters of $\Sp_{2n}$}
In this subsection, let $\bfG=\Sp_{2n}$.
The following lemma is \cite{pan-Lusztig-correspondence} proposition 3.3.
We provide a proof for the sake of completion.

\begin{lem}\label{0302}
Suppose that $\bfG=\Sp_{2n}$, and let $\rho_1,\rho_2\in\cale(\bfG,1)$.
If $\rho_1^\sharp=\rho_2^\sharp$,
then $\rho_1=\rho_2$.
\end{lem}
\begin{proof}
Let $\call_1\colon\cals_\bfG\rightarrow\cale(\bfG,1)$ be a Lusztig parametrization given in
Proposition~\ref{0301}, and write $\rho_1=\rho_{\Lambda_1}, \rho_2=\rho_{\Lambda_2}$
for some $\Lambda_1,\Lambda_2\in\cals_\bfG$.
From (\ref{0332}), the assumption $\rho_{\Lambda_1}^\sharp=\rho_{\Lambda_2}^\sharp$ means that
$\Lambda_1,\Lambda_2\in\cals_Z^\bfG$ for some special symbol $Z$ of rank $n$ and defect $1$.
For any arrangement $\Phi$ of $Z$ and any subset of pairs $\Psi\leq\Phi$,
we know that $\sum_{\Lambda\in C_{\Phi,\Psi}}\rho_\Lambda$ is uniform by Lemma~\ref{0317},
and so
\begin{align*}
\Biggl\langle\sum_{\Lambda\in C_{\Phi,\Psi}}\rho_\Lambda,\rho_{\Lambda_1}\Biggr\rangle_{\!\!\bfG}
=\Biggl\langle\sum_{\Lambda\in C_{\Phi,\Psi}}\rho_\Lambda,\rho_{\Lambda_1}^\sharp\Biggr\rangle_{\!\!\bfG}
=\Biggl\langle\sum_{\Lambda\in C_{\Phi,\Psi}}\rho_\Lambda,\rho_{\Lambda_2}^\sharp\Biggr\rangle_{\!\!\bfG}
=\Biggl\langle\sum_{\Lambda\in C_{\Phi,\Psi}}\rho_\Lambda,\rho_{\Lambda_2}\Biggr\rangle_{\!\!\bfG}.
\end{align*}
Now if $\Lambda_1\neq\Lambda_2$,
by Lemma~\ref{0317}
we can find an arrangement $\Phi$ of $Z$ and a subset of pairs $\Psi\leq\Phi$ such that
$\Lambda_1\in C_{\Phi,\Psi}$ and $\Lambda_2\not\in C_{\Phi,\Psi}$,
and hence
\[
\Biggl\langle\sum_{\Lambda\in C_{\Phi,\Psi}}\rho_\Lambda,\rho_{\Lambda_1}\Biggr\rangle_{\!\!\bfG}=1\quad\text{and}\quad
\Biggl\langle\sum_{\Lambda\in C_{\Phi,\Psi}}\rho_\Lambda,\rho_{\Lambda_2}\Biggr\rangle_{\!\!\bfG}=0.
\]
So we must have $\Lambda_1=\Lambda_2$, i.e., $\rho_1=\rho_2$.
\end{proof}

\begin{prop}\label{0308}
Let $\bfG$ be $\Sp_{2n}$.
Then there is a unique bijection $\call_1\colon\cals_\bfG\rightarrow\cale(\bfG,1)$ satisfying
(\ref{0332}).
\end{prop}
\begin{proof}
Suppose that we have two parametrizations $\Lambda\mapsto\rho_\Lambda$ and
$\Lambda\mapsto\rho'_\Lambda$ from $\cals_{\bfG}$ to $\cale(\bfG,1)$ satisfying (\ref{0332}).
From (\ref{0307}), we see that condition (\ref{0332}) implies that $(\rho_\Lambda)^\sharp=(\rho'_\Lambda)^\sharp$.
Then by Lemma~\ref{0302}, we conclude that $\rho_\Lambda=\rho'_\Lambda$, i.e.,
two parametrizations coincide.
\end{proof}

\begin{cor}\label{0334}
Let $\bfG_n=\Sp_{2n}$.
Then the bijection $\call_1\colon\cals_{\bfG_n}\rightarrow\cale(\bfG_n,1)$ given in Proposition~\ref{0308}
is compatible with the parabolic induction, i.e.,
the diagram analogous to (\ref{0324}) commutes.
\end{cor}
\begin{proof}
The original construction of the bijection $\cals_{\bfG_n}\rightarrow\cale(\bfG_n,1)$ by Lusztig
is compatible with the parabolic induction (\cf.~\cite{lg}).
By Proposition~\ref{0308}, Lusztig's original construction is the only bijection satisfying (\ref{0332})
and hence the corollary is obtained.
\end{proof}

For a non-negative integer $k$, we define the symbol
\begin{equation}\label{0313}
\Lambda_k^\Sp=\begin{cases}
\binom{2k,2k-1,\ldots,0}{-}, & \text{if $k$ is even};\\
\binom{-}{2k,2k-1,\ldots,0}, & \text{if $k$ is odd}.
\end{cases}
\end{equation}
The following are easy to check:
\begin{itemize}
\item ${\rm rk}(\Lambda^\Sp_k)=k(k+1)$ and ${\rm def}(\Lambda^\Sp_k)\equiv 1\pmod 4$,
i.e., $\Lambda_k^\Sp\in\cals_{\Sp_{2k(k+1)}}$,

\item if $\Lambda\in\cals_{\Sp_{2n}}$ with $n<k(k+1)$,
then $|{\rm def}(\Lambda)|<|{\rm def}(\Lambda^\Sp_k)|$,

\item if $\Lambda\in\cals_{\Sp_{2k(k+1)}}$ and $\Lambda\neq\Lambda^\Sp_k$,
then $|{\rm def}(\Lambda)|<|{\rm def}(\Lambda^\Sp_k)|$.
\end{itemize}
Because the defects are preserved by the parabolic induction, we have the following corollary:

\begin{cor}\label{0311}
Let $\call_1\colon\cals_{\Sp_{2n}}\rightarrow\cale(\Sp_{2n},1)$ be the parametrization
in Proposition~\ref{0308}.
Then the unique cuspidal unipotent character $\zeta_k^\Sp$ of\/ $\Sp_{2k(k+1)}(q)$ is parametrized by the symbol
$\Lambda^\Sp_k$, i.e., $\zeta^\Sp_k=\rho_{\Lambda^\Sp_k}$.
\end{cor}

\begin{lem}\label{0312}
Let $\call_1\colon\cals_{\Sp_{2n}}\rightarrow\cale(\Sp_{2n},1)$ be the parametrization
in Proposition~\ref{0308}.
Then
\begin{enumerate}
\item[(i)] ${\bf 1}_{\Sp_{2n}}=\rho_{\binom{n}{-}}$,

\item[(ii)] ${\rm St}_{\Sp_{2n}}=\rho_\Lambda$ where $\Lambda=\binom{n,n-1,\ldots,1,0}{n,n-1,\ldots,1}$.
\end{enumerate}
\end{lem}
\begin{proof}
From \cite{carter-finite} corollary 7.6.5, we know that
$R^{\Sp_{2n}}_{\binom{n}{-}}={\bf 1}_{\Sp_{2n}}$.
Because now $\binom{n}{-}$ is a special symbol of degree $0$, we have
$\cals_{\binom{n}{-}}^{\Sp_{2n}}=\bigl\{\binom{n}{-}\bigr\}$,
$\rho_{\binom{n}{-}}=R^{\Sp_{2n}}_{\binom{n}{-}}$ by (\ref{0332}),
and so (i) is proved.

Write $\Lambda=\binom{n,n-1,\ldots,1,0}{n,n-1,\ldots,1}$.
From \cite{carter-finite} corollary 7.6.6, we see that
$R^{\Sp_{2n}}_{\Lambda}={\rm St}_{\Sp_{2n}}$.
Again, now $\Lambda$ is a special symbol of degree $0$, we have
$\rho_\Lambda=R^{\Sp_{2n}}_\Lambda$, and (ii) is proved.
\end{proof}

\begin{exam}\label{0333}
Let $\bfG=\Sp_4$.
We know that
\begin{align*}
\cals_{\Sp_4} &=\bigl\{\textstyle\binom{2}{-}\bigr\}\cup\bigl\{\binom{2,1,0}{2,1}\bigr\}\cup
\bigl\{\binom{2,1}{0}, \binom{2,0}{1}, \binom{1,0}{2}, \binom{-}{2,1,0}\bigr\}, \\
\cals_{\Sp_4}^\sharp &=\bigl\{\textstyle\binom{2}{-}\bigr\}\cup\bigl\{\binom{2,1,0}{2,1}\bigr\}
\cup\bigl\{\binom{2,1}{0}, \binom{2,0}{1}, \binom{1,0}{2}\bigr\}.
\end{align*}
The character values of irreducible characters in
$\cale(W_2)=\left\{\,\varphi_{\sqbinom{\mu}{\nu}}\mid\sqbinom{\mu}{\nu}\in\calp_2(2)\,\right\}$
are given by the following table:
\[
\begin{tabular}{c|ccccc}
\toprule
$(w)$ & $\{1\}$ & $\{\sigma_2,s_1\sigma_2s_1\}$ & $\{s_1\sigma_2s_1\sigma_2\}$ & $\{s_1,\sigma_2s_1\sigma_2\}$
& $\{s_1\sigma_2,\sigma_2s_1\}$ \\
\midrule
$\varphi_{\sqbinom{2}{0}}$ & $1$ & $1$ & $1$ & $1$ & $1$ \\
$\varphi_{\sqbinom{1,1}{0}}$ & $1$ & $1$ & $1$ & $-1$ & $-1$ \\
$\varphi_{\sqbinom{1}{1}}$ & $2$ & $0$ & $-2$ & $0$ & $0$ \\
$\varphi_{\sqbinom{0}{2}}$ & $1$ & $-1$ & $1$ & $1$ & $-1$ \\
$\varphi_{\sqbinom{0}{1,1}}$ & $1$ & $-1$ & $1$ & $-1$ & $1$ \\
\midrule
$\bfT_w$ & $\bfT_1$ & $\bfT_2$ & $\bfT_3$ & $\bfT_4$ & $\bfT_5$ \\
\bottomrule
\end{tabular}
\]
where $s_1,\sigma_2$ are defined in Subsection~\ref{0209}.
We can check that
\begin{align*}
R_{\bfT_1,1} &= 2\theta_9 +\theta_{11} +\theta_{12} +\theta_{13} +\theta_0,\\
R_{\bfT_2,1} &= \theta_{11} -\theta_{12} -\theta_{13} +\theta_0, \\
R_{\bfT_3,1} &= -2\theta_{10} -\theta_{11} -\theta_{12} +\theta_{13} +\theta_0,\\
R_{\bfT_4,1} &= -\theta_{11} +\theta_{12} -\theta_{13} +\theta_0,\\
R_{\bfT_5,1} &= -\theta_9 +\theta_{10} +\theta_{13} +\theta_0
\end{align*}
where $\theta_i$ is the notion from \cite{srinivasan-sp4}.
Therefore by (\ref{0330}), we have
\begin{align*}
R^\bfG_{\binom{2}{-}} &=\frac{1}{8}\left[R_{\bfT_1,1}+2R_{\bfT_2,1}+R_{\bfT_3,1}+2R_{\bfT_4,1}+2R_{\bfT_5,1}\right]
=\theta_0, \\
R^\bfG_{\binom{2,0}{1}} &=\frac{1}{4}\left[R_{\bfT_1,1}-R_{\bfT_3,1}\right]
=\frac{1}{2}(\theta_9+\theta_{10}+\theta_{11}+\theta_{12}), \\
R^\bfG_{\binom{2,1}{0}} &=\frac{1}{8}\left[R_{\bfT_1,1}+2R_{\bfT_2,1}+R_{\bfT_3,1}-2R_{\bfT_4,1}-2R_{\bfT_5,1}\right]
=\frac{1}{2}(\theta_9-\theta_{10}+\theta_{11}-\theta_{12}), \\
R^\bfG_{\binom{1,0}{2}} &=\frac{1}{8}\left[R_{\bfT_1,1}-2R_{\bfT_2,1}+R_{\bfT_3,1}+2R_{\bfT_4,1}-2R_{\bfT_5,1}\right]
=\frac{1}{2}(\theta_9-\theta_{10}-\theta_{11}+\theta_{12}), \\
R^\bfG_{\binom{2,1,0}{2,1}} &=\frac{1}{8}\left[R_{\bfT_1,1}-2R_{\bfT_2,1}+R_{\bfT_3,1}-2R_{\bfT_4,1}+2R_{\bfT_5,1}\right]
=\theta_{13}.
\end{align*}
It is known that $\theta_0={\bf 1}_{\Sp_4}=\rho_{\binom{2}{-}}$ and 
$\theta_{13}={\rm St}_{\Sp_4}=\rho_{\binom{2,1,0}{2,1}}$ by Lemma~\ref{0312}.
Let $Z=\binom{2,0}{1}$.
The table for $(-1)^{\langle\Sigma,\Lambda\rangle}$ for $\Sigma\in\cals_{Z,1}$ and $\Lambda\in\cals_Z^\bfG$ is
\[
\begin{tabular}{c|cccc}
\toprule
& $\binom{2,0}{1}$ & $\binom{2,1}{0}$ & $\binom{1,0}{2}$ & $\binom{-}{2,1,0}$ \\
\midrule
$\binom{2,0}{1}$ & $\phantom{-}1$ & $\phantom{-}1$ & $\phantom{-}1$ & $\phantom{-}1$ \\
$\binom{2,1}{0}$ & $\phantom{-}1$ & $\phantom{-}1$ & $-1$ & $-1$ \\
$\binom{1,0}{2}$ & $\phantom{-}1$ & $-1$ & $\phantom{-}1$ & $-1$ \\
\bottomrule
\end{tabular}
\]
Then we have
\begin{align*}
{\rho_{\binom{2,0}{1}}}^\sharp &=\frac{1}{2}\left[R_{\binom{2,0}{1}}+R_{\binom{2,1}{0}}+R_{\binom{1,0}{2}}\right], &
{\rho_{\binom{2,1}{0}}}^\sharp &=\frac{1}{2}\left[R_{\binom{2,0}{1}}+R_{\binom{2,1}{0}}-R_{\binom{1,0}{2}}\right], \\
{\rho_{\binom{1,0}{2}}}^\sharp &=\frac{1}{2}\left[R_{\binom{2,0}{1}}-R_{\binom{2,1}{0}}+R_{\binom{1,0}{2}}\right], & {\rho_{\binom{-}{2,1,0}}}^\sharp &=\frac{1}{2}\left[R_{\binom{2,0}{1}}-R_{\binom{2,1}{0}}-R_{\binom{1,0}{2}}\right].
\end{align*}
We know that
\begin{align*}
R^{\Sp_4}_{\Sp_2\times\GL_1}(\rho_{\binom{1}{-}}\otimes{\bf 1})
&={\bf1}_{\Sp_4}+\theta_{11}+\theta_9, &
R^{\Sp_4}_{\Sp_2\times\GL_1}(\rho_{\binom{1,0}{1}}\otimes{\bf 1})
&=\theta_9+\theta_{12}+{\rm St}_{\Sp_4},\\
\textstyle\Omega(\binom{1}{-}) &= \textstyle\left\{\binom{2}{-},\binom{2,1}{0},\binom{2,0}{1}\right\}, &
\textstyle\Omega(\binom{1,0}{1}) &= \textstyle\left\{\binom{2,0}{1},\binom{1,0}{2},\binom{2,1,0}{2,1}\right\}.
\end{align*}
Then by Corollary~\ref{0334}, we conclude that $\theta_9=\rho_{\binom{2,0}{1}}$,
$\theta_{10}=\rho_{\binom{-}{2,1,0}}$, $\theta_{11}=\rho_{\binom{2,1}{0}}$,
and $\theta_{12}=\rho_{\binom{1,0}{2}}$.
\end{exam}

\subsection{Unipotent characters of $\SO_{2n+1}$}

\begin{lem}\label{0322}
Suppose that $\bfG=\SO_{2n+1}$, and let $\rho_1,\rho_2\in\cale(\bfG,1)$.
If $\rho_1^\sharp=\rho_2^\sharp$,
then $\rho_1=\rho_2$.
\end{lem}
\begin{proof}
The proof is similar to that of Lemma~\ref{0302}.
\end{proof}

\begin{prop}\label{0309}
Suppose that $\bfG=\SO_{2n+1}$.
Then there is a unique bijection $\cals_\bfG\rightarrow\cale(\bfG,1)$ satisfying (\ref{0332}).
\end{prop}
\begin{proof}
The proof is similar to that of Proposition~\ref{0308}.
\end{proof}

\begin{cor}
Let $\call_1\colon\cals_{\SO_{2n+1}}\rightarrow\cale(\SO_{2n+1},1)$ be the parametrization
given in Proposition~\ref{0309}.
Then
\begin{itemize}
\item[(i)] ${\bf 1}_{\SO_{2n+1}}=\rho_{\binom{-}{n}}$;

\item[(ii)] ${\rm St}_{\SO_{2n+1}}=\rho_\Lambda$ where $\Lambda=\binom{n,n-1,\ldots,1}{n,n-1,\ldots,1,0}$;

\item[(iii)] if $n=k(k+1)$ for some non-negative integer $k$, then the unique cuspidal unipotent character
$\zeta^{\SO_{\rm odd}}_k$ of\/ $\SO_{2n+1}(q)$ is parametrized by the symbol
\[
\Lambda^{\SO_{\rm odd}}_k=\begin{cases}
\binom{-}{2k,2k-1,\ldots,0}, & \text{if $k$ is even}; \\
\binom{2k,2k-1,\ldots,0}{-}, & \text{if $k$ is odd}.
\end{cases}
\]
\end{itemize}
\end{cor}
\begin{proof}
For (i) and (ii), the proofs are analogous to that of Corollary~\ref{0311},
for (iii) the proof is analogous to that of Lemma~\ref{0312}.
\end{proof}

\subsection{Unipotent characters of $\rmO^\epsilon_{2n}$}\label{0327}
From (\ref{0303}), we know that
$R^{\rmO^\epsilon_{2n}}_{\bfT,1}\cdot\sgn_{\rmO^\epsilon_{2n}}=R^{\rmO^\epsilon_{2n}}_{\bfT,1}$ and then
\[
\langle\rho,R^{\rmO^\epsilon_{2n}}_{\bfT,1}\rangle_{\rmO^\epsilon_{2n}}
=\langle\rho\cdot\sgn_{\rmO^\epsilon_{2n}},R^{\rmO^\epsilon_{2n}}_{\bfT,1}\rangle_{\rmO^\epsilon_{2n}}
\]
for any $\rho\in\cale(\rmO^\epsilon_{2n},1)$.
Therefore,
\begin{equation}\label{0320}
\rho^\sharp=(\rho\cdot\sgn_{\rmO^\epsilon_{2n}})^\sharp,
\end{equation}
i.e., two irreducible characters $\rho,\rho\cdot\sgn_{\rmO^\epsilon_{2n}}$ are not able to be distinguished
by their uniform projections.

The following lemma is \cite{pan-Lusztig-correspondence} proposition 3.5.
We also provide a proof here:
\begin{lem}\label{0321}
Let $\call_1\colon\cals_{\rmO^\epsilon_{2n}}\rightarrow\cale(\rmO^\epsilon_{2n},1)$ by
$\Lambda\mapsto\rho_\Lambda$ be a Lusztig parametrization of unipotent characters.
Then $(\rho_{\Lambda_1})^\sharp=(\rho_{\Lambda_2})^\sharp$ if and only if
$\Lambda_1=\Lambda_2$ or $\Lambda_1=\Lambda_2^\rmt$.
\end{lem}
\begin{proof}
For $\Lambda\in\cals_{\rmO^\epsilon_{2n}}$ and $\Sigma\in\cals_{\rmO^\epsilon_{2n}}^\sharp$,
we have $\langle\Lambda,\Sigma\rangle=\langle\Lambda^\rmt,\Sigma\rangle$ by Lemma~\ref{0203}.
Then we have $(\rho_{\Lambda^\rmt})^\sharp=(\rho_\Lambda)^\sharp$.
Now if $\Lambda_1\neq\Lambda_2,\Lambda_2^\rmt$,
by Lemma~\ref{0318}
we can find an arrangement $\Phi$ of $Z$ and an admissible subset of pairs $\Psi\leq\Phi$ such that
$\Lambda_1\in C_{\Phi,\Psi}$ and $\Lambda_2\not\in C_{\Phi,\Psi}$,
and hence
\[
\Biggl\langle\sum_{\Lambda\in C_{\Phi,\Psi}}\rho_\Lambda,\rho_{\Lambda_1}\Biggr\rangle_{\!\!\bfG}=1\quad\text{and}\quad
\Biggl\langle\sum_{\Lambda\in C_{\Phi,\Psi}}\rho_\Lambda,\rho_{\Lambda_2}\Biggr\rangle_{\!\!\bfG}=0.
\]
As in the proof of Lemma~\ref{0302}, we get a contradiction.
\end{proof}

\begin{cor}\label{0406}
Let $\call_1\colon\cals_{\rmO^\epsilon_{2n}}\rightarrow\cale(\rmO^\epsilon_{2n},1)$ by
$\Lambda\mapsto\rho_\Lambda$ be a Lusztig parametrization of unipotent characters.
Then $\rho_{\Lambda^\rmt}=\rho_\Lambda\cdot\sgn_{\rmO^\epsilon_{2n}}$.
\end{cor}
\begin{proof}
If $\Lambda$ is degenerate, then clearly
$\rho_{\Lambda^\rmt}=\rho_\Lambda=R^{\rmO^\epsilon_{2n}}_\Lambda
=R^{\rmO^\epsilon_{2n}}_\Lambda\cdot\sgn_{\rmO^\epsilon_{2n}}=\rho_\Lambda\cdot\sgn_{\rmO^\epsilon_{2n}}$.
If $\Lambda$ is non-degenerate,
from (\ref{0320}) we know that $(\rho_\Lambda\cdot\sgn_{\rmO^\epsilon_{2n}})^\sharp=(\rho_\Lambda)^\sharp$,
and by Lemma~\ref{0321} we conclude that $\rho_\Lambda\cdot\sgn_{\rmO^\epsilon_{2n}}=\rho_{\Lambda^\rmt}$.
\end{proof}

\begin{cor}\label{0314}
Let $\rho_1,\rho_2\in\cale(\rmO^\epsilon_{2n},1)$.
If $\rho_1^\sharp=\rho_2^\sharp$,
then either $\rho_1=\rho_2$ or $\rho_1=\rho_2\cdot\sgn_{\rmO^\epsilon_{2n}}$.
\end{cor}

\begin{cor}\label{0407}
Let $\bfG=\rmO^\epsilon_{2n}$, and let $\Lambda\mapsto\rho_\Lambda$ be a Lusztig parametrization
of unipotent characters.
Then any bijective mapping $\cals_\bfG\mapsto\cale(\bfG,1)$ such that
$\{\Lambda,\Lambda^\rmt\}\rightarrow\{\rho_\Lambda,\rho_{\Lambda^\rmt}\}$ is also a Lusztig parametrization
of unipotent characters.
\end{cor}
\begin{proof}
Suppose that $\cals_\bfG\mapsto\cale(\bfG,1)$ given by $\Lambda\mapsto\rho'_\Lambda$ is a bijection such that
$\{\Lambda,\Lambda^\rmt\}\rightarrow\{\rho_\Lambda,\rho_{\Lambda^\rmt}\}$, i.e.,
$\{\rho'_\Lambda,\rho'_{\Lambda^\rmt}\}=\{\rho_\Lambda,\rho_{\Lambda^\rmt}\}$.
This implies that $(\rho'_\Lambda)^\sharp=(\rho'_{\Lambda^\rmt})^\sharp=(\rho_\Lambda)^\sharp=(\rho_{\Lambda^\rmt})^\sharp$
and hence the mapping $\Lambda\mapsto\rho'_\Lambda$ satisfies (\ref{0332}),
i.e., $\Lambda\mapsto\rho'_\Lambda$ is also a Lusztig parametrization of unipotent characters for $\rmO^\epsilon_{2n}$.
\end{proof}

\begin{cor}\label{0310}
Let $\call_1\colon\cals_{\rmO^\epsilon_{2n}}\rightarrow\cale(\rmO^\epsilon_{2n},1)$ be a Lusztig parametrization
of unipotent characters.
\begin{enumerate}
\item[(i)] If $\epsilon=+$, then
\begin{enumerate}
\item $\call_1\colon\left\{\binom{n}{0},\binom{0}{n}\right\}\rightarrow\{{\bf1}_{\rmO^+_{2n}},\sgn_{\rmO^+_{2n}}\}$,

\item $\call_1$ maps $\bigl\{\binom{n,n-1,\ldots,1}{n-1,n-2,\ldots,0},\binom{n-1,n-2,\ldots,0}{n,n-1,\ldots,1}\bigr\}$
to the two Steinberg characters of\/ $\rmO^+_{2n}$.
\end{enumerate}
\item[(ii)] If $\epsilon=-$, then
\begin{enumerate}
\item $\call_1\colon\bigl\{\binom{-}{n,0},\binom{n,0}{-}\bigr\}\rightarrow\{{\bf1}_{\rmO^-_{2n}},\sgn_{\rmO^-_{2n}}\}$

\item $\call_1$ maps $\bigl\{\binom{n,n-1,\ldots,1,0}{n-1,n-2,\ldots,1},\binom{n-1,n-2,\ldots,1}{n,n-1,\ldots,1,0}\bigr\}$
to the two Steinberg characters of\/ $\rmO^-_{2n}$.
\end{enumerate}
\end{enumerate}
\end{cor}
\begin{proof}
First suppose that $\epsilon=+$.
From \cite{carter-finite} corollary 7.6.5, we know that
$R_{\binom{n}{0}}^{\SO_{2n}^+}={\bf 1}_{\SO^+_{2n}}$.
Therefore,
\[
R_{\binom{n}{0}}^{\rmO_{2n}^+}=R_{\binom{0}{n}}^{\rmO_{2n}^+}
=\frac{1}{\sqrt 2}({\bf 1}_{\rmO^+_{2n}}+\sgn_{\rmO^+_{2n}})
\]
and (i.a) is proved from (\ref{0332}).
Write $\Lambda=\binom{n-1,n-2,\ldots,0}{n,n-1,\ldots,1}$.
From \cite{carter-finite} corollary 7.6.6, we see that
$R_{\Lambda}^{\SO^+_{2n}}={\rm St}_{\SO^+_{2n}}$.
Then $\sqrt 2 R_{\Lambda}^{\rmO^+_{2n}}=\sqrt 2 R_{\Lambda^\rmt}^{\rmO^+_{2n}}$ 
is the sum of two Steinberg characters of $\rmO^+_{2n}$ and (i.b) is proved.

The proof of (ii) is similar.
\end{proof}

It is known that $\rmO^{\epsilon_k}_{2k^2}$ where $\epsilon_k=(-1)^k$ has two cuspidal unipotent characters,
denoted by $\zeta_k^{\rm I}$ and $\zeta_k^{\rm II}$.
Then from above we see that any Lusztig parametrization
$\call_1\colon\cals_{\rmO^{\epsilon_k}_{2k^2}}\rightarrow\cale(\rmO^{\epsilon_k}_{2k^2},1)$ maps
\[
\left\{\textstyle\binom{2k-1,2k-2,\ldots,1,0}{-},\binom{-}{2k-1,2k-2,\ldots,1,0}\right\}
\rightarrow\{\zeta_k^{\rm I},\zeta_k^{\rm II}\}
\]
bijectively.

\begin{exam}\label{0304}
We consider the unipotent characters of $\rmO^+_4$ or $\rmO^-_4$.
Let $\bfT_i$ be given as in Example~\ref{0333}.
It is know that $\bfT_1,\bfT_3,\bfT_4$ are maximal tori in $\rmO^+_4$,
and $\bfT_2,\bfT_5$ are maximal tori in $\rmO^-_4$.
It is know that
\begin{align*}
R^{\rmO^+_4}_{\bfT_1,1} &={\bf 1}_{\rmO^+_4}+\sgn_{\rmO^+_4}+2\chi^+_{2q}+\chi^+_{q^2}+\chi^+_{q^2}\cdot\sgn_{\rmO^+_4}, \\
R^{\rmO^+_4}_{\bfT_3,1} &={\bf 1}_{\rmO^+_4}+\sgn_{\rmO^+_4}-2\chi^+_{2q}+\chi^+_{q^2}+\chi^+_{q^2}\cdot\sgn_{\rmO^+_4}, \\
R^{\rmO^+_4}_{\bfT_4,1} &={\bf 1}_{\rmO^+_4}+\sgn_{\rmO^+_4}-\chi^+_{q^2}-\chi^+_{q^2}\cdot\sgn_{\rmO^+_4}, \\
R^{\rmO^-_4}_{\bfT_2,1} &={\bf 1}_{\rmO^-_4}+\sgn_{\rmO^-_4}+\chi^-_{q^2}+\chi^-_{q^2}\cdot\sgn_{\rmO^-_4}, \\
R^{\rmO^-_4}_{\bfT_5,1} &={\bf 1}_{\rmO^-_4}+\sgn_{\rmO^-_4}-\chi^-_{q^2}-\chi^-_{q^2}\cdot\sgn_{\rmO^-_4}
\end{align*}
where $\chi_{2q}^+,\chi_{q^2}^+$ are irreducible characters of $\rmO^+_4(q)$ of degrees $2q,q^2$ respectively;
similarly $\chi_{q^2}^-$ is an irreducible character of $\rmO^-_4(q)$ of degree $q^2$.
And so we have
\begin{align*}
\cale(\rmO^+_4,1) &=\{{\bf 1}_{\rmO^+_4},\sgn_{\rmO^+_4},\chi^+_{2q},\chi^+_{q^2},\chi^+_{q^2}\cdot\sgn_{\rmO^+_4}\}, \\
\cale(\rmO^-_4,1) &=\{{\bf 1}_{\rmO^-_4},\sgn_{\rmO^-_4},\chi^-_{q^2},\chi^-_{q^2}\cdot\sgn_{\rmO^-_4}\}.
\end{align*}
We know that
\[
\textstyle
\cals_{\rmO^+_4}=\left\{\binom{2}{0},\binom{0}{2},\binom{2,1}{1,0},\binom{1,0}{2,1},\binom{1}{1}\right\},\quad
\cals_{\rmO^-_4}=\left\{\binom{-}{2,0},\binom{2,0}{-},\binom{1}{2,1,0},\binom{2,1,0}{1}\right\}.
\]
Suppose we choose
$\cals_{\rmO^+_4}^\sharp=\left\{\binom{2}{0},\binom{2,1}{1,0},\binom{1}{1}\right\}$,
and $\cals_{\rmO^-_4}^\sharp=\left\{\binom{0}{2},\binom{1,0}{2,1}\right\}$.
Then we have
\begin{align*}
R^{\rmO^+_4}_{\binom{2}{0}}
&=\frac{1}{4\sqrt 2}\left[R^{\rmO^+_4}_{\bfT_1,1}+R^{\rmO^+_4}_{\bfT_3,1}+2R^{\rmO^+_4}_{\bfT_4,1}\right]
=\frac{1}{\sqrt 2}\left[{\bf 1}_{\rmO^+_4}+\sgn_{\rmO^+_4}\right], \\
R^{\rmO^+_4}_{\binom{1}{1}}
&=\frac{1}{4}\left[R^{\rmO^+_4}_{\bfT_1,1}-R^{\rmO^+_4}_{\bfT_3,1}\right]
=\chi^+_{2q}, \\
R^{\rmO^+_4}_{\binom{2,1}{1,0}}
&=\frac{1}{4\sqrt 2}\left[R^{\rmO^+_4}_{\bfT_1,1}+R^{\rmO^+_4}_{\bfT_3,1}-2R^{\rmO^+_4}_{\bfT_4,1}\right]
=\frac{1}{\sqrt 2}\left[\chi^+_{q^2}+\chi^+_{q^2}\cdot\sgn_{\rmO^+_4}\right], \\
R^{\rmO^-_4}_{\binom{0}{2}}
&=\frac{1}{2\sqrt 2}\left[-R^{\rmO^-_4}_{\bfT_2,1}-R^{\rmO^-_4}_{\bfT_5,1}\right]
=\frac{1}{\sqrt 2}\left[-{\bf 1}_{\rmO^-_4}-\sgn_{\rmO^-_4}\right], \\
R^{\rmO^-_4}_{\binom{1,0}{2,1}}
&=\frac{1}{2\sqrt 2}\left[-R^{\rmO^-_4}_{\bfT_2,1}+R^{\rmO^-_4}_{\bfT_5,1}\right]
=\frac{1}{\sqrt 2}\left[-\chi^-_{q^2}-\chi^-_{q^2}\cdot\sgn_{\rmO^-_4}\right].
\end{align*}
From (\ref{0205}), we know that
\begin{align*}
\textstyle
\bigl\langle\binom{2}{0},\binom{2}{0}\bigr\rangle
=\bigl\langle\binom{0}{2},\binom{2}{0}\bigr\rangle
=\bigl\langle\binom{2,1}{1,0},\binom{2,1}{1,0}\bigr\rangle
=\bigl\langle\binom{1,0}{2,1},\binom{2,1}{1,0}\bigr\rangle
&\equiv 0\pmod 2, \\
\textstyle
\bigl\langle\binom{-}{2,0},\binom{0}{2}\bigr\rangle
=\bigl\langle\binom{2,0}{-},\binom{0}{2}\bigr\rangle
=\bigl\langle\binom{1}{2,1,0},\binom{1,0}{2,1}\bigr\rangle
=\bigl\langle\binom{2,1,0}{1},\binom{1,0}{2,1}\bigr\rangle
&\equiv 1\pmod 2.
\end{align*}
By (\ref{0332}), we have
\begin{align*}
({\bf 1}_{\rmO^+_4})^\sharp=(\sgn_{\rmO^+_4})^\sharp
&=\frac{1}{\sqrt 2}R^{\rmO^+_4}_{\binom{2}{0}}
=\frac{1}{2}({\bf 1}_{\rmO^+_4}+\sgn_{\rmO^+_4}), \\
(\chi^+_{q^2})^\sharp=(\chi^+_{q^2}\cdot\sgn_{\rmO^+_4})^\sharp
&=\frac{1}{\sqrt 2}R^{\rmO^+_4}_{\binom{2,1}{1,0}}
=\frac{1}{2}(\chi^+_{q^2}+\chi^+_{q^2}\cdot\sgn_{\rmO^+_4}), \\
\chi^+_{2q} &=R^{\rmO^+_4}_{\binom{1}{1}}, \\
({\bf 1}_{\rmO^-_4})^\sharp=(\sgn_{\rmO^-_4})^\sharp
&=-\frac{1}{\sqrt 2}R^{\rmO^-_4}_{\binom{0}{2}}
=\frac{1}{2}({\bf 1}_{\rmO^-_4}+\sgn_{\rmO^-_4}), \\
(\chi^-_{q^2})^\sharp=(\chi^-_{q^2}\cdot\sgn_{\rmO^-_4})^\sharp
&=-\frac{1}{\sqrt 2}R^{\rmO^-_4}_{\binom{1,0}{2,1}}
=\frac{1}{2}(\chi^-_{q^2}+\chi^-_{q^2}\cdot\sgn_{\rmO^-_4}).
\end{align*}
Therefore any bijection $\cals_{\rmO^+_4}\rightarrow\cale(\rmO^+_4,1)$ such that
\begin{align*}
\bigl\{\textstyle\binom{2}{0},\binom{0}{2}\bigr\} &\mapsto\{{\bf 1}_{\rmO^+_4},\sgn_{\rmO^+_4}\}, \\
\textstyle\binom{1}{1} &\mapsto \chi^+_{2q}, \\
\bigl\{\textstyle\binom{2,1}{1,0},\binom{1,0}{2,1}\bigr\} &\mapsto\{\chi^+_{q^2},\chi^+_{q^2}\cdot\sgn_{\rmO^+_4}\}
\end{align*}
is a Lusztig parametrization for $\rmO^+_4$,
and any bijection $\cals_{\rmO^-_4}\rightarrow\cale(\rmO^-_4,1)$ such that
\begin{align*}
\bigl\{\textstyle\binom{-}{2,0},\binom{2,0}{-}\bigr\} &\mapsto\{{\bf 1}_{\rmO^-_4},\sgn_{\rmO^-_4}\}, \\
\bigl\{\textstyle\binom{1}{2,1,0},\binom{2,1,0}{1}\bigr\} &\mapsto\{\chi^-_{q^2},\chi^-_{q^2}\cdot\sgn_{\rmO^-_4}\}
\end{align*}
is a Lusztig parametrization for $\rmO^-_4$.
\end{exam}

\section{Finite Theta Correspondence of Unipotent Characters}
In this section we want to purpose several conditions to enforce
the parametrization $\cals_\bfG\rightarrow\cale(\bfG,1)$ for $\bfG=\rmO^\epsilon_{2n}$ to be unique.

\subsection{Finite theta correspondence on unipotent characters}\label{0520}
For a nontrivial additive character $\psi$ of $\bfF_q$,
let $\omega^\psi_{\Sp(W)}$ denote the (character of the ) \emph{Weil representation} of
the finite symplectic group $\Sp(W)$ with respect to $\psi$.

Let $(\bfG,\bfG')$ be a \emph{reductive dual pair} of the form $(\Sp_{2n},\SO_{2n'+1})$ or
$(\rmO^\epsilon_{2n},\Sp_{2n'})$ where $\epsilon=+$ or $-$.
The restriction $\omega^\psi_{\bfG,\bfG'}$ of the Weil character to $G\times G'$ gives a decomposition
\begin{equation}\label{0404}
\omega_{\bfG,\bfG'}^\psi
=\sum_{\rho\in\cale(\bfG),\ \rho'\in\cale(\bfG')}m_{\rho,\rho'}\rho\otimes\rho'
\end{equation}
where the multiplicity $m_{\rho,\rho'}$ is $1$ or $0$.
Then we have a relation
\[
\Theta^\psi_{\bfG,\bfG'}=\{\,(\rho,\rho')\in\cale(\bfG)\times\cale(\bfG')\mid m_{\rho,\rho'}\neq 0\,\}
\]
between $\cale(\bfG)$ and $\cale(\bfG')$ called the \emph{finite theta correspondence}
(or \emph{Howe duality}) for the dual pair $(\bfG,\bfG')$.
We say that an irreducible character $\rho\in\cale(\bfG)$ \emph{occurs in} $\Theta^\psi_{\bfG,\bfG'}$
if there exists $\rho'\in\cale(\bfG')$ such that $(\rho,\rho')\in\Theta^\psi_{\bfG,\bfG'}$.

For a symplectic space $V$ over $\bfF_q$, we have the \emph{symplectic similitude group}
\[
\GSp(V)=\{\,g\in\GL(V)\mid\langle gv,gw\rangle=k_g\langle v,w\rangle\text{ for some }k_g\in\bfF_q^\times
\text{ and any }v,w\in V\,\}.
\]
Note that $\GSp(V)$ normalizes the symplectic group $\Sp(V)$.
Choose an element $h\in\GSp_{2n}(q)$ such that $k_h$ is a non-square element in $\bfF_q^\times$.
For $\rho\in\cale(\Sp_{2n})$, we define the \emph{conjugate character} $\rho^c\in\cale(\Sp_{2n})$ by
$\rho^c(g)=\rho(hgh^{-1})$ for any $g\in\Sp_{2n}(q)$.

\begin{lem}\label{0405}
Suppose that $(\bfG,\bfG')=(\Sp_{2n},\SO_{2n'+1})$.
If $(\rho,\rho')\in\Theta^\psi_{\bfG,\bfG'}$,
then $(\rho^c,\rho')\in\Theta_{\bfG,\bfG'}^{\psi_a}$ where $\psi_a$ is another additive character
of\/ $\bfF_q$ given by $\psi_a(x):=\psi(ax)$ and $a$ is a non-square element in $\bfF_q^\times$.
\end{lem}
\begin{proof}
Suppose that $G=\Sp(V)$ and $G'=\SO(V')$ for a $2n$-dimensional symplectic space $V$ and a $(2n+1)$-dimensional
orthogonal space $V'$, and write $\omega^\psi_{\bfG,\bfG'}$ as in (\ref{0404}).
Choose $h\in\GSp(V)$ such that $k_h$ is a non-square element in $\bfF_q^\times$,
and let $\widetilde h=\iota(h,1)\in\GSp(V\otimes V')$ where
\[
\iota\colon\GSp(V)\times\SO(V')\rightarrow\GSp(V\otimes V').
\]
Now clearly $k_{\widetilde h}=k_h$,
and then by \cite{Szechtman} proposition 11 we have
$\omega^\psi_{\Sp(V\otimes V')}\circ{\rm Ad}_{\widetilde h}=\omega^{\psi_a}_{\Sp(V\otimes V')}$.
Therefore,
\[
\omega^{\psi_a}_{\bfG,\bfG'}
=\sum_{\rho\in\cale(\bfG),\ \rho'\in\cale(\bfG')}m_{\rho,\rho'}(\rho\circ{\rm Ad}_h)\otimes\rho'
=\sum_{\rho\in\cale(\bfG),\ \rho'\in\cale(\bfG')}m_{\rho,\rho'}\rho^c\otimes\rho'.
\]
Thus the lemma is proved.
\end{proof}

For $\rho\in\cale(\rmO^\epsilon_{2n})$,
we can define the \emph{conjugate character} $\rho^c\in\cale(\rmO^\epsilon_{2n})$ as
we did for a symplectic group given above.

\begin{lem}\label{0522}
Suppose that $(\bfG,\bfG')=(\rmO^\epsilon_{2n},\Sp_{2n'})$.
If $(\rho,\rho')\in\Theta^\psi_{\bfG,\bfG'}$,
then $(\rho^c,\rho'^c)\in\Theta_{\bfG,\bfG'}^\psi$,
and $(\rho^c,\rho'),(\rho,\rho'^c)\in\Theta_{\bfG,\bfG'}^{\psi_a}$
where $a$ is a non-square element in $\bfF_q^\times$.
\end{lem}
\begin{proof}
Suppose that $G=\rmO(V)$ and $G'=\Sp(V')$ for a $2n$-dimensional orthogonal space $V$ and
a $2n'$-dimensional symplectic space $V'$,
and write $\omega^\psi_{\bfG,\bfG'}$ as in (\ref{0404}).
Choose $h\in\GO(V)$ (the orthogonal similitude group) and $h'\in\GSp(V')$ such that both $k_h,k_{h'}$ are
non-square elements in $\bfF_q^\times$ and let $\widetilde h=\iota(h,h')\in\GSp(V\otimes V')$
where
\[
\iota\colon\GO(V)\times\GSp(V')\rightarrow\GSp(V\otimes V').
\]
Now $k_{\widetilde h}=k_hk_{h'}$ becomes a square element in $\bfF_q^\times$ and therefore
\begin{align*}
\omega^{\psi}_{\bfG,\bfG'}
=\omega^{\psi}_{\bfG,\bfG'}\circ\Ad_{\widetilde h}
&= \sum_{\rho\in\cale(\bfG),\ \rho'\in\cale(\bfG')}m_{\rho,\rho'}(\rho\circ{\rm Ad}_h)\otimes(\rho'\circ\Ad_{h'}) \\
&= \sum_{\rho\in\cale(\bfG),\ \rho'\in\cale(\bfG')}m_{\rho,\rho'}\rho^c\otimes\rho'^c.
\end{align*}
So we have shown that $(\rho,\rho')\in\Theta^\psi_{\bfG,\bfG'}$ implies that
$(\rho^c,\rho'^c)\in\Theta_{\bfG,\bfG'}^\psi$.
The other assertions can be proved by an analogous argument in the proof of Lemma~\ref{0405}.
\end{proof}

Let $\bfG'_{n'}$ denote $\SO_{2n'+1}$, $\Sp_{2n'}$, or $\rmO^\epsilon_{2n'}$.
For $\rho\in\cale(\bfG)$, it is well known that if $\rho$ occurs in $\Theta^\psi_{\bfG,\bfG'_{n'}}$,
then it also occurs in $\Theta^\psi_{\bfG,\bfG'_{n''}}$ for any $n''\geq n'$.
We say that $\rho$ \emph{first occurs in} $\Theta^\psi_{\bfG,\bfG'_{n'}}$ if it occurs in
$\Theta^\psi_{\bfG,\bfG'_{n'}}$ and does not occur in $\Theta^\psi_{\bfG,\bfG'_{n'-1}}$.

\subsection{Finite theta correspondence on unipotent characters}\label{0403}
If $(\bfG,\bfG')=(\rmO^\epsilon_{2n},\Sp_{2n'})$,
then the unipotent characters are preserved by $\Theta_{\bfG,\bfG'}^\psi$ (\cf.\ \cite{adams-moy} theorem 3.5),
i.e., we can write
\begin{align*}
\omega^\psi_{\bfG,\bfG',1}
&=\sum_{\rho\in\cale(\bfG,1),\ \rho'\in\cale(\bfG',1)}m_{\rho,\rho'}\rho\otimes\rho' \\
\Theta^\psi_{\bfG,\bfG',1} &=\Theta^\psi_{\bfG,\bfG'}\cap(\cale(\bfG,1)\times\cale(\bfG',1))
\end{align*}
where $\omega^\psi_{\bfG,\bfG',1}$ denotes the unipotent part of $\omega^\psi_{\bfG,\bfG'}$.
For the theta correspondence on unipotent characters,
the following are well-known
(\cf.~\cite{adams-moy} theorem 5.2):
\begin{itemize}
\item ${\bf 1}_{\rmO^+_2}$ first occurs in the correspondence for the pair $(\rmO^+_2,\Sp_0)$,

\item $\sgn_{\rmO^+_2}$ first occurs in the correspondence for the pair $(\rmO^+_2,\Sp_2)$,

\item $\zeta_k^\rmI$ first occurs in the correspondence for the pair $(\rmO^{\epsilon_k}_{2k^2},\Sp_{2k(k-1)})$,

\item $\zeta_k^{\rm II}$ first occurs in the correspondence for the pair $(\rmO^{\epsilon_k}_{2k^2},\Sp_{2k(k+1)})$
\end{itemize}
where $\epsilon_k=(-1)^k$, and $\zeta_k^\rmI, \zeta_k^{\rm II}$ are the unipotent cuspidal
characters of $\rmO^{\epsilon_k}_{2k^2}(q)$ given in Subsection~\ref{0327}.

Now we recall some results on $\Theta_{\bfG,\bfG',1}^\psi$ from \cite{pan-finite-unipotent}.
For any two partitions $\lambda=[\lambda_1,\lambda_2,\ldots]$ (with $\lambda_1\geq\lambda_2\geq\cdots$),
$\mu=[\mu_1,\mu_2,\ldots]$ (with $\mu_1\geq\mu_2\geq\cdots$), we define a relation
\begin{equation}\label{0202}
\lambda\preccurlyeq\mu\qquad\text{if }\mu_1\geq\lambda_1\geq\mu_2\geq\lambda_2\geq\mu_3\geq\lambda_3\geq\cdots.
\end{equation}
And then we define a relation $\calb_{\bfG,\bfG'}$ between $\cals_\bfG$ and $\cals_{\bfG'}$ for
$(\bfG,\bfG')=(\rmO^\epsilon_{2n},\Sp_{2n'})$ as follows:
\begin{itemize}
\item If $\epsilon=+$,
let $\calb_{\bfG,\bfG'}$ be the set consisting of pairs $(\Lambda,\Lambda')$ of symbols such that
\begin{itemize}
\item $\Upsilon(\Lambda)_*\preccurlyeq\Upsilon(\Lambda')^*$ and
    $\Upsilon(\Lambda')_*\preccurlyeq\Upsilon(\Lambda)^*$,

\item ${\rm def}(\Lambda')=-{\rm def}(\Lambda)+1$;
\end{itemize}

\item if $\epsilon=-$,
let $\calb_{\bfG,\bfG'}$ be the set consisting of pairs $(\Lambda,\Lambda')$ of symbols such that
\begin{itemize}
\item $\Upsilon(\Lambda)^*\preccurlyeq\Upsilon(\Lambda')_*$ and
    $\Upsilon(\Lambda')^*\preccurlyeq\Upsilon(\Lambda)_*$,

\item ${\rm def}(\Lambda')=-{\rm def}(\Lambda)-1$.
\end{itemize}
\end{itemize}

We say that a symbol $\Lambda\in\cals_\bfG$ \emph{occurs in} $\calb_{\bfG,\bfG'}$
if there is $\Lambda'\in\cals_{\bfG'}$ such that $(\Lambda,\Lambda')\in\calb_{\bfG,\bfG'}$.
For $\Lambda\in\cals_\bfG$, it is not difficult to see that if $\Lambda$ occurs in
$\calb_{\bfG,\bfG'_{n'}}$,
then it also occurs in $\calb_{\bfG,\bfG'_{n''}}$ for any $n''\geq n'$.
We say that $\Lambda$ \emph{first occurs in} $\calb_{\bfG,\bfG'_{n'}}$ if it occurs in
$\calb_{\bfG,\bfG'_{n'}}$ and does not occur in $\calb_{\bfG,\bfG'_{n'-1}}$.

\begin{exam}\label{0519}
\begin{enumerate}
\item We have $\binom{1}{0},\binom{0}{1}\in\cals_{\rmO^+_2}$,
$\binom{0}{-}\in\cals_{\Sp_0}$, and $\binom{1}{-}\in\cals_{\Sp_2}$.
Now $\Upsilon(\binom{1}{0})=\sqbinom{1}{0}$, $\Upsilon(\binom{0}{1})=\sqbinom{0}{1}$,
$\Upsilon(\binom{0}{-})=\sqbinom{0}{0}$, and $\Upsilon(\binom{1}{-})=\sqbinom{1}{0}$,
and so $\binom{1}{0}$ first occurs in $\calb_{\rmO^+_2,\Sp_0}$ and
$\binom{0}{1}$ first occurs in $\calb_{\rmO^+_2,\Sp_2}$.

\item Suppose that $k$ is even, and let $\Lambda_k^{\rm I}=\binom{2k-1,2k-2,\ldots,0}{-}$.
Then $\Lambda_k^{\rm I}\in\cals_{\rmO^+_{2k^2}}$, ${\rm def}(\Lambda_k^{\rm I})=2k$ and
$\Upsilon(\Lambda_k^{\rm I})=\sqbinom{0}{0}$.
Therefore, if $\Lambda_k^{\rm I}$ occurs in $\calb_{\rmO^+_{2k^2},\Sp_{2n'}}$,
then $\cals_{\Sp_{2n'}}$ contains a symbol of defect $-2k+1$.
By (\ref{0204}), we have $n'\geq k(k-1)$, i.e., $\Lambda_k^{\rm I}$ first occurs in
$\calb_{\rmO^+_{2k^2},\Sp_{2(k^2-k)}}$.

\item Suppose that $k$ is odd, and let $\Lambda_k^{\rm I}=\binom{-}{2k-1,2k-2,\ldots,0}$.
Then $\Lambda_k^{\rm I}\in\cals_{\rmO^-_{2k^2}}$, ${\rm def}(\Lambda_k^{\rm I})=-2k$ and
$\Upsilon(\Lambda_k^{\rm I})=\sqbinom{0}{0}$.
This means that if $\Lambda_k^{\rm I}$ occurs in $\calb_{\rmO^-_{2k^2},\Sp_{2n'}}$,
then $\cals_{\Sp_{2n'}}$ contains a symbol of defect $2k-1$.
By (\ref{0204}), we have $n'\geq k(k-1)$, i.e., $\Lambda_k^{\rm I}$ first occurs in
$\calb_{\rmO^-_{2k^2},\Sp_{2(k^2-k)}}$.
\end{enumerate}
\end{exam}

The following proposition is from \cite{pan-finite-unipotent} corollary 5.36:
\begin{prop}
Let $(\bfG,\bfG')=(\rmO^\epsilon_{2n},\Sp_{2n'})$ where $\epsilon=+$ or $-$.
Let $\call_1\colon\cals_\bfG\rightarrow\cale(\bfG,1)$ by $\Lambda\mapsto\rho_\Lambda$ and
$\call'_1\colon\cals_{\bfG'}\rightarrow\cale(\bfG',1)$ by $\Lambda'\mapsto\rho_{\Lambda'}$ be
any Lusztig parametrizations for $\bfG,\bfG'$ respectively.
Then $(\rho_\Lambda,\rho_{\Lambda'})$ or $(\rho_{\Lambda^\rmt},\rho_{\Lambda'})$ occurs in
$\Theta_{\bfG,\bfG',1}^\psi$ if and only if  $(\Lambda,\Lambda')$ or $(\Lambda^\rmt,\Lambda')$ occurs in
$\calb_{\bfG,\bfG'}$.
\end{prop}

\begin{rem}
Note that the parametrization $\call'_1$ is unique by Corollary~\ref{0311} but $\call_1$ is not.
So we want to enforce more conditions on $\call_1$ so that $\call_1$ is unique and
eliminate the ambiguity in the above proposition.
\end{rem}

\subsection{On the uniqueness of Lusztig parametrization for even orthogonal groups}\label{0521}
Now we want to enforce extra conditions on the Lusztig parametrization
$\call_1\colon\cals_{\rmO^\epsilon_{2n}}\rightarrow\cale(\rmO^\epsilon_{2n},1)$
to make it be uniquely determined.

\begin{enumerate}
\item[(I)] We require that $\call_1$ by $\Lambda\mapsto\rho_\Lambda$ is compatible with the parabolic induction
on unipotent characters, i.e.,
we require that $\Omega(\rho_\Lambda)=\{\,\rho_{\Lambda'}\mid\Lambda'\in\Omega(\Lambda)\,\}$ where
$\Omega(\rho_\Lambda)$ and $\Omega(\Lambda)$ are defined as in Subsection~\ref{0325}.

\item[(II)] We require that
\begin{itemize}
\item for $k\geq 1$, $\call_1(\Lambda_k^{\rm I})=\zeta_k^{\rm I}$
and $\call_1(\Lambda_k^{\rm II})=\zeta_k^{\rm II}$, i.e.,
$\zeta_k^\rmI=\rho_{\Lambda_k^\rmI}$ and $\zeta_k^{\rm II}=\rho_{\Lambda_k^{\rm II}}$
where
\begin{align}\label{0402}
\begin{split}
\Lambda_k^\rmI &=\begin{cases}
\binom{2k-1,2k-2,\ldots,1,0}{-}, & \text{if $k$ is even};\\
\binom{-}{2k-1,2k-2,\ldots,1,0}, & \text{if $k$ is odd},
\end{cases} \\
\Lambda_k^{\rm II} &=(\Lambda_k^\rmI)^\rmt
\end{split}
\end{align}
and $\xi_k^\rmI,\xi_k^{\rm II}$ the two cuspidal unipotent characters of $\rmO^{\epsilon_k}_{2k^2}$
given in the previous subsection and $\epsilon_k=(-1)^k$;

\item $\call_1(\binom{1}{0})={\bf 1}_{\rmO^+_2}$ and
$\call_1(\binom{0}{1})=\sgn_{\rmO^+_2}$,
i.e., $\rho_{\binom{1}{0}}={\bf 1}_{\rmO^+_2}$ and $\rho_{\binom{0}{1}}=\sgn_{\rmO^+_2}$.
\end{itemize}
\end{enumerate}
Note that in addition to the specification of $\call_1$ on cuspidal symbols,
due to Remark~\ref{0336} we also need to assign the image of $\call_1$ at $\binom{1}{0}$
or $\binom{0}{1}$.

By Corollary~\ref{0407}, a Lusztig parametrization
$\call_1\colon\cals_{\rmO^\epsilon_{2n}}\rightarrow\cale(\rmO^\epsilon_{2n},1)$
satisfying both (I) and (II) clearly exists.

\begin{prop}\label{0305}
There is a unique bijective parametrization $\cals_{\rmO^\epsilon_{2n}}\rightarrow\cale(\rmO^\epsilon_{2n},1)$
where $\epsilon=+$ or $-$ satisfying (\ref{0332}), and (I), (II) above.
\end{prop}
\begin{proof}
The existence of such a bijection $\call_1$ is obvious,
so now we consider the uniqueness.
Let $\call_1,\call'_1\colon\cals_{\rmO^\epsilon_{2n}}\rightarrow\cale(\rmO^\epsilon_{2n},1)$ be two Lusztig
parametrizations of unipotent characters for $\rmO^\epsilon_{2n}$.
Moreover, suppose that $\call_1,\call'_1$  both satisfy (I) and (II) above.
For $\Lambda\in\cals_{\rmO^\epsilon_{2n}}$,
we know that $\call_1(\Lambda)^\sharp=\call_1'(\Lambda)^\sharp$ by Lemma~\ref{0337},
and hence by Corollary~\ref{0314} either
\[
\call_1(\Lambda)=\call_1'(\Lambda)\quad\text{ or }\quad \call_1(\Lambda)=\call_1'(\Lambda)\cdot\sgn_{\rmO^\epsilon_{2n}}=\call'_1(\Lambda^\rmt),
\]
i.e., if $\call_1(\Lambda)=\call'_1(\Lambda')$, then either $\Lambda'=\Lambda$ or $\Lambda'=\Lambda^\rmt$.
Now we suppose that $\call_1(\Lambda)=\call_1'(\Lambda')$ and consider the following three cases:
\begin{enumerate}
\item Suppose that $\Lambda$ is degenerate, i.e., $\Lambda=\Lambda^\rmt$.
Then $\call_1(\Lambda)=\call_1'(\Lambda')$ implies that $\call_1(\Lambda)=\call_1'(\Lambda)$ immediately.

\item Suppose that ${\rm def}(\Lambda)\neq 0$.
Suppose that the unipotent character $\call_1(\Lambda)=\call'_1(\Lambda')$
where $\Lambda'=\Lambda$ or $\Lambda^\rmt$  is in the Harish-Chandra series initiated
by some unipotent cuspidal character $\zeta$.
Because ${\rm def}(\Lambda)\neq 0$, we have $\zeta\neq\zeta\cdot\sgn$.
By the requirement in (II), we have $\call_1(\Lambda_0)=\zeta=\call'_1(\Lambda_0)$ for some
cuspidal symbol $\Lambda_0$ such that ${\rm def}(\Lambda_0)\neq 0$.
By (I), we must have ${\rm def}(\Lambda)={\rm def}(\Lambda')$, and then we conclude that $\Lambda'=\Lambda$, i.e.,
$\call_1(\Lambda)=\call'_1(\Lambda)$.

\item Suppose that $\Lambda$ is non-degenerate and ${\rm def}(\Lambda)=0$, i.e, $\Lambda\in\cals_{\rmO^+_{2n}}$
for some $n$.
Now we are going to prove this case by induction on $n$.
For $n=1$, the equality $\call_1(\Lambda)=\call'_1(\Lambda)$ is enforced by (II) above.
Now suppose that $n\geq 2$.
Because now $\Lambda^\rmt\neq\Lambda$, by Lemma~\ref{0328}, there exists
$\Lambda_1\in\cals_{\rmO^+_{2(n-1)}}$ such that $\Lambda\in\Omega(\Lambda_1)$
and $\Lambda^\rmt\not\in\Omega(\Lambda_1)$.
By (I) and the induction hypothesis, we have
\[
\call_1(\Lambda)\in\Omega(\call_1(\Lambda_1))=\Omega(\call_1'(\Lambda_1))\not\ni \call'_1(\Lambda^\rmt).
\]
Now $\call'_1(\Lambda^\rmt)\neq\call_1(\Lambda)$ implies that $\call'_1(\Lambda)=\call_1(\Lambda)$.
\end{enumerate}
Hence the proposition is proved.
\end{proof}

\begin{cor}
Let $\bfG=\rmO^\epsilon_{2n}$,
and let $\Lambda\mapsto\rho_\Lambda$ be the Lusztig parametrization in Proposition~\ref{0305}.
Then
\begin{enumerate}
\item[(i)] ${\bf 1}_{\rmO^-_{2n}}=\rho_{\binom{-}{n,0}}$ and $\sgn_{\rmO^-_{2n}}=\rho_{\binom{n,0}{-}}$;

\item[(ii)] ${\bf 1}_{\rmO^+_{2n}}=\rho_{\binom{n}{0}}$ and $\sgn_{\rmO^+_{2n}}=\rho_{\binom{0}{n}}$
\end{enumerate}
\end{cor}
\begin{proof}
Let $\Lambda\mapsto\rho_\Lambda$ be the parametrization for $\rmO^\epsilon_{2n}$
satisfying (\ref{0332}) and (I), (II) above.
We know that ${\bf 1}_{\rmO^-_2}$ (resp.~$\sgn_{\rmO^-_2}$) first occurs in the correspondence for the pair
$(\rmO^-_2,\Sp_0)$ (resp.~$(\rmO^-_2,\Sp_4)$).
Therefore by the requirement in (II) above,
we have ${\bf1}_{\rmO^-_2}=\zeta^\rmI_1=\rho_{\binom{-}{1,0}}$
(resp.\ $\sgn_{\rmO^-_2}=\zeta^{\rm II}_1=\rho_{\binom{1,0}{-}}$).
Write ${\bf 1}_{\rmO^-_{2n}}=\rho_\Lambda$ for some $\Lambda\in\cals_{\rmO^-_{2n}}$.
By Corollary~\ref{0310} we know that $\Lambda$ is either $\binom{-}{n,0}$ or $\binom{n,0}{-}$.
Because ${\bf 1}_{\rmO^-_{2n}}$ is an irreducible constituent of $R^{\rmO^-_{2n}}_{\rmO^-_2\times\GL_1^{n-1}}(\rho_{\binom{-}{1,0}}\otimes{\bf1})$.
By (I) above, we must have ${\rm def}(\Lambda)={\rm def}(\binom{-}{1,0})=-2$,
and so we conclude that ${\bf 1}_{\rmO^-_{2n}}=\rho_{\binom{-}{n,0}}$ and hence
$\sgn_{\rmO^-_{2n}}=\rho_{\binom{n,0}{-}}$.

Now we are going to prove case (ii) by induction on $n$.
By (II) above, we have $\rho_{\binom{1}{0}}={\bf 1}_{\rmO^+_2}$ and $\rho_{\binom{0}{1}}=\sgn_{\rmO^+_2}$.
Now by the induction hypothesis, for $n\geq 2$, we assume that
${\bf 1}_{\rmO^+_{2(n-1)}}=\rho_{\binom{n-1}{0}}$ and $\sgn_{\rmO^+_{2(n-1)}}=\rho_{\binom{0}{n-1}}$.
Suppose that ${\bf 1}_{\rmO^+_{2n}}=\rho_\Lambda$ for some $\Lambda\in\cals_{\rmO^+_{2n}}$.
Then we know that either $\Lambda=\binom{n}{0}$ or $\Lambda=\binom{0}{n}$ by Corollary~\ref{0310}.
Because ${\bf 1}_{\rmO^+_{2n}}\in\Omega({\bf 1}_{\rmO^+_{2(n-1)}})$,
by (II) we see that $\Lambda\in\Omega(\binom{n-1}{0})$ and therefore $\Lambda$ must be $\binom{n}{0}$,
i.e., we conclude that ${\bf 1}_{\rmO^+_{2n}}=\rho_{\binom{n}{0}}$ and $\sgn_{\rmO^+_{2n}}=\rho_{\binom{0}{n}}$.
\end{proof}

\begin{exam}
Keep the notation in Example~\ref{0304}.
Suppose that $\rho_1,\rho_2\in\cale(\rmO^-_4,1)$ are the two irreducible characters of degree $q^2$
satisfying
\[
R_{\rmO^-_2\times\GL_1}^{\rmO^-_4}({\bf1}_{\rmO^-_2}\otimes{\bf1})={\bf 1}_{\rmO^-_4}+\rho_1,\qquad
R_{\rmO^-_2\times\GL_1}^{\rmO^-_4}(\sgn_{\rmO^-_2}\otimes{\bf1})=\sgn_{\rmO^-_4}+\rho_2.
\]
Now (II) implies that ${\bf 1}_{\rmO^-_2}=\rho_{\binom{-}{1,0}}$ and $\sgn_{\rmO^-_2}=\rho_{\binom{1,0}{-}}$.
Therefore the parametrization $\cals_{\rmO^-_4}\rightarrow\cale(\rmO^-_4,1)$ satisfying (I), (II) above must be
\[
\rho_{\binom{-}{2,0}}={\bf1}_{\rmO^-_4},\qquad
\rho_{\binom{2,0}{-}}=\sgn_{\rmO^-_4},\qquad
\rho_{\binom{1}{2,1,0}}=\rho_1,\qquad
\rho_{\binom{2,1,0}{1}}=\rho_2.
\]
\end{exam}

The following proposition which justifies our choice of $\call_1$ in Proposition~\ref{0305}
is from \cite{pan-finite-unipotent} theorem 1.8.

\begin{prop}
Let $(\bfG,\bfG')=(\rmO^\epsilon_{2n'},\Sp_{2n})$ where $\epsilon=+$ or $-$.
Let $\call_1\colon\cals_\bfG\rightarrow\cale(\bfG,1)$ and $\call'_1\colon\cals_{\bfG'}\rightarrow\cale(\bfG',1)$
be the unique Lusztig parametrizations given in Proposition~\ref{0305} and Proposition~\ref{0308} respectively.
Then the diagram
\[
\begin{CD}
\cals_\bfG @> \calb_{\bfG,\bfG'} >> \cals_{\bfG'} \\
@V\call_1 VV @VV\call'_1 V \\
\cale(\bfG,1) @> \Theta^\psi_{\bfG,\bfG',1} >> \cale(\bfG',1) \\
\end{CD}
\]
commutes, i.e., $(\rho_\Lambda,\rho_{\Lambda'})$ occurs in $\Theta_{\bfG,\bfG',1}^\psi$
if and only if $(\Lambda,\Lambda')\in\calb_{\bfG,\bfG'}$.
\end{prop}
When both $\Lambda,\Lambda'$ are cuspidal,
the commutativity of the above diagram can be seen by the requirement (II)
(\cf.~Example~\ref{0519}).
For general $\Lambda,\Lambda'$, the commutativity follows from the fact that both
the correspondence $\Theta_{\bfG,\bfG'}^\psi$ and the parametrizations $\call_1,\call'_1$
are compatible with the parabolic induction.
Details of the proof can be found in \cite{pan-finite-unipotent}.

\begin{exam}
Let $\call_1\colon\cals_{\rmO^\epsilon_{2n}}\rightarrow\cale(\rmO^\epsilon_{2n},1)$ by
$\Lambda\mapsto\rho_\Lambda$ be the parametrization in Proposition~\ref{0305}.
Then by Corollary~\ref{0310} the two Steinberg characters of\/ $\rmO^\epsilon_{2n}$ are parametrized by
the symbols
\[
\begin{cases}
\binom{n,n-1,\ldots,1}{n-1,n-2,\ldots,0}, \binom{n-1,n-2,\ldots,0}{n,n-1,\ldots,1},
& \text{if $\epsilon=+$};\\
\binom{n,n-1,\ldots,0}{n-1,n-2,\ldots,1}, \binom{n-1,n-2,\ldots,1}{n,n-1,\ldots,0},
& \text{if $\epsilon=-$}.
\end{cases}
\]
\begin{enumerate}
\item If $\epsilon=+$,
then one of the Steinberg characters of $\rmO^+_{2n}$ first occurs in the correspondence
$\Theta^\psi_{\rmO^+_{2n},\Sp_{2(n-1)},1}$ and is paired with $\St_{\Sp_{2(n-1)}}$,
and this Steinberg character is parametrized by the symbol
$\binom{n,n-1,\ldots,1}{n-1,n-2,\ldots,0}$.

\item If $\epsilon=-$,
then one of the Steinberg characters of $\rmO^-_{2n}$ first occurs in the correspondence
$\Theta^\psi_{\rmO^-_{2n},\Sp_{2(n-1)},1}$ and is paired with $\St_{\Sp_{2(n-1)}}$,
and this Steinberg character is parametrized by the symbol
$\binom{n-1,n-2,\ldots,1}{n,n-1,\ldots,0}$.
\end{enumerate}
\end{exam}

\section{Lusztig Correspondence and Finite Theta Correspondence}

\subsection{Lusztig correspondences}
Let $\bfG$ be a classical group,
and let $s$ be a semisimple element in the connected component $(G^*)^0$ of $G^*$.
A rational maximal torus $\bfT^*$ in $\bfG^*$ contains $s$ if and only if it is a rational maximal
torus in $C_{\bfG^*}(s)$.
From \cite{carter-finite} theorem~7.3.4, it is known that
\[
\langle R^\bfG_{\bfT^*,s}, R^\bfG_{\bfT'^*,s}\rangle_\bfG
=\langle R^{C_{\bfG^*}(s)}_{\bfT^*,1},R^{C_{\bfG^*}(s)}_{\bfT'^*,1}\rangle_{C_{\bfG^*}(s)}
\]
for any rational maximal tori $\bfT^*,\bfT'^*$ of $\bfG^*$ containing $s$.
Then the mapping
\[
\epsilon_\bfG R^\bfG_{\bfT^*,s}\mapsto \epsilon_{C_{\bfG^*}(s)}R_{\bfT^*,1}^{C_{\bfG^*}(s)}
\]
for any $\bfT^*$ containing $s$ can be extended uniquely to an isometry from $\calv(\bfG,s)^\sharp$
onto $\calv(C_{\bfG^*},1)^\sharp$.
Now a Lusztig correspondence $\grL_s\colon\cale(\bfG,s)\rightarrow\cale(C_{\bfG^*}(s),1)$,
i.e., a bijective mapping satisfying (\ref{0201}), can be extended linearly to be an isometry,
still denoted by $\grL_s$, of inner product spaces
\begin{equation}\label{0502}
\grL_s\colon\calv(\bfG,s)\longrightarrow\calv(C_{\bfG^*}(s),1)
\end{equation}
whose restriction to $\calv(\bfG,s)^\sharp$ is uniquely determined.

Suppose that $\bfG=\prod_{i=1}^k\bfG_k$.
Then $\bfG^*=\prod_{i=1}^k\bfG^*_i$ where $\bfG_i^*$ is the dual group of $\bfG_i$.
If $s\in G^*$ semisimple, the we can write $s=(s_1,\ldots,s_k)$ where each $s_i\in G^*_i$ is semisimple,
and then $C_{\bfG^*}(s)=\prod_{i=1}^k C_{\bfG^*_i}(s_i)$.
Now a rational maximal torus $\bfT^*$ containing $s$ can be written as
$\bfT^*=\prod_{i=1}^k\bfT^*_i$ where $\bfT^*_i$ is a rational maximal torus in $\bfG^*_i$.
Therefore, we have $R^\bfG_{\bfT^*,s}=\bigotimes_{i=1}^k R^{\bfG_i}_{\bfT^*_i,s_i}$.

\begin{cor}
Then a bijection $\grL_s\colon\cale(\bfG,s)\rightarrow\cale(C_{\bfG^*}(s),1)$ is a Lusztig correspondence
if and only if $\grL_s=\prod_{i=1}^k\grL_{s_i}$
where each $\grL_{s_i}\colon\cale(\bfG_i,s_i)\rightarrow\cale(C_{\bfG^*_i}(s_i),1)$
is a Lusztig correspondence.
\end{cor}
\begin{proof}
This is obvious.
\end{proof}

For $s\in (G^*)^0$, we define
\begin{align}
\begin{split}
\bfG^{(0)}=\bfG^{(0)}(s)
&=\prod_{\langle\lambda\rangle\subset\{\lambda_1,\ldots,\lambda_n\},\ \lambda\neq\pm 1}\bfG_{[\lambda]}(s);\\
\bfG^{(1)}=\bfG^{(1)}(s)
&=\bfG_{[-1]}(s); \\
\bfG^{(2)}=\bfG^{(2)}(s)
&=\bfG_{[1]}(s)
\end{split}
\end{align}
where $\bfG_{[\lambda]}(s)$ is given in \cite{amr} subsection 1.B
(see also \cite{pan-Lusztig-correspondence} subsection 2.2).
We know that
\[
C_{\bfG^*}(s)\simeq\bfG^{(0)}\times\bfG^{(1)}\times\bfG^{(2)},
\]
and $\bfG^{(0)}$ is a product of general linear groups or unitary group, and
\begin{equation}\label{0510}
(\bfG^{(1)},\bfG^{(2)})=\begin{cases}
(\Sp_{2n^{(1)}},\Sp_{2n^{(2)}}), & \text{if $\bfG=\SO_{2n+1}$};\\
(\rmO^{\epsilon^{(1)}}_{2n^{(1)}},\SO_{2n^{(2)}+1}), & \text{if $\bfG=\Sp_{2n}$};\\
(\rmO^{\epsilon^{(1)}}_{2n^{(1)}},\rmO^{\epsilon^{(2)}}_{2n^{(2)}}), & \text{if $\bfG=\rmO^\epsilon_{2n}$}
\end{cases}
\end{equation}
for some non-negative integers $n^{(1)},n^{(2)}$ depending on $s$,
and some $\epsilon^{(1)},\epsilon^{(2)}$.
Note that if $\bfG=\rmO^\epsilon_{2n}$, then $\epsilon^{(1)},\epsilon^{(2)}$ also depend on $s$ (and $\epsilon$),
if $\bfG=\Sp_{2n}$, then $\epsilon^{(1)}$ can be $+$ or $-$ for each $s$ such that $n^{(1)}\geq 1$.
The element $s$ can be written as
\begin{equation}\label{0516}
s=s^{(0)}\times s^{(1)}\times s^{(2)}
\end{equation}
where $s^{(1)}$ (resp.~$s^{(2)}$) is the part whose eigenvalues are all equal to $-1$ (resp.~$1$),
and $s^{(0)}$ is the part whose eigenvalues do not contain $1$ or $-1$.
In particular, $s^{(j)}$ is in the center of $\bfG^{(j)}$.
Then a Lusztig correspondence
\begin{equation}\label{0511}
\grL_s\colon\cale(\bfG,s)\rightarrow\cale(\bfG^{(0)}\times\bfG^{(1)}\times\bfG^{(2)},1)
\end{equation}
can be written as
\begin{equation}\label{0505}
\grL_s(\rho)=\rho^{(0)}\otimes\rho^{(1)}\otimes\rho^{(2)}
\end{equation}
where $\rho^{(j)}\in\cale(\bfG^{(j)},1)$ for $j=0,1,2$.
It is known that $\rho$ is cuspidal if and only if $\bfG^{(0)}$ is a product of unitary groups
(i.e., no general linear groups) and each $\rho^{(j)}$ is cuspidal.

Now $\cals_{C_{\bfG^*}(s)}^\sharp=\cals_{\bfG^{(0)}}^\sharp\times \cals_{\bfG^{(1)}}^\sharp\times \cals_{\bfG^{(2)}}^\sharp$ and
\[
R_\Sigma^{C_{\bfG^*}(s)}
=R^{\bfG^{(0)}}_{x}\otimes R^{\bfG^{(1)}}_{\Sigma^{(1)}}\otimes R^{\bfG^{(2)}}_{\Sigma^{(2)}}
\in\calv(C_{\bfG^*}(s),1)^\sharp
\]
for $\Sigma=(x,\Sigma^{(1)},\Sigma^{(2)})\in\cals^\sharp_{C_{\bfG^*}(s)}$.
We define
\begin{equation}\label{0614}
R_\Sigma^\bfG=\grL_s^{-1}(R_\Sigma^{C_{\bfG^*}(s)})\in\calv(\bfG,s)^\sharp.
\end{equation}
Because $\{\,R_\Sigma^{C_{\bfG^*}(s)}\mid\Sigma\in\cals^\sharp_{C_{\bfG^*}(s)}\,\}$ is an orthonormal basis
for $\calv(C_{\bfG^*}(s),1)^\sharp$ and
\[
\grL_s\colon\calv(\bfG,s)\rightarrow\calv(C_{\bfG^*}(s),1)
\]
is an isometry which maps $\calv(\bfG,s)^\sharp$ onto $\calv(C_{\bfG^*}(s),1)^\sharp$,
we see that $\{\,R_\Sigma^\bfG\mid\Sigma\in\cals^\sharp_{C_{\bfG^*}(s)}\,\}$ forms an orthonormal basis for
the space $\calv(\bfG,s)^\sharp$.
For $\rho\in\cale(\bfG,s)$, we have
\begin{align}\label{0514}
\begin{split}
\rho^\sharp &=\sum_{\Sigma\in \cals_{C_{\bfG^*}(s)}^\sharp}\langle\rho,R_\Sigma^\bfG\rangle_\bfG R_\Sigma^\bfG, \\
\grL_s(\rho)^\sharp &=\sum_{\Sigma\in \cals_{C_{\bfG^*}(s)}^\sharp}
\langle\grL_s(\rho),R_\Sigma^{C_{\bfG^*}(s)}\rangle_{C_{\bfG^*}(s)} R_\Sigma^{C_{\bfG^*}(s)}.
\end{split}
\end{align}
Therefore we have $\grL_s(\rho^\sharp)=\grL_s(\rho)^\sharp$ for any $\rho\in\cale(\bfG,s)$.

\begin{lem}\label{0504}
Let $\grL_s,\grL'_s\colon\cale(\bfG,s)\rightarrow\cale(C_{\bfG^*}(s),1)$ be two Lusztig correspondences,
and write $\grL_s(\rho)=\rho^{(0)}\otimes\rho^{(1)}\otimes\rho^{(2)}$,
$\grL'_s(\rho)=\rho'^{(0)}\otimes\rho'^{(1)}\otimes\rho'^{(2)}$.
Then
\[
\rho^{(0)}=\rho'^{(0)},\qquad (\rho^{(1)})^\sharp=(\rho'^{(1)})^\sharp,\qquad
(\rho^{(2)})^\sharp=(\rho'^{(2)})^\sharp.
\]
\end{lem}
\begin{proof}
Because $\grL_s(\rho^\sharp)=\grL_s(\rho)^\sharp$,
we have
\[
\grL_s(\rho^\sharp)=(\rho^{(0)})^\sharp\otimes(\rho^{(1)})^\sharp\otimes(\rho^{(2)})^\sharp,\qquad
\grL'_s(\rho^\sharp)=(\rho'^{(0)})^\sharp\otimes(\rho'^{(1)})^\sharp\otimes(\rho'^{(2)})^\sharp.
\]
Because the restrictions of $\grL_s$ and $\grL'_s$ to $\calv(\bfG,s)^\sharp$ are the same, i.e.,
$\grL_s(\rho^\sharp)=\grL_s'(\rho^\sharp)$,
we have
$(\rho^{(0)})^\sharp=(\rho'^{(0)})^\sharp$,
$(\rho^{(1)})^\sharp=(\rho'^{(1)})^\sharp$,
$(\rho^{(2)})^\sharp=(\rho'^{(2)})^\sharp$.
Now $\bfG^{(0)}$ is a product of general linear groups or unitary groups,
so we have $\rho^{(0)}=(\rho^{(0)})^\sharp=(\rho'^{(0)})^\sharp=\rho'^{(0)}$.
\end{proof}

\subsection{Lusztig correspondence and parabolic induction}\label{0615}
Let $\bfG_n=\SO_{2n+1}$, $\Sp_{2n}$ or $\rmO^\epsilon_{2n}$ where $\epsilon=+$ or $-1$.
The group $\bfG_n\times\GL_l$ is the Levi factor of a parabolic subgroup of $\bfG_{n+l}$.
Let $\sigma$ be an irreducible cuspidal character of $\GL_l$,
and so $\sigma\in\cale(\GL_l,t)$ for some semisimple element $t\in\GL_l(q)$.
For $\rho\in\cale(\bfG_n,s)$, an irreducible constituent of
$R^{\bfG_{n+l}}_{\bfG_n\times\GL_l}(\rho\otimes\sigma)$ is in $\cale(\bfG_{n+l},s')$
where $s'=(s,t)$ is regarded as an element in $(G_{n+l}^*)^0$.
Then we define a relation $\Omega_t\colon\cale(\bfG_n,s)\rightarrow\cale(\bfG_{n+l},s')$ by
\[
\Omega_t(\rho)=\left\{\,\rho'\in\cale(\bfG_{n+l},s')\mid
\left\langle\rho',R^{\bfG_{n+l}}_{\bfG_n\times\GL_l}(\rho\otimes\sigma)\right\rangle_{\!\bfG_{n+l}}\neq 0\,\right\}.
\]
Because we assume that $\sigma$ is cuspidal, one can see that
$C_{\bfG_n^*}(s)\times\GL_1^\dag$ is a Levi subgroup of $C_{\bfG_{n+l}^*}(s')$
where $\GL_1^\dag$ denotes the restriction to $\bfF_q$ of $\GL_1$ defined over a finite extension
(depending on $t$) of $\bfF_q$.
Suppose that $\rho'\in\Omega_t(\rho)$,
it is clear that $\rho'^c\in\Omega_t(\rho^c)$ for $\bfG_n=\Sp_{2n},\rmO^\epsilon_{2n}$
(\cf.~\cite{waldspurger} \S 4.4),
and $\rho'\cdot\sgn_{\bfG_{n+l}}\in\Omega_t(\rho\cdot\sgn_{\bfG_n})$ for $\bfG_n=\rmO^\epsilon_{2n}$.

We define
the relation $\Omega\colon\cale(C_{\bfG_n^*}(s),1)\rightarrow\cale(C_{\bfG_{n+l}^*}(s'),1)$
as in (\ref{0323}), i.e.,
\[
\Omega(\rho)=\Biggl\{\,\rho'\in\cale(C_{\bfG_{n+l}^*}(s'),1)\mid \left\langle\rho',
R^{C_{\bfG_{n+l}^*}(s')}_{C_{\bfG_n^*}(s)\times\GL_1^\dag}(\rho\otimes{\bf1})\right\rangle_{\!C_{\bfG_{n+l}^*}(s')}
\neq 0\,\Biggr\}
\]
for $\rho\in\cale(C_{\bfG_n^*}(s),1)$.
Then we say that a Lusztig correspondence
\[
\grL_s\colon\cale(\bfG,s)\rightarrow\cale(C_{\bfG^*}(s),1)
\]
is \emph{compatible with the parabolic induction} if the following diagram
\begin{equation}\label{0512}
\begin{CD}
\cale(\bfG_n,s) @> \Omega_t >> \cale(\bfG_{n+l},s') \\
@V \grL_s VV @VV \grL_{s'} V \\
\cale(C_{\bfG^*_n}(s),1)
@> \Omega >> \cale(C_{\bfG_{n+l}^*}(s'),1)
\end{CD}
\end{equation}
commutes for any $s$ and $s'=(s,t)$ given as above.
Now write
\[
C_{\bfG^*_n}(s)=\bfG^{(0)}(s)\times\bfG^{(1)}(s)\times\bfG^{(2)}(s),\qquad
C_{\bfG_{n+l}^*}(s')=\bfG^{(0)}(s')\times\bfG^{(1)}(s')\times\bfG^{(2)}(s')
\]
and $\grL_s(\rho)=\rho^{(0)}\otimes\rho^{(1)}\otimes\rho^{(2)}$ as in (\ref{0505}).
Then the diagram (\ref{0512}) can be described more precisely according to the following three cases:
\begin{enumerate}
\item If $t=1$ and so $l=1$, then $\bfG^{(0)}(s)=\bfG^{(0)}(s')$ and $\bfG^{(1)}(s)=\bfG^{(1)}(s')$,
and then the relation $\Omega$ is given by
\[
\Omega(\rho^{(0)}\otimes\rho^{(1)}\otimes\rho^{(2)})
=\rho^{(0)}\otimes\rho^{(1)}\otimes\Omega(\rho^{(2)})
\]
where $\Omega(\rho^{(2)})$ is defined as in (\ref{0323}).

\item If $t=-1$ and so $l=1$, then $\bfG^{(0)}(s)=\bfG^{(0)}(s')$ and $\bfG^{(2)}(s)=\bfG^{(2)}(s')$,
and then
\[
\Omega(\rho^{(0)}\otimes\rho^{(1)}\otimes\rho^{(2)})
=\rho^{(0)}\otimes\Omega(\rho^{(1)})\otimes\rho^{(2)}
\]
where $\Omega(\rho^{(1)})$ is defined as in (\ref{0323}).

\item If $t\neq \pm 1$, then $\bfG^{(1)}(s)=\bfG^{(1)}(s')$ and $\bfG^{(2)}(s)=\bfG^{(2)}(s')$,
and then
\[
\Omega(\rho^{(0)}\otimes\rho^{(1)}\otimes\rho^{(2)})
=\Omega(\rho^{(0)})\otimes\rho^{(1)}\otimes\rho^{(2)}
\]
where $\Omega(\rho^{(0)})$ is defined as in (\ref{0323}).
\end{enumerate}

\subsection{Modified Lusztig correspondence}
For a semisimple element $s\in G^*$, we define
\begin{equation}\label{0517}
\bfG^{(-)}=\bfG^{(1)},\qquad
\bfG^{(+)}=\begin{cases}
\bfG^{(2)}, & \text{if $\bfG=\SO_{2n+1}$ or $\rmO^\epsilon_{2n}$};\\
(\bfG^{(2)})^*, & \text{if $\bfG=\Sp_{2n}$}.
\end{cases}
\end{equation}
Combining $\call_1$ in Proposition~\ref{0301} (for $\bfG^{(0)}\times\bfG^{(-)}\times\bfG^{(+)}$)
and the inverse of $\grL_s$ in (\ref{0511}), we obtain
a bijection
\begin{align}\label{0509}
\begin{split}
\call_s\colon\cals_{\bfG^{(0)}(s)}\times\cals_{\bfG^{(-)}(s)}\times\cals_{\bfG^{(+)}(s)} &\rightarrow\cale(\bfG,s) \\
(x,\Lambda_1,\Lambda_2) &\mapsto  \rho_{x,\Lambda_1,\Lambda_2}=\rho^\bfG_{x,\Lambda_1,\Lambda_2}.
\end{split}
\end{align}
Note that from Proposition~\ref{0308}, Lemma~\ref{0321} and Lemma~\ref{0504}, we have
\begin{enumerate}
\item $(\rho^{\Sp_{2n}}_{x,\Lambda_1,\Lambda_2})^\sharp=(\rho^{\Sp_{2n}}_{x',\Lambda_1',\Lambda_2'})^\sharp$
if and only if
\begin{itemize}
\item $x'=x$,

\item $\Lambda_1'=\Lambda_1,\Lambda_1^\rmt$,

\item $\Lambda'_2=\Lambda_2$;
\end{itemize}

\item $(\rho^{\rmO^\epsilon_{2n}}_{x,\Lambda_1,\Lambda_2})^\sharp=(\rho^{\rmO^\epsilon_{2n}}_{x',\Lambda_1',\Lambda_2'})^\sharp$
if and only if
\begin{itemize}
\item $x'=x$,

\item $\Lambda_1'=\Lambda_1,\Lambda_1^\rmt$,

\item $\Lambda'_2=\Lambda_2,\Lambda_2^\rmt$.
\end{itemize}
\end{enumerate}
Moreover, diagram (\ref{0512}) becomes
\begin{equation}\label{0513}
\begin{CD}
\cals_{\bfG^{(0)}(s)}\times\cals_{\bfG^{(-)}(s)}\times\cals_{\bfG^{(+)}(s)}
@> \Omega>>
\cals_{\bfG^{(0)}(s')}\times\cals_{\bfG^{(-)}(s')}\times\cals_{\bfG^{(+)}(s')} \\
@V \call_{s} VV @VV \call_{s'} V \\
\cale(\bfG_n,s) @> \Omega_t >> \cale(\bfG_{n+l},s')
\end{CD}
\end{equation}
where the relation $\Omega$ is given as in Subsection~\ref{0326} or Subsection~\ref{0325}.

\begin{rem}
If $s$ is a semisimple element in $(G^*)^0$ such that $\bfG^{(0)}(s)$ is trivial,
then an irreducible character $\rho\in\cale(\bfG,s)$ is called \emph{quadratic unipotent}.
A Lusztig correspondence $\call_s\colon (\Lambda_1,\Lambda_2)\mapsto\rho_{\Lambda_1,\Lambda_2}$
of quadratic unipotent characters for $\bfG=\Sp_{2n}$, $\SO_{2n+1}$ or $\rmO^\epsilon_{2n}$ is
given in \cite{waldspurger} \S 4.11, \S 4.8, \S 4.4 respectively.
\end{rem}

\subsection{Theta correspondence and modified Lusztig correspondence}
First suppose that $(\bfG,\bfG')=(\Sp_{2n},\SO_{2n'+1})$, and $s\in G^*$, $s'\in G'^*$ semisimple.
Then $\bfG^{(-)}=\rmO^{\epsilon^{(-)}}_{2n^{(-)}}$,
$\bfG^{(+)}=\Sp_{2n^{(+)}}$,
$\bfG'^{(-)}=\Sp_{2n'^{(-)}}$,
and $\bfG'^{(+)}=\Sp_{2n'^{(+)}}$,
for some $\epsilon^{(-)}$ and some $n^{(-)}$, $n^{(+)}$, $n'^{(-)}$, $n'^{(+)}$.
The following proposition is from \cite{pan-Lusztig-correspondence} Proposition 8.3:

\begin{prop}\label{0518}
Let $(\bfG,\bfG')=(\Sp_{2n},\SO_{2n'+1})$, and $s\in G^*$, $s'\in G'^*$ semisimple.
Let
\begin{align*}
\call_s\colon\cals_{\bfG^{(0)}(s)}\times\cals_{\bfG^{(-)}(s)}\times\cals_{\bfG^{(+)}(s)}
&\rightarrow\cale(\bfG,s) \\
\call_{s'}\colon\cals_{\bfG'^{(0)}(s')}\times\cals_{\bfG'^{(-)}(s')}\times\cals_{\bfG'^{(+)}(s')}
&\rightarrow\cale(\bfG',s')
\end{align*}
be any modified Lusztig correspondences for $\bfG$ and $\bfG'$ respectively.
Then one of
\[
(\rho_{x,\Lambda_1,\Lambda_2},\rho_{x',\Lambda_1',\Lambda_2'}),\quad
(\rho_{x,\Lambda_1^\rmt,\Lambda_2},\rho_{x',\Lambda_1',\Lambda_2'})
\]
occurs in $\Theta_{\bfG,\bfG'}^\psi$ if and only if
\begin{itemize}
\item $s^{(0)}=-s'^{(0)}$ (up to conjugation) and $x=x'$,

\item $\Lambda_2=\Lambda'_1$, and

\item $(\Lambda_1,\Lambda'_2)$ or $(\Lambda_1^\rmt,\Lambda'_2)$ is in $\calb_{\bfG^{(-)}(s),\bfG'^{(+)}(s')}$.
\end{itemize}
\end{prop}

\begin{rem}
We shall see in Theorem~\ref{0507} that the modified Lusztig correspondence $\call_{s'}$
for $\SO_{2n+1}$ is unique.
\end{rem}

Next suppose that $(\bfG,\bfG')=(\Sp_{2n},\rmO_{2n'}^\epsilon)$ where $\epsilon=+$ or $-$,
and $s\in G^*$, $s'\in(G'^*)^0$ semisimple.
Then $\bfG^{(-)}=\rmO^{\epsilon^{(-)}}_{2n^{(-)}}$,
and $\bfG^{(+)}=\Sp_{2n^{(+)}}$,
$\bfG'^{(-)}=\rmO^{\epsilon'^{(-)}}_{2n'^{(-)}}$,
$\bfG'^{(+)}=\rmO^{\epsilon'^{(+)}}_{2n'^{(+)}}$,
for some $\epsilon^{(-)}$, $\epsilon'^{(-)}$, $\epsilon'^{(+)}$,
and some $n^{(-)}$, $n^{(+)}$, $n'^{(-)}$, $n'^{(+)}$.
The following proposition is from \cite{pan-Lusztig-correspondence} proposition~8.1:

\begin{prop}\label{0515}
Let $(\bfG,\bfG')=(\Sp_{2n},\rmO_{2n'}^\epsilon)$ where $\epsilon=+$ or $-$,
and $s\in G^*$, $s'\in(G'^*)^0$ semisimple.
Let
\begin{align*}
\call_s\colon\cals_{\bfG^{(0)}(s)}\times\cals_{\bfG^{(-)}(s)}\times\cals_{\bfG^{(+)}(s)}
&\rightarrow\cale(\bfG,s) \\
\call_{s'}\colon\cals_{\bfG'^{(0)}(s')}\times\cals_{\bfG'^{(-)}(s')}\times\cals_{\bfG'^{(+)}(s')}
&\rightarrow\cale(\bfG',s')
\end{align*}
be any modified Lusztig correspondences for $\bfG$ and $\bfG'$ respectively.
Then one of
\[
(\rho_{x,\Lambda_1,\Lambda_2},\rho_{x',\Lambda_1',\Lambda_2'}),\quad
(\rho_{x,\Lambda_1,\Lambda_2},\rho_{x',\Lambda_1'^\rmt,\Lambda_2'}),\quad
(\rho_{x,\Lambda_1,\Lambda_2},\rho_{x',\Lambda_1',\Lambda_2'^\rmt}),\quad
(\rho_{x,\Lambda_1,\Lambda_2},\rho_{x',\Lambda_1'^\rmt,\Lambda_2'^\rmt})
\]
occurs in $\Theta_{\bfG,\bfG'}^\psi$ if and only if
\begin{itemize}
\item $s^{(0)}=s'^{(0)}$ (up to conjugation) and $x=x'$,

\item $\Lambda_1=\Lambda'_1,\Lambda_1'^\rmt$, and

\item $(\Lambda_2,\Lambda'_2)$ or $(\Lambda_2,\Lambda_2'^\rmt)$ is in $\calb_{\bfG^{(+)}(s),\bfG'^{(+)}(s')}$.
\end{itemize}
\end{prop}

\section{Lusztig Correspondences for $\SO_{2n+1}$}

\subsection{Lusztig correspondence for $\SO_{2n+1}$}
Let $\bfG=\SO_{2n+1}$.
For a semisimple element $s\in G^*$, recall that (\cf.~(\ref{0511}), (\ref{0505})) we have
\begin{align*}
\grL_s\colon \cale(\bfG,s) &\rightarrow \cale(\bfG^{(0)}(s)\times\bfG^{(1)}(s)\times\bfG^{(2)}(s),1). \\
\rho &\mapsto \rho^{(0)}\otimes\rho^{(1)}\otimes\rho^{(2)}.
\end{align*}
Now we know that
$\bfG^{(1)}(s)=\Sp_{2n^{(1)}}$ and
$\bfG^{(2)}(s)=\Sp_{2n^{(2)}}$ for some non-negative integers $n^{(1)},n^{(2)}$ depending on $s$.

\begin{thm}\label{0507}
Let $\bfG=\SO_{2n+1}$ and $s\in G^*$.
There exists a unique bijection $\grL_s\colon\cale(\bfG,s)\rightarrow\cale(C_{\bfG^*}(s),1)$
satisfying (\ref{0201}), i.e., the Lusztig correspondence is unique.
\end{thm}
\begin{proof}
Let $\bfG=\SO_{2n+1}$ and $s\in G^*$,
and let $\grL_s,\grL'_s\colon\cale(\bfG,s)\rightarrow\cale(C_{\bfG^*}(s),1)$ be two Lusztig correspondences,
and write
\[
\grL_s(\rho)=\rho^{(0)}\otimes\rho^{(1)}\otimes\rho^{(2)},\qquad
\grL'_s(\rho)=\rho'^{(0)}\otimes\rho'^{(1)}\otimes\rho'^{(2)}
\]
as in (\ref{0505}).
Then by Lemma~\ref{0504}, we know that
$\rho^{(0)}=\rho'^{(0)}$,
$(\rho^{(1)})^\sharp=(\rho'^{(1)})^\sharp$,
$(\rho^{(2)})^\sharp=(\rho'^{(2)})^\sharp$.
Because now $\bfG=\SO_{2n+1}$,
we know that both $\bfG^{(1)}$ and $\bfG^{(2)}$ are symplectic groups.
Then we have $\rho^{(1)}=\rho'^{(1)}$ and $\rho^{(2)}=\rho'^{(2)}$ by Lemma~\ref{0302},
i.e., $\grL_s(\rho)=\grL'_s(\rho)$ for any $\rho\in\cale(\bfG,s)$.
\end{proof}

\begin{rem}
$\SO_{2n+1}$ is a connected group with connected center,
so the above result is considered in \cite{DM-Lusztig} theorem 7.1.
However, no extra condition is other than (\ref{0201}) in fact needed to make the 
Lusztig correspondence $\grL_s$ unique for this case.
\end{rem}

\begin{rem}
The Lusztig correspondence $\grL_s$ given in the theorem satisfies
the commutativity of the diagram in (\ref{0512}).
This is a special case of \cite{GM-guide} theorem~4.7.5.
\end{rem}

\begin{cor}\label{0603}
Let $\bfG=\SO_{2n+1}$ and $s\in G^*$.
Then there is a unique modified Lusztig correspondence
$\call_s\colon\cals_{\bfG^{(0)}(s)}\times\cals_{\bfG^{(-)}(s)}\times\cals_{\bfG^{(+)}(s)}\rightarrow\cale(\bfG,s)$.
\end{cor}

\begin{cor}\label{0607}
For $\rho,\rho'\in\cale(\SO_{2n+1})$,
then $\rho'^\sharp=\rho^\sharp$ if and only if $\rho'=\rho$.
\end{cor}
\begin{proof}
Let $\bfG=\SO_{2n+1}$, and let $\rho,\rho'\in\cale(\bfG)$.
It is obvious that $\rho=\rho'$ implies that $\rho^\sharp=\rho'^\sharp$.
Now we suppose that $\rho^\sharp=\rho'^\sharp$ and $\rho$ is in $\cale(\bfG,s)$ for some $s\in G^*$.
Because $\calv(\bfG)=\bigoplus_{(s)}\calv(\bfG,s)^\sharp$ and $\rho^\sharp=\rho'^\sharp$,
we see that $\rho'$ is also in $\cale(\bfG,s)$.
Write $\grL_s(\rho)=\rho^{(0)}\otimes\rho^{(1)}\otimes\rho^{(2)}$ and
$\grL_s(\rho')=\rho'^{(0)}\otimes\rho'^{(1)}\otimes\rho'^{(2)}$
where $\grL_s$ is given as in Theorem~\ref{0507}.
Then we know that
\[
(\rho^{(0)})^\sharp\otimes(\rho^{(1)})^\sharp\otimes(\rho^{(2)})^\sharp
=\grL_s(\rho^\sharp)
=\grL_s(\rho'^\sharp)
=(\rho'^{(0)})^\sharp\otimes(\rho'^{(1)})^\sharp\otimes(\rho'^{(2)})^\sharp,
\]
i.e., $(\rho^{(0)})^\sharp=(\rho'^{(0)})^\sharp$,
$(\rho^{(1)})^\sharp=(\rho'^{(1)})^\sharp$,
$(\rho^{(2)})^\sharp=(\rho'^{(2)})^\sharp$.
As in the proof of Theorem~\ref{0507}, we have $\rho^{(0)}=\rho'^{(0)}$,
$\rho^{(1)}=\rho'^{(1)}$, $\rho^{(2)}=\rho'^{(2)}$.
Because $\grL_s$ is a bijection, we conclude that $\rho=\rho'$.
\end{proof}

\begin{cor}\label{0611}
The bijection $\grL_1\colon\cale(\SO_{2n+1},1)\rightarrow\cale(\Sp_{2n},1)$ is given by
$\rho_\Lambda\mapsto\rho_{\Lambda^\rmt}$ for $\Lambda\in\cals_{\SO_{2n+1}}$,
i.e., the diagram
\[
\begin{CD}
\cals_{\SO_{2n+1}} @>>> \cals_{\Sp_{2n}} \\
@V \call_1 VV @VV \call_1 V \\
\cale(\SO_{2n+1},1) @> \grL_1 >> \cale(\Sp_{2n},1)
\end{CD}
\]
commutes where the mapping on the top is given by $\Lambda\mapsto\Lambda^\rmt$.
\end{cor}
\begin{proof}
Recall that $W_n=W_{\SO_{2n+1}}=W_{\Sp_{2n}}$ and
$\grL_1\colon R_{\bfT^*_w,1}^{\SO_{2n+1}}\mapsto R_{\bfT^*_w,1}^{\Sp_{2n}}$ for $w\in W_n$.
From (\ref{0330}), we see that the isometry
\[
\grL_1\colon\calv(\SO_{2n+1},1)\rightarrow\calv(\Sp_{2n},1)
\]
maps $R_\Sigma^{\SO_{2n+1}}\mapsto R_{\Sigma^\rmt}^{\Sp_{2n}}$
for $\Sigma\in\cals_{\SO_{2n+1}}^\sharp$.
Then we see that $\grL_1(\rho_\Lambda)^\sharp=(\rho_{\Lambda^\rmt})^\sharp$
for any $\Lambda\in\cals_{\SO_{2n+1}}$.
By Proposition~\ref{0507}, we conclude that $\grL_1(\rho_\Lambda)=\rho_{\Lambda^\rmt}$
and hence the lemma is proved.
\end{proof}

\begin{lem}\label{0708}
Let $\bfG=\SO_{2n+1}$,
and let $\call_s$ be the Lusztig correspondence given in Corollary~\ref{0603}.
Then
\begin{equation}\label{0606}
\rho_{x,\Lambda_1,\Lambda_2}\chi_\bfG=\rho_{x,\Lambda_2,\Lambda_1}
\end{equation}
where $\chi_\bfG$ denotes the spinor character of\/ $\bfG$.
\end{lem}
\begin{proof}
It is known that $\chi_\bfG R^\bfG_{\bfT^*,s}=R^\bfG_{\bfT^*,-s}$ for each pair $(\bfT^*,s)$,
and $\bfG^{(0)}(-s)=\bfG^{(0)}(s)$, $\bfG^{(-)}(-s)=\bfG^{(+)}(s)$ and $\bfG^{(+)}(-s)=\bfG^{(-)}(s)$.
Then the mapping $R^\bfG_{\bfT^*,s}\mapsto \chi_\bfG R^\bfG_{\bfT^*,s}$ induces an isometry
$\calv(\bfG,s)^\sharp\rightarrow\calv(\bfG,-s)^\sharp$ such that
$R^\bfG_{x,\Sigma_1,\Sigma_2}\mapsto R^\bfG_{x,\Sigma_2,\Sigma_1}$
where $(x,\Sigma_1,\Sigma_2)\in\cals_{\bfG^*(s)}^\sharp$ (\cf.~(\ref{0614})).
This means that
$(\rho_{x,\Lambda_1,\Lambda_2}\chi_\bfG)^\sharp=(\rho_{x,\Lambda_2,\Lambda_1})^\sharp$
for any $x\in\cals_{\bfG^{(0)}(s)}$, $\Lambda_1\in\cals_{\bfG^{(-)}(s)}$
and $\Lambda_2\in\cals_{\bfG^{(+)}(s)}$.
Then we conclude that $\rho_{x,\Lambda_1,\Lambda_2}\chi_\bfG=\rho_{x,\Lambda_2,\Lambda_1}$
by Corollary~\ref{0607}.
\end{proof}

\begin{rem}
For the case that $\bfG^{(0)}(s)$ is trivial, (\ref{0606}) is also given in \cite{waldspurger} proposition~4.8.
\end{rem}

\begin{exam}
Let $\bfG=\SO_3$.
Then the unique modified Lusztig correspondence
\begin{align*}
\call_s\colon\cals_{\bfG^{(0)}}\times\cals_{\bfG^{(-)}}\times\cals_{\bfG^{(+)}} &\rightarrow\cale(\bfG,s) \\
(x,\Lambda_1,\Lambda_2) &\mapsto\rho_{x,\Lambda_1,\Lambda_2}
\end{align*}
are given in the following table:
\[
\begin{tabular}{ccc|cccc}
\toprule
$G^{(0)}$ & $G^{(-)}$ & $G^{(+)}$ & $\cale(\bfG,s)$ & $\rho_{x,\Lambda_1,\Lambda_2}$ & number & cuspidality \\
\midrule
& & $\Sp_2(q)$ & ${\bf 1}_{\SO_3}$ & $\rho_{-,-,\binom{1}{-}}$ \\
& & & ${\rm St}_{\SO_3}$ & $\rho_{-,-,\binom{1,0}{1}}$ \\
\midrule
& $\Sp_2(q)$ & & $\chi_{\SO_3}$ & $\rho_{-,\binom{1}{-},-}$ \\
& & & ${\rm St}_{\SO_3}\chi_{\SO_3}$ & $\rho_{-,\binom{1,0}{1},-}$ \\
\midrule
$\GL_1(q)$ & & & & $\rho_{[1],-,-}$ & $\frac{q-3}{2}$ \\
\midrule
$\rmU_1(q)$ & & & & $\rho_{[\bar 1],-,-}$ & $\frac{q-1}{2}$ & $\checkmark$ \\
\bottomrule
\end{tabular}
\]
\end{exam}

\section{Lusztig Correspondence for $\rmO^\epsilon_{2n}$}

\subsection{Lusztig correspondence for $\rmO^\epsilon_{2n}$}
Let $\bfG=\rmO^\epsilon_{2n}$ where $\epsilon=+$ or $-$.
For a semisimple element $s\in(G^*)^0$, recall that (\cf.~(\ref{0511}), (\ref{0505})) we have
\begin{align*}
\grL_s\colon \cale(\bfG,s) &\rightarrow \cale(\bfG^{(0)}(s)\times\bfG^{(1)}(s)\times\bfG^{(2)}(s),1). \\
\rho &\mapsto \rho^{(0)}\otimes\rho^{(1)}\otimes\rho^{(2)}.
\end{align*}
Now we know that
$\bfG^{(1)}(s)=\rmO^{\epsilon^{(1)}}_{2n^{(1)}}$ and
$\bfG^{(2)}(s)=\rmO^{\epsilon^{(2)}}_{2n^{(2)}}$ for some non-negative integers $n^{(1)},n^{(2)}$
and some $\epsilon^{(1)},\epsilon^{(2)}$ depending on $s$ such that
$\epsilon^{(1)}\epsilon^{(2)}=\epsilon$.

\begin{lem}\label{0701}
Let $\bfG=\rmO_{2n}^\epsilon$ where $\epsilon=+$ or $-$, and let $s\in(G^*)^0$.
Let $\grL_s\colon\cale(\bfG,s)\rightarrow\cale(C_{\bfG^*}(s),1)$ be a Lusztig correspondence
and write $\grL_s(\rho)=\rho^{(0)}\otimes\rho^{(1)}\otimes\rho^{(2)}$.
Moreover, let $\grL'_s\colon\cale(\bfG,s)\rightarrow\cale(C_{\bfG^*}(s),1)$ be a bijective mapping
and write $\grL'_s(\rho)=\rho'^{(0)}\otimes\rho'^{(1)}\otimes\rho'^{(2)}$.
Then $\grL'_s$ is a Lusztig correspondence if and only if
\begin{itemize}
\item $\rho'^{(0)}=\rho^{(0)}$;

\item $\rho'^{(1)}=\rho^{(1)},\rho^{(1)}\cdot\sgn$;

\item $\rho'^{(2)}=\rho^{(2)},\rho^{(2)}\cdot\sgn$.
\end{itemize}
\end{lem}
\begin{proof}
First suppose that $\grL'_s$ is a Lusztig correspondence.
By Lemma~\ref{0504}, we know that $\grL_s(\rho)^\sharp=\grL'_s(\rho)^\sharp$, i.e.,
$\rho^{(0)}=\rho'^{(0)}$,
$(\rho^{(1)})^\sharp=(\rho'^{(1)})^\sharp$,
$(\rho^{(2)})^\sharp=(\rho'^{(2)})^\sharp$.
Now $\bfG^{(0)}$ is a product of general linear groups or unitary groups;
$\bfG^{(1)}$ and $\bfG^{(2)}$ are even orthogonal groups.
Then by Corollary~\ref{0314}, we have $\rho'^{(0)}=\rho^{(0)}$,
and $\rho'^{(i)}=\rho^{(i)}$ or $\rho'^{(i)}=\rho^{(i)}\cdot\sgn$ for $i=1,2$.

Next suppose that $\rho'^{(0)}=\rho^{(0)}$,
and $\rho'^{(i)}=\rho^{(i)}$ or $\rho'^{(i)}=\rho^{(i)}\cdot\sgn$ for $i=1,2$.
Then we have $(\rho^{(1)})^\sharp=(\rho'^{(1)})^\sharp$,
$(\rho^{(2)})^\sharp=(\rho'^{(2)})^\sharp$, i.e., $\grL_s(\rho)^\sharp=\grL'_s(\rho)^\sharp$ for any
$\rho\in\cale(\bfG,s)$.
This means that $\grL'_s$ also satisfies (\ref{0201}), i.e.,
$\grL'_s$ is also a Lusztig correspondence.
\end{proof}

Recall that for $\rho\in\cale(\bfG)$ the character $\rho^c$ is defined in Subsection~\ref{0520}.

\begin{lem}\label{0802}
Let $\bfG=\rmO^\epsilon_{2n}$, $s\in(G^*)^0$,
and let $\grL_s\colon\cale(\bfG,s)\rightarrow\cale(C_{\bfG^*}(s),1)$ be a Lusztig correspondence.
Suppose that $\rho\in\cale(\bfG,s)$ and write $\grL_s(\rho)=\rho^{(0)}\otimes\rho^{(1)}\otimes\rho^{(2)}$.
Then we have $\rho^c\in\cale(\bfG,s)$ and
\[
\grL_s(\rho^c)=\rho^{(0)}\otimes(\rho^{(1)}\cdot\sgn)\otimes\rho^{(2)}.
\]
\end{lem}
\begin{proof}
When $\bfG^{(0)}$ is trivial, i.e., when $\rho$ is quadratic unipotent,
the result is proved in \cite{waldspurger} \S 4.4.
The same argument still works here.
\end{proof}

\begin{cor}\label{0803}
Let $\bfG=\rmO^\epsilon_{2n}$, $s\in(G^*)^0$.
Suppose that $\rho\in\cale(\bfG,s)$ and $\rho=\rho_{x,\Lambda_1,\Lambda_2}$ under a modified Lusztig
correspondence $\call_s$.
Then $\rho^c=\rho_{x,\Lambda_1^\rmt,\Lambda_2}$.
\end{cor}
\begin{proof}
This follows from Lemma~\ref{0802} and Corollary~\ref{0406} immediately.
\end{proof}

\begin{cor}\label{0702}
For $\rho,\rho'\in\cale(\rmO^\epsilon_{2n})$,
then $\rho'^\sharp=\rho^\sharp$ if and only if
$\rho'=\rho,\rho^c,\rho\cdot\sgn,\rho^c\cdot\sgn$.
\end{cor}
\begin{proof}
Let $\bfG=\rmO^\epsilon_{2n}$, and let $\rho,\rho'\in\cale(\bfG)$.
By Lemma~\ref{0802}, Lemma~\ref{0504} and (\ref{0320}), we have $(\rho^c)^\sharp=\rho^\sharp$.
Moreover, because of (\ref{0303}), we have
$R^{\rmO^\epsilon_{2n}}_{\bfT^*,s}\cdot\sgn=R^{\rmO^\epsilon_{2n}}_{\bfT^*,s}$
and then $(\rho\cdot\sgn)^\sharp=\rho^\sharp$.
Therefore, if $\rho'=\rho$, $\rho^c$, $\rho\cdot\sgn$, or $\rho^c\cdot\sgn$,
we have $\rho'^\sharp=\rho^\sharp$.

Next we suppose that $\rho^\sharp=\rho'^\sharp$ and $\rho$ is in $\cale(\bfG,s)$ for some
$s\in(G^*)^0$.
As in proof of Corollary~\ref{0607}, we also have $\rho'\in\cale(\bfG,s)$.
Write $\grL_s(\rho)=\rho^{(0)}\otimes\rho^{(1)}\otimes\rho^{(2)}$ and
$\grL_s(\rho')=\rho'^{(0)}\otimes\rho'^{(1)}\otimes\rho'^{(2)}$.
By (\ref{0514}), we have
\[
(\rho^{(0)})^\sharp\otimes(\rho^{(1)})^\sharp\otimes(\rho^{(2)})^\sharp
=\grL_s(\rho^\sharp)
=\grL_s(\rho'^\sharp)
=(\rho'^{(0)})^\sharp\otimes(\rho'^{(1)})^\sharp\otimes(\rho'^{(2)})^\sharp,
\]
i.e., $\rho'^{(0)}=\rho^{(0)}$,
$\rho'^{(1)}=\rho^{(1)},\rho^{(1)}\cdot\sgn$, and $\rho'^{(2)}=\rho^{(2)},\rho^{(2)}\cdot\sgn$.
Now two sets $\{\rho,\rho^c,\rho\cdot\sgn,\rho^c\cdot\sgn\}$ and
\begin{multline*}
\left\{\rho^{(0)}\otimes\rho^{(1)}\otimes\rho^{(2)},
\rho^{(0)}\otimes(\rho^{(1)}\cdot\sgn)\otimes\rho^{(2)},\right. \\
\left.\rho^{(0)}\otimes\rho^{(1)}\otimes(\rho^{(2)}\cdot\sgn),
\rho^{(0)}\otimes(\rho^{(1)}\cdot\sgn)\otimes(\rho^{(2)}\cdot\sgn)\right\}
\end{multline*}
have the same cardinality.
We see that if $\rho'^\sharp=\rho^\sharp$,
then $\rho'$ must be one of $\rho$, $\rho^c$, $\rho\cdot\sgn$, or $\rho^c\cdot\sgn$.
\end{proof}

\subsection{Basic characters of $\rmO^\epsilon_{2n}$}
Let $\bfG=\rmO^\epsilon_{2n}$.
For a semisimple element $s\in(G^*)^0$,
let $\rho\in\cale(\bfG,s)$ and write
$\grL_s(\rho)=\rho^{(0)}\otimes\rho^{(1)}\otimes\rho^{(2)}$ for some Lusztig correspondence $\grL_s$.
Assume that both $\bfG^{(1)}(s),\bfG^{(2)}(s)$ are not trivial,
by Lemma~\ref{0701} and Corollary~\ref{0702} we know that any Lusztig correspondence gives a bijection between
$\{\rho,\rho^c,\rho\cdot\sgn,\rho^c\cdot\sgn\}$ and
\begin{multline*}
\left\{\rho^{(0)}\otimes\rho^{(1)}\otimes\rho^{(2)},
\rho^{(0)}\otimes(\rho^{(1)}\cdot\sgn)\otimes\rho^{(2)},\right. \\
\left.\rho^{(0)}\otimes\rho^{(1)}\otimes(\rho^{(2)}\cdot\sgn),
\rho^{(0)}\otimes(\rho^{(1)}\cdot\sgn)\otimes(\rho^{(2)}\cdot\sgn)\right\}.
\end{multline*}
Now we consider the situation where
\begin{itemize}
\item $\bfG^{(0)}(s)$ is a product of unitary groups, and $\rho^{(0)}$ is cuspidal; and

\item each of $\rho^{(1)},\rho^{(2)}$ is either cuspidal (i.e., is $\zeta_k^{\rm I},\zeta_k^{\rm II}$ for some $k$),
or is ${\bf 1}_{\rmO^+_2},\sgn_{\rmO^+_2}$.
\end{itemize}
An irreducible character $\rho$ satisfies the above conditions is call \emph{basic}.
Note that the class of basic characters is slightly larger than the class of cuspidal characters.

Now we denote the set $\{\rho,\rho^c,\rho\cdot\sgn,\rho^c\cdot\sgn\}$ by $\{\rho_1,\rho_2,\rho_3,\rho_4\}$.
From Proposition~\ref{0515} and the result in Subsection~\ref{0403},
we know that exactly two elements (says $\rho_1,\rho_2$) in $\{\rho_1,\rho_2,\rho_3,\rho_4\}$ first occur
in the correspondence for the pair $(\bfG,\bfG')=(\rmO^\epsilon_{2n},\Sp_{2(n-k')})$
where
\begin{equation}\label{0901}
k'=\begin{cases}
k, & \text{if $\rho^{(2)}=\zeta_k^{\rm I},\zeta_k^{\rm II}$};\\
1, & \text{if $\rho^{(2)}={\bf 1}_{\rmO^+_2},\sgn_{\rmO^+_2}$}.
\end{cases}
\end{equation}
Then we know that $\{\rho_3,\rho_4\}=\{\rho_1\cdot\sgn,\rho_2\cdot\sgn\}$ and $\rho_2=\rho_1^c$.

We know that $\rho_i\chi_\bfG\in\cale(\bfG,-s)$ for $=1,2,3,4$
where $\chi_\bfG$ denotes the spinor character (\cf.~Lemma~\ref{0708}),
and any Lusztig correspondence
\[
\grL_{-s}\colon\cale(\bfG,-s)\rightarrow
\cale(C_{\bfG^*}(-s),1)=\cale(\bfG^{(0)}(s)\times\bfG^{(2)}(s)\times\bfG^{(1)}(s),1)
\]
gives a bijection between $\{\rho_1\chi_\bfG,\rho_2\chi_\bfG,\rho_3\chi_\bfG,\rho_4\chi_\bfG\}$ and
\begin{multline*}
\left\{\rho^{(0)}\otimes\rho^{(2)}\otimes\rho^{(1)},
\rho^{(0)}\otimes(\rho^{(2)}\cdot\sgn)\otimes\rho^{(1)},\right. \\
\left.\rho^{(0)}\otimes\rho^{(2)}\otimes(\rho^{(1)}\cdot\sgn),
\rho^{(0)}\otimes(\rho^{(2)}\cdot\sgn)\otimes(\rho^{(1)}\cdot\sgn)\right\}.
\end{multline*}
Again, we know that there are exactly two elements in
$\{\rho_1\chi_\bfG,\rho_2\chi_\bfG,\rho_3\chi_\bfG,\rho_4\chi_\bfG\}$
first occur in the correspondence for the pair $(\rmO^\epsilon_{2n},\Sp_{2(n-h')})$
where
\begin{equation}\label{0902}
h'=\begin{cases}
h, & \text{if $\rho^{(1)}=\zeta_h^{\rm I},\zeta_h^{\rm II}$};\\
1, & \text{if $\rho^{(1)}={\bf 1}_{\rmO^+_2},\sgn_{\rmO^+_2}$}.
\end{cases}
\end{equation}

\begin{lem}\label{0703}
Keep the above settings.
There exists a unique character $\rho$ in $\{\rho_1,\rho_2,\rho_3,\rho_4\}$ above such that
\begin{itemize}
\item $\rho$ first occurs in the correspondence for $(\rmO^\epsilon_{2n},\Sp_{2(n-k')})$, and

\item $\rho\chi_\bfG$ first occurs in the correspondence for $(\rmO^\epsilon_{2n},\Sp_{2(n-h')})$
\end{itemize}
where $k',h'$ are given as in (\ref{0901}) and (\ref{0902}).
\end{lem}
\begin{proof}
We know that there exists $\rho',\rho''$ in $\{\rho_1,\rho_2,\rho_3,\rho_4\}$ such that
\begin{itemize}
\item $\rho',\rho'^c$ first occurs in the correspondence for $(\rmO^\epsilon_{2n},\Sp_{2(n-k')})$, and

\item $\rho''\chi_\bfG, (\rho''\chi_\bfG)^c$ first occurs in the correspondence for
$(\rmO^\epsilon_{2n},\Sp_{2(n-h')})$.
\end{itemize}
Moreover, we have
\begin{align*}
\{\rho_1,\rho_2,\rho_3,\rho_4\}
& =\{\rho',\rho'^c,\rho'\cdot\sgn,\rho'^c\cdot\sgn\}, \\
\{\rho_1\chi_\bfG,\rho_2\chi_\bfG,\rho_3\chi_\bfG,\rho_4\chi_\bfG\}
& =\{\rho''\chi_\bfG,(\rho''\chi_\bfG)^c,\rho''\chi_\bfG\cdot\sgn,(\rho''\chi_\bfG)^c\cdot\sgn\}.
\end{align*}
By \cite{waldspurger} (1) in \S 4.3, we know that
$(\rho''\chi_\bfG)^c=\rho''^c\chi_\bfG\cdot\sgn$,
and so the intersection
\[
\{\rho',\rho'^c\}\cap\{\rho'',\rho''^c\cdot\sgn\}
\]
clearly contains exactly one element.
\end{proof}

\begin{rem}\label{0709}
Let $\rho$ be the character given in Lemma~\ref{0703}.
Then
\begin{itemize}
\item $\rho^c$ first occurs in the correspondence for $(\rmO^\epsilon_{2n},\Sp_{2(n-k')})$,

\item $\rho^c\chi_\bfG$ first occurs in the correspondence for $(\rmO^\epsilon_{2n},\Sp_{2(n+h')})$,

\item $\rho^c\cdot\sgn$ first occurs in the correspondence for $(\rmO^\epsilon_{2n},\Sp_{2(n+k')})$,

\item $(\rho^c\cdot\sgn)\chi_\bfG$ first occurs in the correspondence for $(\rmO^\epsilon_{2n},\Sp_{2(n-h')})$,

\item $\rho\cdot\sgn$ first occurs in the correspondence for $(\rmO^\epsilon_{2n},\Sp_{2(n+k')})$,

\item $(\rho\cdot\sgn)\chi_\bfG$ first occurs in the correspondence for $(\rmO^\epsilon_{2n},\Sp_{2(n+h')})$
\end{itemize}
where $k',h'$ are given in (\ref{0901}) and (\ref{0902}) respectively.
\end{rem}

\begin{rem}\label{0705}
If $\bfG^{(1)}(s)$ or $\bfG^{(2)}(s)$ is trivial,
the situation is easier as follows.
\begin{enumerate}
\item If $\bfG^{(1)}(s)$ is trivial and $\bfG^{(2)}(s)$ is not,
then $\rho=\rho^c$, and $\rho,\rho\cdot\sgn$ are the only two irreducible characters whose uniform
projection is equal to $\rho^\sharp$.
And it is clear that exactly one of them first occurs in the correspondence for the pair
$(\rmO^\epsilon_{2n},\Sp_{2(n-k')})$.

\item If $\bfG^{(2)}(s)$ is trivial and $\bfG^{(1)}(s)$ is not,
then $\rho^c=\rho\cdot\sgn$,
and $\rho,\rho^c$ are the only two irreducible characters whose uniform projection is equal to $\rho^\sharp$.
Moreover, there is a unique element $\rho_1$ in $\{\rho,\rho^c\}$ such that
$\rho_1\chi_\bfG$ first occurs in the correspondence for the pair $(\rmO^\epsilon_{2n},\Sp_{2(n-h')})$.

\item If both $\bfG^{(1)}(s),\bfG^{(2)}(s)$ are trivial, then $\rho=\rho^c=\rho\cdot\sgn$,
and so $\rho$ is uniquely determined by its uniform projection.
\end{enumerate}
\end{rem}

\subsection{A uniqueness choice of $\call_s$}\label{0707}
To make a modified Lusztig correspondence
\[
\call_s\colon\cals_{\bfG^{(0)}(s)}\times\cals_{\bfG^{(-)}(s)}\times\cals_{\bfG^{(+)}(s)}
\rightarrow\cale(\bfG,s)
\]
uniquely determined for $\bfG=\rmO^\epsilon_{2n}$,
we follow the same idea used in Subsection~\ref{0521} and consider the following requirements:

\begin{enumerate}
\item[(I)] $\call_s$ is compatible with the parabolic induction as in (\ref{0513}).

\item[(II)] Suppose $\rho\in\cale(\bfG,s)$ is the unique basic character given in Lemma~\ref{0703}
(or in Remark~\ref{0705}),
then $\call_s$ is required that
\begin{equation}\label{0704}
\rho=\rho_{x,\Lambda_1,\Lambda_2},\quad
\rho^c=\rho_{x,\Lambda_1^\rmt,\Lambda_2},\quad
\rho^c\cdot\sgn=\rho_{x,\Lambda_1,\Lambda_2^\rmt},\quad
\rho\cdot\sgn=\rho_{x,\Lambda_1^\rmt,\Lambda_2^\rmt}
\end{equation}
where
\begin{itemize}
\item $x\in\cals_{\bfG^{(0)}(s)}$ is determined by $\rho^{(0)}$ (\cf.~Subsection~\ref{0326}),

\item $\Lambda_1\in\cals_{\bfG^{(-)}(s)}$ is given by
\[
\Lambda_1 =\begin{cases}
\Lambda_h^{\rm I}, & \text{if $\rho^{(1)}=\zeta_h^{\rm I},\zeta_h^{\rm II}$};\\
\binom{1}{0}, & \text{if $\rho^{(1)}={\bf1}_{\rmO^+_2},\sgn_{\rmO^+_2}$},
\end{cases}
\]

\item $\Lambda_2\in\cals_{\bfG^{(+)}(s)}$ is given by
\[
\Lambda_2 =\begin{cases}
\Lambda_k^{\rm I}, & \text{if $\rho^{(2)}=\zeta_k^{\rm I},\zeta_k^{\rm II}$};\\
\binom{1}{0}, & \text{if $\rho^{(2)}={\bf1}_{\rmO^+_2},\sgn_{\rmO^+_2}$}.
\end{cases}
\]
\end{itemize}
\end{enumerate}

\begin{exam}
Let $\bfG=\rmO^-_2$.
By (\ref{0704}), we have the following parametrization of irreducible characters
in various Lusztig series of $\rmO^-_2$:
\[
\begin{tabular}{ccc|cccc}
\toprule
$G^{(0)}$ & $G^{(-)}$ & $G^{(+)}$ & $\rho\in\cale(\bfG,s)$ & $\rho_{x,\Lambda_1,\Lambda_2}$ & number & cuspidality \\
\midrule
& & $\rmO^-_2(q)$ & ${\bf 1}_{\rmO^-_2}$ & $\rho_{-,-,\binom{-}{1,0}}$ & & $\checkmark$ \\
& & & $\sgn_{\rmO^-_2}$ & $\rho_{-,-,\binom{1,0}{-}}$ & & $\checkmark$ \\
\midrule
& $\rmO^-_2(q)$ & & $\chi_{\rmO^-_2}$ & $\rho_{-,\binom{-}{1,0},-}$ & & $\checkmark$ \\
& & & $\chi_{\rmO^-_2}\cdot\sgn_{\rmO^-_2}$ & $\rho_{-,\binom{1,0}{-},-}$ & & $\checkmark$ \\
\midrule
$\rmU_1(q)$ & & & $\chi^{(k)}$ & $\rho_{[\bar 1],-,-}$ & $\frac{q-1}{2}$ & $\checkmark$ \\
\bottomrule
\end{tabular}
\]
Here each $\chi^{(k)}$ is an irreducible character of degree $2$.
\end{exam}

\begin{thm}\label{0508}
Let $\bfG=\rmO^\epsilon_{2n}$ where $\epsilon=+$ or $-$, and let $s\in(G^*)^0$ be semisimple.
There exists a unique bijection $\grL_s\colon\cale(\bfG,s)\rightarrow\cale(C_{\bfG^*}(s),1)$
satisfying (\ref{0201}) and (I), (II) above.
\end{thm}
\begin{proof}
For $s\in(G^*)^0$,
let $\grL_s,\grL'_s\colon\cale(\bfG,s)\rightarrow\cale(C_{\bfG^*}(s),1)$ be two Lusztig
correspondences satisfying (I) and (II) above.
Let
\[
\call_s,\call'_s\colon\cals_{\bfG^{(0)}(s)}\times\cals_{\bfG^{(-)}(s)}\times\cals_{\bfG^{(+)}(s)}
\rightarrow\cale(\bfG,s)
\]
be the corresponding modified Lusztig correspondences.
Suppose that $\call_s(x,\Lambda_1,\Lambda_2)=\call'_s(x',\Lambda'_1,\Lambda'_2)$
for some $x,x'\in\cals_{\bfG^{(0)}(s)}$, $\Lambda_1,\Lambda_1'\in\cals_{\bfG^{(-)}(s)}$,
$\Lambda_2,\Lambda_2'\in\cals_{\bfG^{(+)}(s)}$.
Then by Lemma~\ref{0701} we know that $x=x'$, $\Lambda_1=\Lambda'_1,\Lambda_1'^\rmt$,
and $\Lambda_2=\Lambda'_2,\Lambda_2'^\rmt$.
So we need to show that $\Lambda_1=\Lambda_1'$ and $\Lambda_2=\Lambda_2'$.

For $\Lambda_1$, we consider the following two situations:
\begin{enumerate}
\item Suppose that $\Lambda_1$ is degenerate, i.e., $\Lambda_1=\Lambda_1^\rmt$.
This implies that $\Lambda_1=\Lambda'_1$ immediately.

\item Next we consider the case that $\Lambda_1$ is non-degenerate.
Suppose that $\rho$ is an irreducible constituent of the parabolic induced character
$R^{\rmO^\epsilon_{2n}}_{\rmO^\epsilon_{2n_0}\times\bfL}(\zeta\otimes\sigma)$
where $\zeta$ is a cuspidal character of $\rmO^\epsilon_{2n_0}(q)$ for some $n_0\leq n$,
$\bfL$ is a product of general linear groups, and $\sigma$ is a cuspidal character of $\bfL$.
Suppose that $\zeta\in\cale(\rmO^\epsilon_{2n_0},s_0)$ and write
\[
\zeta
=\call_{s_0}(x_0,\Lambda_{0,1},\Lambda_{0,2})
=\call'_{s_0}(x'_0,\Lambda'_{0,1},\Lambda'_{0,2})
\]
for some $x_0,x'_0\in\cals_{\bfG^{(0)}(s_0)}$, $\Lambda_{0,1},\Lambda_{0,1}'\in\cals_{\bfG^{(-)}(s_0)}$
and $\Lambda_{0,2},\Lambda_{0,2}'\in\cals_{\bfG^{(+)}(s_0)}$.
By condition (II) above, we know that $x_0=x'_0$, $\Lambda_{0,1}=\Lambda'_{0,1}$,
and $\Lambda_{0,2}=\Lambda'_{0,2}$.
\begin{enumerate}
\item Suppose that ${\rm def}(\Lambda_1)\neq 0$.
By condition (I), we have
\[
{\rm def}(\Lambda_1)={\rm def}(\Lambda_{0,1})={\rm def}(\Lambda'_{0,1})={\rm def}(\Lambda'_1)
\neq{\rm def}(\Lambda_1'^\rmt).
\]
This means that $\Lambda_1=\Lambda'_1$.

\item Suppose that $\Lambda_1$ is non-degenerate and ${\rm def}(\Lambda_1)=0$, i.e.,
$\Lambda_1,\Lambda_1'\in\cals_{\rmO^+_{2n^{(-)}}}$ for some $n^{(-)}$.
Note that for this case, $\Lambda_{0,1}=\Lambda_{0,1}'=\binom{-}{-}$.
Now we are going to prove this case by induction on $n^{(-)}$.
For $n^{(-)}=1$, the equality $\Lambda_1=\Lambda'_1$ is enforced by (II) above.
Next suppose that $n^{(-)}\geq 2$.
Because now $\Lambda_1^\rmt\neq\Lambda_1$, by Lemma~\ref{0328}, there exists
$\Lambda_{1,1}\in\cals_{\rmO^+_{2(n^{(-)}-1)}}$ such that $\Lambda_1\in\Omega(\Lambda_{1,1})$
and $\Lambda_1^\rmt\not\in\Omega(\Lambda_{1,1})$.
By the induction hypothesis and condition (I) above, we have
\[
\call_s(x,\Lambda_1,\Lambda_2)\in\Omega(\call_{s_1}(x,\Lambda_{1,1},\Lambda_2))
=\Omega(\call_{s_1}'(x,\Lambda_{1,1},\Lambda_2))
\not\ni \call'_s(x,\Lambda_1^\rmt,\Lambda_2).
\]
Now $\call'_s(x,\Lambda_1^\rmt,\Lambda_2)\neq\call_s(x,\Lambda_1,\Lambda_2)$ implies that
$\Lambda_1=\Lambda'_1$.
\end{enumerate}
\end{enumerate}

By the same argument we can also show that $\Lambda_2=\Lambda'_2$.
And then we conclude that $\call_s$ (and hence $\grL_s$) is uniquely determined
by (\ref{0201}) and (I),(II).
\end{proof}

\begin{cor}\label{0706}
Let $\bfG=\rmO^\epsilon_{2n}$,
and let $\call_s$ be the modified Lusztig correspondence given in Theorem~\ref{0508}.
Then $\rho_{x,\Lambda_1,\Lambda_2}\cdot\sgn=\rho_{x,\Lambda_1^\rmt,\Lambda_2^\rmt}$.
\end{cor}
\begin{proof}
From (\ref{0704}) we see that the assertion is true if $\rho$ is basic, i.e., if
\begin{itemize}
\item $\bfG^{(0)}$ is a product of unitary groups, and $x$ is cuspidal; and

\item each of $\Lambda_1,\Lambda_2$ is either cuspidal,
or is $\binom{1}{0},\binom{0}{1}$.
\end{itemize}
In general, suppose that $\rho'\in\Omega_t(\rho)$ for some $t$ corresponding a cuspidal character
of a general linear group (\cf.~Subsection~\ref{0615}).
Then we have $\rho'\cdot\sgn\in\Omega_t(\rho\cdot\sgn)$.
Then the corollary can be proved by induction on the rank of $\bfG$ via the similar argument in the proof of
Theorem~\ref{0508}.
\end{proof}

\begin{cor}
Let $\bfG=\rmO^\epsilon_{2n}$, and let $\call_s$ be the modified Lusztig correspondence given in
Theorem~\ref{0508}.
Then $\rho_{x,\Lambda_1,\Lambda_2}\chi_\bfG=\rho_{x,\Lambda_2,\Lambda_1}$
where $\chi_\bfG$ denotes the spinor character.
\end{cor}
\begin{proof}
Suppose that $\rho_{x,\Lambda_1,\Lambda_2}\in\cale(\bfG,s)$.
By the same argument in the proof of Lemma~\ref{0708},
we know that $\rho_{x,\Lambda_1,\Lambda_2}\chi_\bfG\in\cale(\bfG,-s)$ and
$(\rho_{x,\Lambda_1,\Lambda_2}\chi_\bfG)^\sharp=(\rho_{x,\Lambda_2,\Lambda_1})^\sharp$.
From Lemma~\ref{0703}, Remark~\ref{0709} and (\ref{0704}),
we see that the assertion is true if $\rho$ is basic.
Then the remaining proof is similar to that of Corollary~\ref{0706}.
\end{proof}

\section{Lusztig Correspondences for $\Sp_{2n}$}

\subsection{Lusztig correspondence for $\Sp_{2n}$}
Let $\bfG=\Sp_{2n}$.
For a semisimple element $s\in G^*=\SO_{2n+1}(q)$, recall that (\cf.~(\ref{0511}), (\ref{0505})) we have
\begin{align*}
\grL_s\colon \cale(\bfG,s) &\rightarrow \cale(\bfG^{(0)}(s)\times\bfG^{(1)}(s)\times\bfG^{(2)}(s),1). \\
\rho &\mapsto \rho^{(0)}\otimes\rho^{(1)}\otimes\rho^{(2)}.
\end{align*}
Now we know that $\bfG^{(1)}(s)=\rmO^{\epsilon^{(1)}}_{2n^{(1)}}$ and $\bfG^{(2)}(s)=\SO_{2n^{(2)}+1}$
for some $\epsilon^{(1)}=+$ or $-$, and some non-negative integers $n^{(1)},n^{(2)}$
depending on $s$.

\begin{lem}\label{0601}
Let $\bfG=\Sp_{2n}$ and $s\in G^*$.
Let $\grL_s\colon\cale(\bfG,s)\rightarrow\cale(C_{\bfG^*}(s),1)$ be a Lusztig correspondence
and write $\grL_s(\rho)=\rho^{(0)}\otimes\rho^{(1)}\otimes\rho^{(2)}$.
Moreover, let $\grL'_s\colon\cale(\bfG,s)\rightarrow\cale(C_{\bfG^*}(s),1)$ be a bijective mapping
and write $\grL'_s(\rho)=\rho'^{(0)}\otimes\rho'^{(1)}\otimes\rho'^{(2)}$.
Then $\grL'_s$ is a Lusztig correspondence if and only if
\begin{itemize}
\item $\rho'^{(0)}=\rho^{(0)}$;

\item $\rho'^{(1)}=\rho^{(1)},\rho^{(1)}\cdot\sgn$;

\item $\rho'^{(2)}=\rho^{(2)}$.
\end{itemize}
\end{lem}
\begin{proof}
The proof is similar to that of Lemma~\ref{0701}.
\end{proof}

\begin{lem}\label{0612}
Let $\bfG=\Sp_{2n}$, $s\in G^*$,
and let $\grL_s\colon\cale(\bfG,s)\rightarrow\cale(C_{\bfG^*}(s),1)$ be a Lusztig correspondence.
Suppose that $\rho\in\cale(\bfG,s)$ and write $\grL_s(\rho)=\rho^{(0)}\otimes\rho^{(1)}\otimes\rho^{(2)}$.
Then we have $\rho^c\in\cale(\bfG,s)$ and
\[
\grL_s(\rho^c)=\rho^{(0)}\otimes(\rho^{(1)}\cdot\sgn)\otimes\rho^{(2)}.
\]
\end{lem}
\begin{proof}
When $\bfG^{(0)}$ is trivial, i.e., when $\rho$ is quadratic unipotent,
the result is proved in \cite{waldspurger} \S 4.11.
The same argument still works here.
\end{proof}

\begin{cor}
Let $\bfG=\Sp_{2n}$, $s\in G^*$.
Suppose that $\rho\in\cale(\bfG,s)$ and $\rho=\rho_{x,\Lambda_1,\Lambda_2}$ under a modified Lusztig
correspondence $\call_s$.
Then $\rho^c=\rho_{x,\Lambda_1^\rmt,\Lambda_2}$.
\end{cor}
\begin{proof}
The proof is similar to that of Corollary~\ref{0803}.
\end{proof}

\begin{cor}\label{0608}
For $\rho,\rho'\in\cale(\Sp_{2n})$,
then $\rho'^\sharp=\rho^\sharp$ if and only if $\rho'=\rho,\rho^c$.
\end{cor}
\begin{proof}
Let $\bfG=\Sp_{2n}$, and let $\rho,\rho'\in\cale(\bfG)$.
By Lemma~\ref{0612}, $\rho'=\rho,\rho^c$ implies that $\rho'^\sharp=\rho^\sharp$.
Then remaining proof is similar to that of Corollary~\ref{0702}.
\end{proof}

\begin{exam}\label{0401}
Let $\bfG=\Sp_4$.
Now we follow the notations in \cite{srinivasan-sp4}.
\begin{enumerate}
\item Let $s_1\in\SO_5(q)$ such that $C_{\SO_5}(s_1)\simeq\rmO^+_4$.
Then we can check that
\[
\cale(\Sp_4,s_1)=\{\theta_1,\theta_2,\Phi_9,\theta_3,\theta_4\},\qquad
\cale(\rmO^+_4,1)=\bigl\{\rho_{\binom{2}{0}},\rho_{\binom{0}{2}},\rho_{\binom{1}{1}},
\rho_{\binom{2,1}{1,0}},\rho_{\binom{1,0}{2,1}}\bigr\}
\]
where $\rho_{\binom{2}{0}}={\bf1}_{\rmO^+_4}$, $\rho_{\binom{0}{2}}=\sgn_{\rmO^+_4}$.
Now $\grL_{s_1}\colon\cale(\Sp_4,s_1)\rightarrow\cale(\rmO^+_4,1)$ is a bijection such that
\begin{align*}
\{\theta_3,\theta_4\} &\rightarrow \bigl\{\rho_{\binom{2}{0}},\rho_{\binom{0}{2}}\bigr\}, \\
\{\Phi_9\} &\rightarrow \bigl\{\rho_{\binom{1}{1}}\bigr\}, \\
\{\theta_1,\theta_2\} &\rightarrow \bigl\{\rho_{\binom{2,1}{1,0}},\rho_{\binom{1,0}{2,1}}\bigr\}.
\end{align*}

\item Let $s_2\in\SO_5(q)$ such that $C_{\SO_5}(s_2)\simeq\rmO^-_4$.
Then we can check that
\[
\cale(\Sp_4,s_2)=\{\theta_5,\theta_6,\theta_7,\theta_8\},\qquad
\cale(\rmO^-_4,1)=\bigl\{\rho_{\binom{-}{2,0}},\rho_{\binom{2,0}{-}},\rho_{\binom{1}{2,1,0}},\rho_{\binom{2,1,0}{1}}\bigr\}
\]
where $\rho_{\binom{-}{2,0}}={\bf1}_{\rmO^-_4}$, $\rho_{\binom{2,0}{-}}=\sgn_{\rmO^-_4}$.
Now
$\grL_{s_2}\colon\cale(\Sp_4,s_2)\rightarrow\cale(\rmO^-_4,1)$ is a bijection such that
\begin{align*}
\{\theta_7,\theta_8\} &\rightarrow \bigl\{\rho_{\binom{-}{2,0}},\rho_{\binom{2,0}{-}}\bigr\}, \\
\{\theta_5,\theta_6\} &\rightarrow \bigl\{\rho_{\binom{1}{2,1,0}},\rho_{\binom{2,1,0}{1}}\bigr\}.
\end{align*}\end{enumerate}
\end{exam}

\subsection{Basic characters of $\Sp_{2n}$}\label{0904}
Let $\bfG=\Sp_{2n}$.
For a semisimple element $s$ in $G^*$,
let $\rho\in\cale(\bfG,s)$ and write $\grL_s(\rho)=\rho^{(0)}\otimes\rho^{(1)}\otimes\rho^{(2)}$
for a Lusztig correspondence $\grL_s$.
Assume that $\bfG^{(1)}(s)$ is not trivial,
we know that any Lusztig correspondence gives a bijection between $\{\rho,\rho^c\}$ and
\[
\left\{\rho^{(0)}\otimes\rho^{(1)}\otimes\rho^{(2)},
\rho^{(0)}\otimes(\rho^{(1)}\cdot\sgn)\otimes\rho^{(2)}\right\}.
\]
Now we consider the situation where
\begin{itemize}
\item $\bfG^{(0)}$ is a product of unitary groups,
and $\rho^{(0)}$ is cuspidal; and

\item $\rho^{(1)}$ is either cuspidal (i.e., is $\zeta_k^{\rm I},\zeta_k^{\rm II}$ for some $k$),
or is ${\bf 1}_{\rmO^+_2},\sgn_{\rmO^+_2}$; and

\item $\rho^{(2)}$ is cuspidal.
\end{itemize}
An irreducible character $\rho$ satisfies the above conditions is call \emph{basic}.
From Proposition~\ref{0515} and the result in Subsection~\ref{0403},
we know that exactly one element in $\{\rho,\rho^c\}$ first occurs in the correspondence for the pair
$(\Sp_{2n},\SO_{2(n-k')+1})$
where
\[
k'=\begin{cases}
k, & \text{if $\rho^{(1)}=\zeta_k^{\rm I},\zeta_k^{\rm II}$};\\
1, & \text{if $\rho^{(1)}={\bf 1}_{\rmO^+_2},\sgn_{\rmO^+_2}$}.
\end{cases}
\]

\subsection{The choice of $\call_s$ with respect to $(\Sp_{2n},\SO_{2n'+1})$}\label{0613}

To make a modified Lusztig correspondence
\[
\call_s\colon\cals_{\bfG^{(0)}(s)}\times\cals_{\bfG^{(-)}(s)}\times\cals_{\bfG^{(+)}(s)}
\rightarrow\cale(\bfG,s)
\]
for $\Sp_{2n}$ with respect to the dual pair $(\Sp_{2n},\SO_{2n'+1})$
we consider the following requirements:
\begin{enumerate}
\item[(I)] We require that $\call_s$ is compatible with the parabolic induction as in (\ref{0513}).

\item[(II)] Suppose $\rho\in\cale(\bfG,s)$ is the unique irreducible basic character which first occurs
in the correspondence for $(\Sp_{2n},\SO_{2(n-k')+1})$.
Then $\call_s$ is required that
\begin{equation}\label{0609}
\rho=\rho_{x,\Lambda_1,\Lambda_2},\quad
\rho^c=\rho_{x,\Lambda_1^\rmt,\Lambda_2}
\end{equation}
where
\begin{itemize}
\item $x\in\cals_{\bfG^{(0)}(s)}$ is uniquely determined by $\rho^{(0)}$
(\cf.~Subsection~\ref{0326}),

\item $\Lambda_1\in\cals_{\bfG^{(-)}(s)}$ given by
\[
\Lambda_1 =\begin{cases}
\Lambda_k^{\rm I}, & \text{if $\rho^{(1)}=\zeta_k^{\rm I},\zeta_k^{\rm II}$};\\
\binom{1}{0}, & \text{if $\rho^{(1)}={\bf1}_{\rmO^+_2},\sgn_{\rmO^+_2}$},
\end{cases}
\]

\item $\Lambda_2\in\cals_{\bfG^{(+)}(s)}$ is uniquely determined by $\rho^{(2)}$,
i.e., $\rho^{(2)}=\rho_{\Lambda_2^\rmt}$
(\cf.~(\ref{0517}) and Corollary~\ref{0611}).
\end{itemize}
\end{enumerate}

\begin{thm}\label{0506}
Let $\bfG=\Sp_{2n}$ and $s\in G^*$ semisimple.
There exists a unique bijection
$\grL_s\colon\cale(\bfG,s)\rightarrow\cale(C_{\bfG^*}(s),1)$
satisfying (\ref{0201}) and (I), (II) above.
\end{thm}
\begin{proof}
The proof is similar to that of Proposition~\ref{0305}.
For $s\in G^*$,
let $\grL_s,\grL'_s\colon\cale(\bfG,s)\rightarrow\cale(C_{\bfG^*}(s),1)$ be two Lusztig
correspondence satisfying (I) and (II) above.
Let
\[
\call_s,\call'_s\colon\cals_{\bfG^{(0)}(s)}\times\cals_{\bfG^{(-)}(s)}\times\cals_{\bfG^{(+)}(s)}
\rightarrow\cale(\bfG,s)
\]
be the corresponding modified Lusztig correspondences.
Let $\rho\in\cale(\bfG,s)$ and suppose that
\[
\rho=\call_s(x,\Lambda_1,\Lambda_2)=\call'_s(x',\Lambda'_1,\Lambda'_2)
\]
for some $x,x'\in\cals_{\bfG^{(0)}(s)}$, $\Lambda_1,\Lambda_1'\in\cals_{\bfG^{(-)}(s)}$,
$\Lambda_2,\Lambda_2'\in\cals_{\bfG^{(+)}(s)}$.
Then by Lemma~\ref{0601} we know that $x=x'$, $\Lambda_2=\Lambda'_2$, and $\Lambda_1=\Lambda'_1,\Lambda_1'^\rmt$.
So our goal is to prove that $\Lambda_1=\Lambda_1'$.
Now we consider the following situations:
\begin{enumerate}
\item Suppose that $\Lambda_1$ is degenerate, i.e., $\Lambda_1=\Lambda_1^\rmt$.
This of course implies that $\Lambda_1=\Lambda'_1$.

\item Next suppose that $\Lambda_1$ is non-degenerate.
Suppose that $\rho$ is an irreducible constituent of the parabolic induced character
$R^{\Sp_{2n}}_{\Sp_{2n_0}\times\bfL}(\zeta\otimes\sigma)$
where $\zeta$ is a cuspidal character of $\Sp_{2n_0}(q)$ for some $n_0\leq n$,
$\bfL$ is a product of general linear groups, and $\sigma$ is a cuspidal character of $\bfL$.
Suppose that $\zeta\in\cale(\Sp_{2n_0},s_0)$ for some $s_0$ and write
\[
\zeta
=\call_{s_0}(x_0,\Lambda_{0,1},\Lambda_{0,2})
=\call'_{s_0}(x'_0,\Lambda'_{0,1},\Lambda'_{0,2})
\]
for some $x_0,x'_0\in\cals_{\bfG^{(0)}(s_0)}$, $\Lambda_{0,1},\Lambda_{0,1}'\in\cals_{\bfG^{(-)}(s_0)}$
and $\Lambda_{0,2},\Lambda_{0,2}'\in\cals_{\bfG^{(+)}(s_0)}$.
By (II) above, we know that $\Lambda_{0,1}=\Lambda_{0,1}'$.
By the same argument as in the proof of Theorem~\ref{0508},
we conclude that $\Lambda_1=\Lambda'_1$.
\end{enumerate}
Therefore the theorem is proved.
\end{proof}

\begin{rem}
Note that the modified Lusztig correspondence $\call_s$ in the theorem depends on
the theta correspondence $\Theta_{\bfG,\bfG'}^\psi$,
in particular, it depends on the choice of $\psi$.
Let $\call'_s$ be the corresponding modified Lusztig correspondence with respect to
another character $\psi'=\psi_a$ where $a\in\bfF_q^\times$ is a non-square element.
Then by Lemma~\ref{0405} and Lemma~\ref{0612} we see that
\[
\call_s(x,\Lambda_1,\Lambda_2)=\call'_s(x,\Lambda_1^\rmt,\Lambda_2)
\]
for $x\in\cals_{\bfG^{(0)}(s)}$,
$\Lambda_1\in\cals_{\bfG^{(-)}(s)}$,
$\Lambda_2\in\cals_{\bfG^{(+)}(s)}$.
\end{rem}

To justify the choice of $\call_s$ in the above theorem,
we have the following result which refines Proposition~\ref{0518} .

\begin{thm}\label{0610}
Let $(\bfG,\bfG')=(\Sp_{2n},\SO_{2n'+1})$, and $s\in G^*$, $s'\in G'^*$ semisimple.
Let
\begin{align*}
\call_s\colon\cals_{\bfG^{(0)}(s)}\times\cals_{\bfG^{(-)}(s)}\times\cals_{\bfG^{(+)}(s)}
& \rightarrow\cale(\bfG,s) \\
\call_{s'}\colon\cals_{\bfG'^{(0)}(s')}\times\cals_{\bfG'^{(-)}(s')}\times\cals_{\bfG'^{(+)}(s')}
& \rightarrow\cale(\bfG',s')
\end{align*}
be the modified Lusztig correspondence for $\bfG$ given in Theorem~\ref{0506},
and the Lusztig correspondence for $\bfG'$ given by Theorem~\ref{0507} respectively.
Then $(\rho_{x,\Lambda_1,\Lambda_2},\rho_{x',\Lambda_1',\Lambda_2'})\in\Theta_{\bfG,\bfG'}^\psi$
if and only if
\begin{itemize}
\item $s^{(0)}=-s'^{(0)}$ (up to conjugation) and $x=x'$,

\item $\Lambda_2=\Lambda'_1$, and

\item $(\Lambda_1,\Lambda'_2)\in\calb_{\bfG^{(-)}(s),\bfG^{(+)}(s')}$.
\end{itemize}
\end{thm}
\begin{proof}
Suppose that $(\rho_{x,\Lambda_1,\Lambda_2},\rho_{x',\Lambda_1',\Lambda_2'})\in\Theta_{\bfG,\bfG'}^\psi$.
Then by Proposition~\ref{0518},
we have
\begin{itemize}
\item $s^{(0)}=-s'^{(0)}$ and $x=x'$,

\item $\Lambda_2=\Lambda'_1$, and

\item $(\Lambda_1,\Lambda'_2)$ or $(\Lambda_1^\rmt,\Lambda'_2)$ is in $\calb_{\bfG^{(1)}(s),\bfG'^{(2)}(s')}$.
\end{itemize}
Now we want to show that in fact $(\Lambda_1,\Lambda'_2)\in\calb_{\bfG^{(1)}(s),\bfG'^{(2)}(s')}$.
Note that $(\bfG^{(1)}(s),\bfG'^{(2)}(s'))=(\rmO^{\epsilon^{(1)}}_{2n^{(1)}},\Sp_{2n'^{(2)}})$
for some $\epsilon^{(1)}$, and some $n^{(1)},n'^{(2)}$.
Now we consider the following situations:
\begin{enumerate}
\item Suppose that ${\rm def}(\Lambda_1)\neq 0$.
First we consider the case that both $\Lambda_1,\Lambda_2'$ are cuspidal.
\begin{enumerate}
\item Suppose that $(\rho_{x,\Lambda_1,\Lambda_2},\rho_{x',\Lambda_1',\Lambda_2'})$
first occurs in the correspondence for the pair $(\Sp_{2n},\SO_{2(n-k)+1})$ for some $k$,
i.e., $(\bfG^{(1)}(s),\bfG'^{(2)}(s'))=(\rmO^{\epsilon_k}_{2k^2},\Sp_{2k(k-1)})$.
From our choice of $\call_s$ (\cf.~(\ref{0609})) and $\call_{s'}$,
we know that $\Lambda_1=\Lambda_k^\rmI$ (\cf.~(\ref{0402})) and $\Lambda_2'=\Lambda_{k-1}^{\Sp}$ (\cf.~(\ref{0313})),
and it is clearly that $(\Lambda_1,\Lambda'_2)\in\calb_{\bfG^{(1)}(s),\bfG'^{(2)}(s')}$.

\item Suppose that $(\rho_{x,\Lambda_1,\Lambda_2},\rho_{x',\Lambda_1',\Lambda_2'})$
first occurs in the correspondence for the pair $(\Sp_{2n},\SO_{2(n+k)+1})$,
i.e., $(\bfG^{(1)}(s),\bfG'^{(2)}(s'))=(\rmO^{\epsilon_k}_{2k^2},\Sp_{2k(k+1)})$.
Now we have $\Lambda_1=\Lambda_k^{\rm II}$ and $\Lambda_2'=\Lambda_k^{\Sp}$,
and again $(\Lambda_1,\Lambda'_2)\in\calb_{\bfG^{(1)}(s),\bfG'^{(2)}(s')}$.
\end{enumerate}
Now if $\Lambda_1$ or $\Lambda_2'$ is not cuspidal,
by the same argument in the proof of \cite{pan-finite-unipotent} proposition 6.4 we
still conclude that $(\Lambda_1,\Lambda'_2)\in\calb_{\bfG^{(1)}(s),\bfG'^{(2)}(s')}$.

\item Suppose that ${\rm def}(\Lambda_1)=0$.
Then ${\rm def}(\Lambda_2')=1$.
\begin{enumerate}
\item Suppose that $\Lambda_1=\binom{-}{-}$, i.e.,
$(\bfG^{(1)}(s),\bfG'^{(2)}(s'))=(\rmO^+_0,\Sp_{2n'^{(2)}})$.
This case is obvious.

\item Suppose that $\Lambda_1=\binom{1}{0}$ or $\binom{0}{1}$.
If $(\rho_{x,\Lambda_1,\Lambda_2},\rho_{x',\Lambda_1',\Lambda_2'})$
first occurs in the correspondence for the pair $(\Sp_{2n},\SO_{2n-1})$,
then $\Lambda_1=\binom{1}{0}$, $(\bfG^{(1)}(s),\bfG'^{(2)}(s'))=(\rmO^+_2,\Sp_0)$,
and $\Lambda_2'=\binom{0}{-}$.
If $(\rho_{x,\Lambda_1,\Lambda_2},\rho_{x',\Lambda_1',\Lambda_2'})$
first occurs in the correspondence for $(\Sp_{2n},\SO_{2n+1})$,
then $\Lambda_1=\binom{0}{1}$, $(\bfG^{(1)}(s),\bfG'^{(2)}(s'))=(\rmO^+_2,\Sp_2)$,
and $\Lambda_2'=\binom{1}{-}$.
We have $(\Lambda_1,\Lambda'_2)\in\calb_{\bfG^{(1)}(s),\bfG'^{(2)}(s')}$ for
both situations.
\end{enumerate}
Now for general $\Lambda_1$ and $\Lambda_2'$,
we can use the same argument in \cite{pan-finite-unipotent} \S 6
(in particular the proof of proposition 6.20) and conclude that
$(\Lambda_1,\Lambda'_2)\in\calb_{\bfG^{(1)}(s),\bfG'^{(2)}(s')}$.
\end{enumerate}
Therefore, the theorem is proved.
\end{proof}

\subsection{The choice of $\call_s$ with respect to $(\Sp_{2n},\rmO^\epsilon_{2n'})$}
Keep the settings in Subsection~\ref{0904}.
Now any Lusztig correspondence gives a bijection
\[
\{\rho,\rho^c\}\rightarrow
\left\{\rho^{(0)}\otimes\rho^{(1)}\otimes\rho^{(2)},
\rho^{(0)}\otimes(\rho^{(1)}\cdot\sgn)\otimes\rho^{(2)}\right\}
\]
where $\rho\in\cale(\bfG,s)$ is a basic character.
Now $\rho^{(2)}$ is a cuspidal character of $\SO_{2n^{(2)}+1}$,
so we assume that $\rho^{(2)}=\zeta^{\SO_{\rm odd}}_k$ for some $k$.
By Proposition~\ref{0515}, we know that both $\rho,\rho^c$
first occur in the correspondence for the pair $(\bfG,\bfG')=(\Sp_{2n},\rmO^\epsilon_{2(n-k)})$
for some $\epsilon$.
Suppose that $(\rho,\rho')\in\Theta^\psi_{\bfG,\bfG'}$ for some unique $\rho'\in\cale(\bfG',s')$
and write $\rho'=\rho_{x',\Lambda_1',\Lambda'_2}$ where $\call_{s'}$ is given in Theorem~\ref{0508}.
Note that $(\rho^c,\rho'^c)$ also occurs in $\Theta^\psi_{\bfG,\bfG'}$ by Lemma~\ref{0522}.

Now besides the condition (I) given in Subsection~\ref{0613}, we also consider the following requirement:
\begin{enumerate}
\item[(III)] Suppose $\rho\in\cale(\bfG,s)$ is a basic character given above.
Then $\call_s$ is required that
\begin{equation}\label{0905}
\rho=\rho_{x,\Lambda_1,\Lambda_2},\quad
\rho^c=\rho_{x,\Lambda_1^\rmt,\Lambda_2}
\end{equation}
where
\begin{itemize}
\item $x\in\cals_{\bfG^{(0)}(s)}$ is uniquely determined by $\rho^{(0)}$,

\item $\Lambda_1\in\cals_{\bfG^{(-)}(s)}$ given by $\Lambda_1=\Lambda'_1$ where $\Lambda'_1$ is determined by
$\rho'$ as above,

\item $\Lambda_2\in\cals_{\bfG^{(+)}(s)}$ is uniquely determined by $\rho^{(2)}$, i.e.,
$\rho^{(2)}=\rho_{\Lambda_2^\rmt}$ (\cf.~(\ref{0517}) and Corollary~\ref{0611}).
\end{itemize}
\end{enumerate}

\begin{thm}\label{0903}
Let $\bfG=\Sp_{2n}$ and $s\in G^*$ semisimple.
There exists a unique bijection $\grL_s\colon\cale(\bfG,s)\rightarrow\cale(C_{\bfG^*}(s),1)$
satisfying (\ref{0201}), (I) in Subsection~\ref{0613}, and (III) above.
\end{thm}
\begin{proof}
The proof is analogous to that of Theorem~\ref{0506}.
\end{proof}

Now we have the following result which refines Proposition~\ref{0515} .

\begin{thm}\label{0801}
Let $(\bfG,\bfG')=(\Sp_{2n},\rmO_{2n'}^\epsilon)$ where $\epsilon=+$ or $-$,
and $s\in G^*$, $s'\in(G'^*)^0$ semisimple.
Let
\begin{align*}
\call_s\colon\cals_{\bfG^{(0)}(s)}\times\cals_{\bfG^{(-)}(s)}\times\cals_{\bfG^{(+)}(s)}
& \rightarrow\cale(\bfG,s) \\
\call_{s'}\colon\cals_{\bfG^{(0)}(s')}\times\cals_{\bfG^{(-)}(s')}\times\cals_{\bfG^{(+)}(s')}
& \rightarrow\cale(\bfG',s')
\end{align*}
be the modified Lusztig correspondence for $\bfG$ given in Theorem~\ref{0903},
and the modified Lusztig correspondence for $\bfG'$ given by Theorem~\ref{0508} respectively.
Then $(\rho_{x,\Lambda_1,\Lambda_2},\rho_{x',\Lambda_1',\Lambda_2'})\in \Theta_{\bfG,\bfG'}^\psi$
if and only if
\begin{itemize}
\item $s^{(0)}=s'^{(0)}$ (up to conjugation) and $x=x'$,

\item $\Lambda_1=\Lambda'_1$, and

\item $(\Lambda_2,\Lambda'_2)\in\calb_{\bfG^{(+)}(s),\bfG^{(+)}(s')}$.
\end{itemize}
\end{thm}
\begin{proof}
Suppose that $(\rho_{x,\Lambda_1,\Lambda_2},\rho_{x',\Lambda_1',\Lambda_2'})\in \Theta_{\bfG,\bfG'}^\psi$.
Then by Proposition~\ref{0515}, we have
\begin{itemize}
\item $s^{(0)}=s'^{(0)}$ and $x=x'$,

\item $\Lambda_1=\Lambda'_1,\Lambda_1'^\rmt$, and

\item $(\Lambda_2,\Lambda'_2)$ or $(\Lambda_2^\rmt,\Lambda'_2)$ is in $\calb_{\bfG^{(+)}(s),\bfG'^{(+)}(s')}$.
\end{itemize}
Note that $(\bfG^{(+)}(s),\bfG'^{(+)}(s'))=(\Sp_{2n^{(+)}},\rmO^{\epsilon'^{(+)}}_{2n'^{(+)}})$
for some $n^{(+)},\epsilon'^{(+)},n'^{(+)}$.
By the same argument in the proof of Theorem~\ref{0610}, we can conclude that
$(\Lambda_2,\Lambda'_2)\in\calb_{\bfG^{(+)}(s),\bfG'^{(+)}(s')}$.

Next we want to show that $\Lambda_1=\Lambda'_1$.
If $\Lambda_1$ is cuspidal or is equal to $\binom{1}{0},\binom{0}{1}$, and $\Lambda_2$ is cuspidal,
we have $\Lambda_1=\Lambda_1'$ by the requirement of $\call_s$ in (\ref{0905}).
Then by the same argument in the proof of Proposition~\ref{0305}, we can conclude
that $\Lambda_1=\Lambda_1'$ for general $\Lambda_1$.
\end{proof}

\bibliography{refer}
\bibliographystyle{amsalpha}

\end{document}